\documentclass[10pt]{amsart}
\usepackage[cp1251]{inputenc}
\usepackage[english]{babel}
\usepackage{amsmath}
\usepackage{amssymb}
\usepackage{amsfonts}

\usepackage[linktocpage=true, colorlinks=true, linkcolor=blue, citecolor=blue, urlcolor=blue]{hyperref}

\begin{document}
\renewcommand{\refname}{References}

\thispagestyle{empty}

\title[Expansion of Iterated Stratonovich Stochastic Integrals 
of Multiplicity 3]
{Expansion of Iterated Stratonovich Stochastic Integrals of Multiplicity 3
Based on Generalized Multiple Fourier Series Converging in the Mean: 
General Case of Series Summation}
\author[D.F. Kuznetsov]{Dmitriy F. Kuznetsov}
\address{Dmitriy Feliksovich Kuznetsov
\newline\hphantom{iii} Peter the Great Saint-Petersburg Polytechnic University,
\newline\hphantom{iii} Polytechnicheskaya ul., 29,
\newline\hphantom{iii} 195251, Saint-Petersburg, Russia}%
\email{sde\_kuznetsov@inbox.ru}

\maketitle{\small
\begin{quote}
\noindent{\sc Abstract.}
The article is devoted to the development of the method of expansion
and mean-square approximation of iterated Ito stochastic integrals
based on generalized multiple Fourier series converging in the mean.
We adapt this method for iterated Stratonovich stochastic
integrals of multiplicity 3 
from the Taylor--Stratonovich expansion.
The main result of the article has been derived 
using the triple Fourier--Legendre series and
triple trigonometric Fourier series for the 
general case of series summation.
Some recent results on the expansion of iterated Stratonovich 
stochastic integrals of multiplicities 3 to 6 are given.
The results of the article can be applied to
the numerical integration of Ito stochastic differential equations
in accordance with the strong criterion of convergence.
\medskip
\end{quote}
}


\setlength{\baselineskip}{1.8em}

\tableofcontents

\setlength{\baselineskip}{1.2em}


\section{Introduction}

\vspace{5mm}

Let $(\Omega,$ ${\rm F},$ ${\sf P})$ be a complete probability space, let 
$\{{\rm F}_t, t\in[0,T]\}$ be a nondecreasing right-continous 
family of $\sigma$-algebras of ${\rm F},$
and let ${\bf f}_t$ be a standard $m$-dimensional Wiener stochastic 
process, which is
${\rm F}_t$-measurable for any $t\in[0, T].$ We assume that the components
${\bf f}_{t}^{(i)}$ $(i=1,\ldots,m)$ of this process 
are independent. Consider
an Ito stochastic differential equation (SDE) 
in the integral form

\vspace{-1mm}
\begin{equation}
\label{1.5.2}
{\bf x}_t={\bf x}_0+\int\limits_0^t {\bf a}({\bf x}_{\tau},\tau)d\tau+
\int\limits_0^t B({\bf x}_{\tau},\tau)d{\bf f}_{\tau},\ \ \
{\bf x}_0={\bf x}(0,\omega).
\end{equation}

\vspace{2mm}
\noindent
Here ${\bf x}_t$ is some $n$-dimensional stochastic process 
satisfying the equation (\ref{1.5.2}). 
The non-random functions ${\bf a}: \mathbb{R}^n\times[0, T]\to\mathbb{R}^n$,
$B: \mathbb{R}^n\times[0, T]\to\mathbb{R}^{n\times m}$
guarantee the existence and uniqueness up to stochastic equivalence 
of a solution
of (\ref{1.5.2}) \cite{1}. The second integral on 
the right-hand side of (\ref{1.5.2}) is 
interpreted as an Ito stochastic integral.
Let ${\bf x}_0$ be an $n$-dimensional random variable, which is 
${\rm F}_0$-measurable and 
${\sf M}\{\left|{\bf x}_0\right|^2\}<\infty$ 
(${\sf M}$ denotes a mathematical expectation).
We assume that
${\bf x}_0$ and ${\bf f}_t-{\bf f}_0$ are independent when $t>0.$

It is well known that one of the effective approaches 
to the numerical integration of 
Ito SDEs is an approach based on the Taylor--Ito and 
Taylor--Stratonovich expansions
\cite{KlPl2}-\cite{Mi3}. The most important feature of such 
expansions is a presence in them of the so-called iterated
Ito and Stratonovich stochastic integrals, which play the key 
role for solving the 
problem of numerical integration of Ito SDEs and have the 
following form

\vspace{-1mm}
\begin{equation}                    
\label{ito}
J[\psi^{(k)}]_{T,t}=\int\limits_t^T\psi_k(t_k) \ldots \int\limits_t^{t_{2}}
\psi_1(t_1) d{\bf w}_{t_1}^{(i_1)}\ldots
d{\bf w}_{t_k}^{(i_k)},
\end{equation}

\begin{equation}
\label{str}
J^{*}[\psi^{(k)}]_{T,t}=
\int\limits_t^{*T}\psi_k(t_k) \ldots \int\limits_t^{*t_{2}}
\psi_1(t_1) d{\bf w}_{t_1}^{(i_1)}\ldots
d{\bf w}_{t_k}^{(i_k)},
\end{equation}

\vspace{2mm}
\noindent
where every $\psi_l(\tau)$ $(l=1,\ldots,k)$ is a
non-random function 
on $[t,T],$ ${\bf w}_{\tau}^{(i)}={\bf f}_{\tau}^{(i)}$
for $i=1,\ldots,m$ and
${\bf w}_{\tau}^{(0)}=\tau,$
and 
$$
\int\limits\ \hbox{and}\ \int\limits^{*}
$$ 

\vspace{2mm}
\noindent
denote Ito and 
Stratonovich stochastic integrals,
respectively; $i_1,\ldots,i_k = 0, 1,\ldots,m.$
In this paper we use the definition of the Stratonovich 
stochastic integral from \cite{KlPl2}.

Note that $\psi_l(\tau)\equiv 1$ $(l=1,\ldots,k)$ and
$i_1,\ldots,i_k = 0, 1,\ldots,m$ in  
\cite{KlPl2}-\cite{Mi3}. At the same time
$\psi_l(\tau)\equiv (t-\tau)^{q_l}$ ($l=1,\ldots,k$; 
$q_1,\ldots,q_k=0, 1, 2,\ldots $) and $i_1,\ldots,i_k = 1,\ldots,m$ in
\cite{3}-\cite{20xx2}.

The construction of 
effective expansions (that converge in the mean-square sense)
for the iterated Stra\-to\-no\-vich stochastic integrals
(\ref{str}) of multiplicity 3
composes the subject of this article.

The problem of effective jointly numerical modeling 
(in accordance with the mean-square convergence criterion) of iterated 
Ito and Stratonovich stochastic integrals 
(\ref{ito}) and (\ref{str}) is 
difficult from 
theoretical and computing point of view \cite{KlPl2}-\cite{rr}.
The only exception is connected with a narrow particular case, when 
$i_1=\ldots=i_k\ne 0$ and
$\psi_1(\tau),\ldots,\psi_k(\tau)\equiv \psi(\tau)$.
This case allows 
the investigation with using the Ito formula 
\cite{KlPl2}-\cite{Mi3}.

Seems that iterated stochastic integrals can be approximated by multiple 
integral sums of different types \cite{Mi2}, \cite{Mi3}, \cite{Al}. 
However, this approach implies partitioning of the integration interval 
$[t, T]$ of iterated stochastic integrals 
(the length $T-t$ of this interval is a rather small 
value, because it is the integration step of numerical methods for 
Ito SDEs) and according to numerical 
experiments this additional partitioning leads to 
unacceptably high computational cost and accumulation of computation errors
\cite{7}.

In \cite{Mi2} (also see \cite{KlPl2}, \cite{KPS}, \cite{Mi3}, 
\cite{KPW}, \cite{Zapad-9}) 
Milstein G.N. proposed to expand (\ref{ito}), (\ref{str})
into iterated series of products
of standard Gaussian random variables by representing the Brownian bridge
process as the trigonometric Fourier series with random coefficients 
(version of the so-called Karhunen--Loeve expansion).
To obtain the Milstein expansion of (\ref{str}), the truncated Fourier
expansions of components of the Wiener process ${\bf f}_s$ must be
iteratively substituted in the single integrals, and the integrals
must be calculated, starting from the innermost integral.
This is a complicated procedure that does not lead to a general
expansion of (\ref{str}) valid for an arbitrary multiplicity $k.$
For this reason, only expansions of simplest single, double, and triple
stochastic integrals (\ref{str}) were presented 
in \cite{KlPl2}, \cite{KPS}, \cite{KPW}, \cite{Zapad-9} ($k=1, 2, 3$)
and in \cite{Mi2}, \cite{Mi3} ($k=1, 2$) 
for the case $\psi_1(\tau), \psi_2(\tau), \psi_3(\tau)\equiv 1;$ 
$i_1, i_2, i_3=0,1,\ldots,m.$

Moreover, the authors of the works
\cite{KlPl2}
(Sect.~5.8, pp.~202--204), \cite{KPS} (pp.~82-84),
\cite{KPW} (pp.~438-439),  
\cite{Zapad-9} (pp.~263-264) use 
the Wong--Zakai approximation 
\cite{W-Z-1}-\cite{Watanabe} (without rigorous proof) within the frames
of the Milstein approach 
\cite{Mi2} based on the series expansion 
of the Brownian bridge process. See discussion in Sect.~9 of 
this paper for detail.

Note that in \cite{rr} the method (similar to the Milstein
approach) of expansion
of double Ito stochastic integrals (\ref{ito}) 
($k=2;$ $\psi_1(\tau), \psi_2(\tau) \equiv 1;$ $i_1, i_2 =1,\ldots,m$) 
based on the series expansion
of the Wiener process \cite{Lipt}
using Haar basis functions and 
trigonometric basis functions has been considered.

It is necessary to note that the approach based 
on the Karhunen--Loeve expansion \cite{Mi2} excelled 
in several times (or even in several orders) 
the methods of integral sums \cite{Mi2}, \cite{Mi3}, \cite{Al}
considering computational costs in the sense 
of their diminishing.

An alternative strong approximation method was 
proposed for (\ref{str}) in \cite{3}, \cite{4} (also see
\cite{11}-\cite{16}, \cite{19}, \cite{20}, \cite{20xx}-\cite{20xx2}),
where $J^{*}[\psi^{(k)}]_{T,t}$ was represented as the multiple stochastic 
integral
from the certain discontinuous non-random function of $k$ variables, and the 
function
was then expressed as the generalized 
iterated Fourier series by complete systems of 
continuously
differentiable functions that are orthonormal in the space 
$L_2([t, T])$. 
As a result,
the general iterated series expansion of products
of standard Gaussian random variables was obtained in 
\cite{3}, \cite{4} (also see
\cite{11}-\cite{16}, \cite{19}, \cite{20}, \cite{20xx}-\cite{20xx2}) 
for (\ref{str}) with an 
arbitrary multiplicity $k.$
Hereinafter, this method is referred to as the method of 
generalized iterated Fourier series.
It was shown \cite{3}, \cite{4} (also see
\cite{11}-\cite{16}, \cite{19}, \cite{20}, \cite{20xx}-\cite{20xx2})
that the method of 
generalized iterated Fourier series leads to the approach
based on the Karhunen--Loeve expansion \cite{Mi2}
in the case of trigonometric system 
of functions
and to a substantially simpler expansion of (\ref{str}) in the case
of Legendre polynomial system.

Obviously, the 
approach based on the Karhunen--Loeve expansion \cite{Mi2}
and the method of generalized 
iterated Fourier series \cite{3}, \cite{4} (also see
\cite{11}-\cite{16}, \cite{19}, \cite{20}, \cite{20xx}-\cite{20xx2}) lead
to iterated application of the operation of limit transition.
So, these methods may not converge in 
the mean-square sense 
to appropriate integrals (\ref{str}) 
for some methods of series summation.
The mentioned problem not appears in the method, which 
is proposed for (\ref{ito}) in Theorems 1, 2 (see below).

\vspace{5mm}

\section{Method of Expansion of Iterated Ito Stochastic Integrals
Based on Generalized Multiple Fourier Series}

\vspace{5mm}

Let us consider another approach to the expansion of iterated
Ito stochastic integrals (\ref{ito}) 
\cite{7}-\cite{19}, \cite{20}-\cite{31} (the so-called
method of generalized
multiple Fourier series). 
The idea of this method is as follows: 
the iterated Ito stochastic 
integral (\ref{ito}) of the multiplicity $k$ is represented as 
the multiple stochastic 
integral from the certain discontinuous non-random function of $k$ variables 
defined on the hypercube $[t, T]^k$, where $[t, T]$ is the interval of 
integration of the iterated Ito stochastic 
integral (\ref{ito}). Then, 
the indicated
non-random function of $k$ variables is expanded in the hypercube $[t, T]^k$ 
into the generalized 
multiple Fourier series that converges
in the mean-square sense
in the space 
$L_2([t,T]^k)$. After a number of nontrivial transformations we come 
(see Theorems 1, 2 below) to the 
mean-square convergening expansion of 
the iterated Ito stochastic 
integral (\ref{ito})
into the multiple 
series of products
of standard  Gaussian random 
variables. Coefficients of this 
series are coefficients of 
the generalized multiple Fourier series for the mentioned non-random function 
of $k$ variables, which can be calculated using the explicit formula 
regardless of the multiplicity $k$ of 
the iterated Ito stochastic 
integral (\ref{ito}).

Suppose that every $\psi_l(\tau)$ $(l=1,\ldots,k)$ is a 
non-random function from the space $L_2([t, T])$.
Define the following function on the hypercube $[t, T]^k$

\vspace{-1mm}
\begin{equation}
\label{ppp}
K(t_1,\ldots,t_k)=
\begin{cases}
\psi_1(t_1)\ldots \psi_k(t_k)\ &\hbox{for}\ \ t_1<\ldots<t_k\\
~\\
~\\
0\ &\hbox{otherwise}
\end{cases},\ \ \ \ t_1,\ldots,t_k\in[t, T],\ \ \ \ k\ge 2,
\end{equation}

\vspace{3mm}
\noindent
and 
$K(t_1)\equiv\psi_1(t_1)$ for $t_1\in[t, T].$

Suppose that $\{\phi_j(x)\}_{j=0}^{\infty}$
is a complete orthonormal system of functions in the space
$L_2([t, T])$. 
The function $K(t_1,\ldots,t_k)$ belongs to the space 
$L_2([t, T]^k).$
At this situation it is well known that the generalized 
multiple Fourier series 
of $K(t_1,\ldots,t_k)\in L_2([t, T]^k)$ is converging 
to $K(t_1,\ldots,t_k)$ in the hypercube $[t, T]^k$ in 
the mean-square sense, i.e.

\vspace{1mm}
$$
\hbox{\vtop{\offinterlineskip\halign{
\hfil#\hfil\cr
{\rm lim}\cr
$\stackrel{}{{}_{p_1,\ldots,p_k\to \infty}}$\cr
}} }\Biggl\Vert
K(t_1,\ldots,t_k)-
\sum_{j_1=0}^{p_1}\ldots \sum_{j_k=0}^{p_k}
C_{j_k\ldots j_1}\prod_{l=1}^{k} \phi_{j_l}(t_l)\Biggr\Vert_{L_2([t, T]^k)}=0,
$$

\vspace{4mm}
\noindent
where
\begin{equation}
\label{ppppa}
C_{j_k\ldots j_1}=\int\limits_{[t,T]^k}
K(t_1,\ldots,t_k)\prod_{l=1}^{k}\phi_{j_l}(t_l)dt_1\ldots dt_k,
\end{equation}

\vspace{1mm}
$$
\left\Vert f\right\Vert_{L_2([t, T]^k)}=\left(\int\limits_{[t,T]^k}
f^2(t_1,\ldots,t_k)dt_1\ldots dt_k\right)^{1/2}.
$$

\vspace{4mm}

Consider the partition $\{\tau_j\}_{j=0}^N$ of $[t,T]$ such that

\begin{equation}
\label{1111}
t=\tau_0<\ldots <\tau_N=T,\ \ \
\Delta_N=
\hbox{\vtop{\offinterlineskip\halign{
\hfil#\hfil\cr
{\rm max}\cr
$\stackrel{}{{}_{0\le j\le N-1}}$\cr
}} }\Delta\tau_j\to 0\ \ \hbox{if}\ \ N\to \infty,\ \ \
\Delta\tau_j=\tau_{j+1}-\tau_j.
\end{equation}

\vspace{4mm}

{\bf Theorem 1}\ \cite{7} (2006) (also see \cite{8}-\cite{19}, 
\cite{20}-\cite{31}). 
{\it Suppose that
every $\psi_l(\tau)$ $(l=1,\ldots, k)$ is a continuous non-random function on 
$[t, T]$ and
$\{\phi_j(x)\}_{j=0}^{\infty}$ is a complete orthonormal system  
of continuous functions in the space $L_2([t,T]).$ Then

\vspace{-1mm}
$$
J[\psi^{(k)}]_{T,t}\  =\ 
\hbox{\vtop{\offinterlineskip\halign{
\hfil#\hfil\cr
{\rm l.i.m.}\cr
$\stackrel{}{{}_{p_1,\ldots,p_k\to \infty}}$\cr
}} }\sum_{j_1=0}^{p_1}\ldots\sum_{j_k=0}^{p_k}
C_{j_k\ldots j_1}\Biggl(
\prod_{l=1}^k\zeta_{j_l}^{(i_l)}\ -
\Biggr.
$$

\vspace{2mm}
\begin{equation}
\label{tyyy}
-\ \Biggl.
\hbox{\vtop{\offinterlineskip\halign{
\hfil#\hfil\cr
{\rm l.i.m.}\cr
$\stackrel{}{{}_{N\to \infty}}$\cr
}} }\sum_{(l_1,\ldots,l_k)\in {\rm G}_k}
\phi_{j_{1}}(\tau_{l_1})
\Delta{\bf w}_{\tau_{l_1}}^{(i_1)}\ldots
\phi_{j_{k}}(\tau_{l_k})
\Delta{\bf w}_{\tau_{l_k}}^{(i_k)}\Biggr),
\end{equation}

\vspace{5mm}
\noindent
where $J[\psi^{(k)}]_{T,t}$ is defined by {\rm (\ref{ito}),}

$$
{\rm G}_k={\rm H}_k\backslash{\rm L}_k,\ \ \ 
{\rm H}_k=\{(l_1,\ldots,l_k):\ l_1,\ldots,l_k=0,\ 1,\ldots,N-1\},
$$

$$
{\rm L}_k=\{(l_1,\ldots,l_k):\ l_1,\ldots,l_k=0,\ 1,\ldots,N-1;\
l_g\ne l_r\ (g\ne r);\ g, r=1,\ldots,k\},
$$

\vspace{4mm}
\noindent
${\rm l.i.m.}$ is a limit in the mean-square sense$,$
$i_1,\ldots,i_k=0,1,\ldots,m,$

\vspace{-1mm}
\begin{equation}
\label{rr23}
\zeta_{j}^{(i)}=
\int\limits_t^T \phi_{j}(s) d{\bf w}_s^{(i)}
\end{equation} 

\vspace{2mm}
\noindent
are independent standard Gaussian random variables
for various
$i$ or $j$ {\rm(}if $i\ne 0${\rm),}
$C_{j_k\ldots j_1}$ is the Fourier coefficient {\rm(\ref{ppppa}),}
$\Delta{\bf w}_{\tau_{j}}^{(i)}=
{\bf w}_{\tau_{j+1}}^{(i)}-{\bf w}_{\tau_{j}}^{(i)}$
$(i=0, 1,\ldots,m),$
$\left\{\tau_{j}\right\}_{j=0}^{N}$ is the partition of
the interval $[t, T],$ which satisfies the condition {\rm (\ref{1111})}.
}

\vspace{2mm}

It was shown 
that Theorem 1 is valid for convergence 
in the mean of degree $2n$ ($n\in \mathbb{N}$) 
\cite{20xx} (Sect.~1.1.9, 1.11, 1.12), \cite{26a} (Sect.~6, 15, 16)
and for convergence with probablity 1 (w.~p.~1) \cite{20xx} (Sect.~1.7.2), \cite{new-new-2}.
(the case of
complete
orthonormal systems of Legendre polynomials and trigonometric
functions in the space $L_2([t, T])$).
Moreover, the complete orthonormal systems of Haar and 
Rademacher--Walsh functions in the space $L_2([t,T])$ 
can also be applied in Theorem 1
\cite{7}-\cite{16}, \cite{19}, \cite{20}-\cite{20xx2}, \cite{26a}.
The modification of Theorem 1 for 
complete orthonormal with weigth $r(t_1)\ldots r(t_k)\ge 0$ systems
of functions in the space $L_2([t,T]^k)$ can be found in 
\cite{20}-\cite{20xx2}, \cite{26b}.
Note that Theorem 1 and Theorem 2 (see below) 
have been applied to the approximation of iterated
stochastic integrals with respect to the infinite-dimensional
$Q$-Wiener process in \cite{20xx}-\cite{20xx2} (Chapter 7), \cite{30a1}-\cite{31}.

In order to evaluate the significance of Theorem 1 for practice we will
demonstrate its transformed particular cases for 
$k=1,\ldots,6$ \cite{7}-\cite{19}, \cite{20}-\cite{31}

\begin{equation}
\label{a1}
J[\psi^{(1)}]_{T,t}
=\hbox{\vtop{\offinterlineskip\halign{
\hfil#\hfil\cr
{\rm l.i.m.}\cr
$\stackrel{}{{}_{p_1\to \infty}}$\cr
}} }\sum_{j_1=0}^{p_1}
C_{j_1}\zeta_{j_1}^{(i_1)},
\end{equation}

\vspace{2mm}
\begin{equation}
\label{a2}
J[\psi^{(2)}]_{T,t}
=\hbox{\vtop{\offinterlineskip\halign{
\hfil#\hfil\cr
{\rm l.i.m.}\cr
$\stackrel{}{{}_{p_1,p_2\to \infty}}$\cr
}} }\sum_{j_1=0}^{p_1}\sum_{j_2=0}^{p_2}
C_{j_2j_1}\Biggl(\zeta_{j_1}^{(i_1)}\zeta_{j_2}^{(i_2)}
-{\bf 1}_{\{i_1=i_2\ne 0\}}
{\bf 1}_{\{j_1=j_2\}}\Biggr),
\end{equation}

\vspace{5mm}
$$
J[\psi^{(3)}]_{T,t}=
\hbox{\vtop{\offinterlineskip\halign{
\hfil#\hfil\cr
{\rm l.i.m.}\cr
$\stackrel{}{{}_{p_1,\ldots,p_3\to \infty}}$\cr
}} }\sum_{j_1=0}^{p_1}\sum_{j_2=0}^{p_2}\sum_{j_3=0}^{p_3}
C_{j_3j_2j_1}\Biggl(
\zeta_{j_1}^{(i_1)}\zeta_{j_2}^{(i_2)}\zeta_{j_3}^{(i_3)}
-\Biggr.
$$

\vspace{-1mm}
\begin{equation}
\label{a3}
\Biggl.-{\bf 1}_{\{i_1=i_2\ne 0\}}
{\bf 1}_{\{j_1=j_2\}}
\zeta_{j_3}^{(i_3)}
-{\bf 1}_{\{i_2=i_3\ne 0\}}
{\bf 1}_{\{j_2=j_3\}}
\zeta_{j_1}^{(i_1)}-
{\bf 1}_{\{i_1=i_3\ne 0\}}
{\bf 1}_{\{j_1=j_3\}}
\zeta_{j_2}^{(i_2)}\Biggr),
\end{equation}

\vspace{6mm}
$$
J[\psi^{(4)}]_{T,t}
=
\hbox{\vtop{\offinterlineskip\halign{
\hfil#\hfil\cr
{\rm l.i.m.}\cr
$\stackrel{}{{}_{p_1,\ldots,p_4\to \infty}}$\cr
}} }\sum_{j_1=0}^{p_1}\ldots\sum_{j_4=0}^{p_4}
C_{j_4\ldots j_1}\Biggl(
\prod_{l=1}^4\zeta_{j_l}^{(i_l)}
\Biggr.
-
$$
$$
-
{\bf 1}_{\{i_1=i_2\ne 0\}}
{\bf 1}_{\{j_1=j_2\}}
\zeta_{j_3}^{(i_3)}
\zeta_{j_4}^{(i_4)}
-
{\bf 1}_{\{i_1=i_3\ne 0\}}
{\bf 1}_{\{j_1=j_3\}}
\zeta_{j_2}^{(i_2)}
\zeta_{j_4}^{(i_4)}-
$$
$$
-
{\bf 1}_{\{i_1=i_4\ne 0\}}
{\bf 1}_{\{j_1=j_4\}}
\zeta_{j_2}^{(i_2)}
\zeta_{j_3}^{(i_3)}
-
{\bf 1}_{\{i_2=i_3\ne 0\}}
{\bf 1}_{\{j_2=j_3\}}
\zeta_{j_1}^{(i_1)}
\zeta_{j_4}^{(i_4)}-
$$
$$
-
{\bf 1}_{\{i_2=i_4\ne 0\}}
{\bf 1}_{\{j_2=j_4\}}
\zeta_{j_1}^{(i_1)}
\zeta_{j_3}^{(i_3)}
-
{\bf 1}_{\{i_3=i_4\ne 0\}}
{\bf 1}_{\{j_3=j_4\}}
\zeta_{j_1}^{(i_1)}
\zeta_{j_2}^{(i_2)}+
$$
$$
+
{\bf 1}_{\{i_1=i_2\ne 0\}}
{\bf 1}_{\{j_1=j_2\}}
{\bf 1}_{\{i_3=i_4\ne 0\}}
{\bf 1}_{\{j_3=j_4\}}
+
$$
$$
+
{\bf 1}_{\{i_1=i_3\ne 0\}}
{\bf 1}_{\{j_1=j_3\}}
{\bf 1}_{\{i_2=i_4\ne 0\}}
{\bf 1}_{\{j_2=j_4\}}+
$$
\begin{equation}
\label{a4}
+\Biggl.
{\bf 1}_{\{i_1=i_4\ne 0\}}
{\bf 1}_{\{j_1=j_4\}}
{\bf 1}_{\{i_2=i_3\ne 0\}}
{\bf 1}_{\{j_2=j_3\}}\Biggr),
\end{equation}

\vspace{9mm}
$$
J[\psi^{(5)}]_{T,t}
=\hbox{\vtop{\offinterlineskip\halign{
\hfil#\hfil\cr
{\rm l.i.m.}\cr
$\stackrel{}{{}_{p_1,\ldots,p_5\to \infty}}$\cr
}} }\sum_{j_1=0}^{p_1}\ldots\sum_{j_5=0}^{p_5}
C_{j_5\ldots j_1}\Biggl(
\prod_{l=1}^5\zeta_{j_l}^{(i_l)}
-\Biggr.
$$
$$
-
{\bf 1}_{\{i_1=i_2\ne 0\}}
{\bf 1}_{\{j_1=j_2\}}
\zeta_{j_3}^{(i_3)}
\zeta_{j_4}^{(i_4)}
\zeta_{j_5}^{(i_5)}-
{\bf 1}_{\{i_1=i_3\ne 0\}}
{\bf 1}_{\{j_1=j_3\}}
\zeta_{j_2}^{(i_2)}
\zeta_{j_4}^{(i_4)}
\zeta_{j_5}^{(i_5)}-
$$
$$
-
{\bf 1}_{\{i_1=i_4\ne 0\}}
{\bf 1}_{\{j_1=j_4\}}
\zeta_{j_2}^{(i_2)}
\zeta_{j_3}^{(i_3)}
\zeta_{j_5}^{(i_5)}-
{\bf 1}_{\{i_1=i_5\ne 0\}}
{\bf 1}_{\{j_1=j_5\}}
\zeta_{j_2}^{(i_2)}
\zeta_{j_3}^{(i_3)}
\zeta_{j_4}^{(i_4)}-
$$
$$
-
{\bf 1}_{\{i_2=i_3\ne 0\}}
{\bf 1}_{\{j_2=j_3\}}
\zeta_{j_1}^{(i_1)}
\zeta_{j_4}^{(i_4)}
\zeta_{j_5}^{(i_5)}-
{\bf 1}_{\{i_2=i_4\ne 0\}}
{\bf 1}_{\{j_2=j_4\}}
\zeta_{j_1}^{(i_1)}
\zeta_{j_3}^{(i_3)}
\zeta_{j_5}^{(i_5)}-
$$
$$
-
{\bf 1}_{\{i_2=i_5\ne 0\}}
{\bf 1}_{\{j_2=j_5\}}
\zeta_{j_1}^{(i_1)}
\zeta_{j_3}^{(i_3)}
\zeta_{j_4}^{(i_4)}
-{\bf 1}_{\{i_3=i_4\ne 0\}}
{\bf 1}_{\{j_3=j_4\}}
\zeta_{j_1}^{(i_1)}
\zeta_{j_2}^{(i_2)}
\zeta_{j_5}^{(i_5)}-
$$
$$
-
{\bf 1}_{\{i_3=i_5\ne 0\}}
{\bf 1}_{\{j_3=j_5\}}
\zeta_{j_1}^{(i_1)}
\zeta_{j_2}^{(i_2)}
\zeta_{j_4}^{(i_4)}
-{\bf 1}_{\{i_4=i_5\ne 0\}}
{\bf 1}_{\{j_4=j_5\}}
\zeta_{j_1}^{(i_1)}
\zeta_{j_2}^{(i_2)}
\zeta_{j_3}^{(i_3)}+
$$
$$
+
{\bf 1}_{\{i_1=i_2\ne 0\}}
{\bf 1}_{\{j_1=j_2\}}
{\bf 1}_{\{i_3=i_4\ne 0\}}
{\bf 1}_{\{j_3=j_4\}}\zeta_{j_5}^{(i_5)}+
{\bf 1}_{\{i_1=i_2\ne 0\}}
{\bf 1}_{\{j_1=j_2\}}
{\bf 1}_{\{i_3=i_5\ne 0\}}
{\bf 1}_{\{j_3=j_5\}}\zeta_{j_4}^{(i_4)}+
$$
$$
+
{\bf 1}_{\{i_1=i_2\ne 0\}}
{\bf 1}_{\{j_1=j_2\}}
{\bf 1}_{\{i_4=i_5\ne 0\}}
{\bf 1}_{\{j_4=j_5\}}\zeta_{j_3}^{(i_3)}+
{\bf 1}_{\{i_1=i_3\ne 0\}}
{\bf 1}_{\{j_1=j_3\}}
{\bf 1}_{\{i_2=i_4\ne 0\}}
{\bf 1}_{\{j_2=j_4\}}\zeta_{j_5}^{(i_5)}+
$$
$$
+
{\bf 1}_{\{i_1=i_3\ne 0\}}
{\bf 1}_{\{j_1=j_3\}}
{\bf 1}_{\{i_2=i_5\ne 0\}}
{\bf 1}_{\{j_2=j_5\}}\zeta_{j_4}^{(i_4)}+
{\bf 1}_{\{i_1=i_3\ne 0\}}
{\bf 1}_{\{j_1=j_3\}}
{\bf 1}_{\{i_4=i_5\ne 0\}}
{\bf 1}_{\{j_4=j_5\}}\zeta_{j_2}^{(i_2)}+
$$
$$
+
{\bf 1}_{\{i_1=i_4\ne 0\}}
{\bf 1}_{\{j_1=j_4\}}
{\bf 1}_{\{i_2=i_3\ne 0\}}
{\bf 1}_{\{j_2=j_3\}}\zeta_{j_5}^{(i_5)}+
{\bf 1}_{\{i_1=i_4\ne 0\}}
{\bf 1}_{\{j_1=j_4\}}
{\bf 1}_{\{i_2=i_5\ne 0\}}
{\bf 1}_{\{j_2=j_5\}}\zeta_{j_3}^{(i_3)}+
$$
$$
+
{\bf 1}_{\{i_1=i_4\ne 0\}}
{\bf 1}_{\{j_1=j_4\}}
{\bf 1}_{\{i_3=i_5\ne 0\}}
{\bf 1}_{\{j_3=j_5\}}\zeta_{j_2}^{(i_2)}+
{\bf 1}_{\{i_1=i_5\ne 0\}}
{\bf 1}_{\{j_1=j_5\}}
{\bf 1}_{\{i_2=i_3\ne 0\}}
{\bf 1}_{\{j_2=j_3\}}\zeta_{j_4}^{(i_4)}+
$$
$$
+
{\bf 1}_{\{i_1=i_5\ne 0\}}
{\bf 1}_{\{j_1=j_5\}}
{\bf 1}_{\{i_2=i_4\ne 0\}}
{\bf 1}_{\{j_2=j_4\}}\zeta_{j_3}^{(i_3)}+
{\bf 1}_{\{i_1=i_5\ne 0\}}
{\bf 1}_{\{j_1=j_5\}}
{\bf 1}_{\{i_3=i_4\ne 0\}}
{\bf 1}_{\{j_3=j_4\}}\zeta_{j_2}^{(i_2)}+
$$
$$
+
{\bf 1}_{\{i_2=i_3\ne 0\}}
{\bf 1}_{\{j_2=j_3\}}
{\bf 1}_{\{i_4=i_5\ne 0\}}
{\bf 1}_{\{j_4=j_5\}}\zeta_{j_1}^{(i_1)}+
{\bf 1}_{\{i_2=i_4\ne 0\}}
{\bf 1}_{\{j_2=j_4\}}
{\bf 1}_{\{i_3=i_5\ne 0\}}
{\bf 1}_{\{j_3=j_5\}}\zeta_{j_1}^{(i_1)}+
$$
\begin{equation}
\label{a5}
+\Biggl.
{\bf 1}_{\{i_2=i_5\ne 0\}}
{\bf 1}_{\{j_2=j_5\}}
{\bf 1}_{\{i_3=i_4\ne 0\}}
{\bf 1}_{\{j_3=j_4\}}\zeta_{j_1}^{(i_1)}\Biggr),
\end{equation}

\vspace{9mm}

$$
J[\psi^{(6)}]_{T,t}
=\hbox{\vtop{\offinterlineskip\halign{
\hfil#\hfil\cr
{\rm l.i.m.}\cr
$\stackrel{}{{}_{p_1,\ldots,p_6\to \infty}}$\cr
}} }\sum_{j_1=0}^{p_1}\ldots\sum_{j_6=0}^{p_6}
C_{j_6\ldots j_1}\Biggl(
\prod_{l=1}^6
\zeta_{j_l}^{(i_l)}
-\Biggr.
$$
$$
-
{\bf 1}_{\{i_1=i_6\ne 0\}}
{\bf 1}_{\{j_1=j_6\}}
\zeta_{j_2}^{(i_2)}
\zeta_{j_3}^{(i_3)}
\zeta_{j_4}^{(i_4)}
\zeta_{j_5}^{(i_5)}-
{\bf 1}_{\{i_2=i_6\ne 0\}}
{\bf 1}_{\{j_2=j_6\}}
\zeta_{j_1}^{(i_1)}
\zeta_{j_3}^{(i_3)}
\zeta_{j_4}^{(i_4)}
\zeta_{j_5}^{(i_5)}-
$$
$$
-
{\bf 1}_{\{i_3=i_6\ne 0\}}
{\bf 1}_{\{j_3=j_6\}}
\zeta_{j_1}^{(i_1)}
\zeta_{j_2}^{(i_2)}
\zeta_{j_4}^{(i_4)}
\zeta_{j_5}^{(i_5)}-
{\bf 1}_{\{i_4=i_6\ne 0\}}
{\bf 1}_{\{j_4=j_6\}}
\zeta_{j_1}^{(i_1)}
\zeta_{j_2}^{(i_2)}
\zeta_{j_3}^{(i_3)}
\zeta_{j_5}^{(i_5)}-
$$
$$
-
{\bf 1}_{\{i_5=i_6\ne 0\}}
{\bf 1}_{\{j_5=j_6\}}
\zeta_{j_1}^{(i_1)}
\zeta_{j_2}^{(i_2)}
\zeta_{j_3}^{(i_3)}
\zeta_{j_4}^{(i_4)}-
{\bf 1}_{\{i_1=i_2\ne 0\}}
{\bf 1}_{\{j_1=j_2\}}
\zeta_{j_3}^{(i_3)}
\zeta_{j_4}^{(i_4)}
\zeta_{j_5}^{(i_5)}
\zeta_{j_6}^{(i_6)}-
$$
$$
-
{\bf 1}_{\{i_1=i_3\ne 0\}}
{\bf 1}_{\{j_1=j_3\}}
\zeta_{j_2}^{(i_2)}
\zeta_{j_4}^{(i_4)}
\zeta_{j_5}^{(i_5)}
\zeta_{j_6}^{(i_6)}-
{\bf 1}_{\{i_1=i_4\ne 0\}}
{\bf 1}_{\{j_1=j_4\}}
\zeta_{j_2}^{(i_2)}
\zeta_{j_3}^{(i_3)}
\zeta_{j_5}^{(i_5)}
\zeta_{j_6}^{(i_6)}-
$$
$$
-
{\bf 1}_{\{i_1=i_5\ne 0\}}
{\bf 1}_{\{j_1=j_5\}}
\zeta_{j_2}^{(i_2)}
\zeta_{j_3}^{(i_3)}
\zeta_{j_4}^{(i_4)}
\zeta_{j_6}^{(i_6)}-
{\bf 1}_{\{i_2=i_3\ne 0\}}
{\bf 1}_{\{j_2=j_3\}}
\zeta_{j_1}^{(i_1)}
\zeta_{j_4}^{(i_4)}
\zeta_{j_5}^{(i_5)}
\zeta_{j_6}^{(i_6)}-
$$
$$
-
{\bf 1}_{\{i_2=i_4\ne 0\}}
{\bf 1}_{\{j_2=j_4\}}
\zeta_{j_1}^{(i_1)}
\zeta_{j_3}^{(i_3)}
\zeta_{j_5}^{(i_5)}
\zeta_{j_6}^{(i_6)}-
{\bf 1}_{\{i_2=i_5\ne 0\}}
{\bf 1}_{\{j_2=j_5\}}
\zeta_{j_1}^{(i_1)}
\zeta_{j_3}^{(i_3)}
\zeta_{j_4}^{(i_4)}
\zeta_{j_6}^{(i_6)}-
$$
$$
-
{\bf 1}_{\{i_3=i_4\ne 0\}}
{\bf 1}_{\{j_3=j_4\}}
\zeta_{j_1}^{(i_1)}
\zeta_{j_2}^{(i_2)}
\zeta_{j_5}^{(i_5)}
\zeta_{j_6}^{(i_6)}-
{\bf 1}_{\{i_3=i_5\ne 0\}}
{\bf 1}_{\{j_3=j_5\}}
\zeta_{j_1}^{(i_1)}
\zeta_{j_2}^{(i_2)}
\zeta_{j_4}^{(i_4)}
\zeta_{j_6}^{(i_6)}-
$$
$$
-
{\bf 1}_{\{i_4=i_5\ne 0\}}
{\bf 1}_{\{j_4=j_5\}}
\zeta_{j_1}^{(i_1)}
\zeta_{j_2}^{(i_2)}
\zeta_{j_3}^{(i_3)}
\zeta_{j_6}^{(i_6)}+
$$
$$
+
{\bf 1}_{\{i_1=i_2\ne 0\}}
{\bf 1}_{\{j_1=j_2\}}
{\bf 1}_{\{i_3=i_4\ne 0\}}
{\bf 1}_{\{j_3=j_4\}}
\zeta_{j_5}^{(i_5)}
\zeta_{j_6}^{(i_6)}+
{\bf 1}_{\{i_1=i_2\ne 0\}}
{\bf 1}_{\{j_1=j_2\}}
{\bf 1}_{\{i_3=i_5\ne 0\}}
{\bf 1}_{\{j_3=j_5\}}
\zeta_{j_4}^{(i_4)}
\zeta_{j_6}^{(i_6)}+
$$
$$
+
{\bf 1}_{\{i_1=i_2\ne 0\}}
{\bf 1}_{\{j_1=j_2\}}
{\bf 1}_{\{i_4=i_5\ne 0\}}
{\bf 1}_{\{j_4=j_5\}}
\zeta_{j_3}^{(i_3)}
\zeta_{j_6}^{(i_6)}
+
{\bf 1}_{\{i_1=i_3\ne 0\}}
{\bf 1}_{\{j_1=j_3\}}
{\bf 1}_{\{i_2=i_4\ne 0\}}
{\bf 1}_{\{j_2=j_4\}}
\zeta_{j_5}^{(i_5)}
\zeta_{j_6}^{(i_6)}+
$$
$$
+
{\bf 1}_{\{i_1=i_3\ne 0\}}
{\bf 1}_{\{j_1=j_3\}}
{\bf 1}_{\{i_2=i_5\ne 0\}}
{\bf 1}_{\{j_2=j_5\}}
\zeta_{j_4}^{(i_4)}
\zeta_{j_6}^{(i_6)}
+{\bf 1}_{\{i_1=i_3\ne 0\}}
{\bf 1}_{\{j_1=j_3\}}
{\bf 1}_{\{i_4=i_5\ne 0\}}
{\bf 1}_{\{j_4=j_5\}}
\zeta_{j_2}^{(i_2)}
\zeta_{j_6}^{(i_6)}+
$$
$$
+
{\bf 1}_{\{i_1=i_4\ne 0\}}
{\bf 1}_{\{j_1=j_4\}}
{\bf 1}_{\{i_2=i_3\ne 0\}}
{\bf 1}_{\{j_2=j_3\}}
\zeta_{j_5}^{(i_5)}
\zeta_{j_6}^{(i_6)}
+
{\bf 1}_{\{i_1=i_4\ne 0\}}
{\bf 1}_{\{j_1=j_4\}}
{\bf 1}_{\{i_2=i_5\ne 0\}}
{\bf 1}_{\{j_2=j_5\}}
\zeta_{j_3}^{(i_3)}
\zeta_{j_6}^{(i_6)}+
$$
$$
+
{\bf 1}_{\{i_1=i_4\ne 0\}}
{\bf 1}_{\{j_1=j_4\}}
{\bf 1}_{\{i_3=i_5\ne 0\}}
{\bf 1}_{\{j_3=j_5\}}
\zeta_{j_2}^{(i_2)}
\zeta_{j_6}^{(i_6)}
+
{\bf 1}_{\{i_1=i_5\ne 0\}}
{\bf 1}_{\{j_1=j_5\}}
{\bf 1}_{\{i_2=i_3\ne 0\}}
{\bf 1}_{\{j_2=j_3\}}
\zeta_{j_4}^{(i_4)}
\zeta_{j_6}^{(i_6)}+
$$
$$
+
{\bf 1}_{\{i_1=i_5\ne 0\}}
{\bf 1}_{\{j_1=j_5\}}
{\bf 1}_{\{i_2=i_4\ne 0\}}
{\bf 1}_{\{j_2=j_4\}}
\zeta_{j_3}^{(i_3)}
\zeta_{j_6}^{(i_6)}
+
{\bf 1}_{\{i_1=i_5\ne 0\}}
{\bf 1}_{\{j_1=j_5\}}
{\bf 1}_{\{i_3=i_4\ne 0\}}
{\bf 1}_{\{j_3=j_4\}}
\zeta_{j_2}^{(i_2)}
\zeta_{j_6}^{(i_6)}+
$$
$$
+
{\bf 1}_{\{i_2=i_3\ne 0\}}
{\bf 1}_{\{j_2=j_3\}}
{\bf 1}_{\{i_4=i_5\ne 0\}}
{\bf 1}_{\{j_4=j_5\}}
\zeta_{j_1}^{(i_1)}
\zeta_{j_6}^{(i_6)}
+
{\bf 1}_{\{i_2=i_4\ne 0\}}
{\bf 1}_{\{j_2=j_4\}}
{\bf 1}_{\{i_3=i_5\ne 0\}}
{\bf 1}_{\{j_3=j_5\}}
\zeta_{j_1}^{(i_1)}
\zeta_{j_6}^{(i_6)}+
$$
$$
+
{\bf 1}_{\{i_2=i_5\ne 0\}}
{\bf 1}_{\{j_2=j_5\}}
{\bf 1}_{\{i_3=i_4\ne 0\}}
{\bf 1}_{\{j_3=j_4\}}
\zeta_{j_1}^{(i_1)}
\zeta_{j_6}^{(i_6)}
+
{\bf 1}_{\{i_6=i_1\ne 0\}}
{\bf 1}_{\{j_6=j_1\}}
{\bf 1}_{\{i_3=i_4\ne 0\}}
{\bf 1}_{\{j_3=j_4\}}
\zeta_{j_2}^{(i_2)}
\zeta_{j_5}^{(i_5)}+
$$
$$
+
{\bf 1}_{\{i_6=i_1\ne 0\}}
{\bf 1}_{\{j_6=j_1\}}
{\bf 1}_{\{i_3=i_5\ne 0\}}
{\bf 1}_{\{j_3=j_5\}}
\zeta_{j_2}^{(i_2)}
\zeta_{j_4}^{(i_4)}
+
{\bf 1}_{\{i_6=i_1\ne 0\}}
{\bf 1}_{\{j_6=j_1\}}
{\bf 1}_{\{i_2=i_5\ne 0\}}
{\bf 1}_{\{j_2=j_5\}}
\zeta_{j_3}^{(i_3)}
\zeta_{j_4}^{(i_4)}+
$$
$$
+
{\bf 1}_{\{i_6=i_1\ne 0\}}
{\bf 1}_{\{j_6=j_1\}}
{\bf 1}_{\{i_2=i_4\ne 0\}}
{\bf 1}_{\{j_2=j_4\}}
\zeta_{j_3}^{(i_3)}
\zeta_{j_5}^{(i_5)}
+
{\bf 1}_{\{i_6=i_1\ne 0\}}
{\bf 1}_{\{j_6=j_1\}}
{\bf 1}_{\{i_4=i_5\ne 0\}}
{\bf 1}_{\{j_4=j_5\}}
\zeta_{j_2}^{(i_2)}
\zeta_{j_3}^{(i_3)}+
$$
$$
+
{\bf 1}_{\{i_6=i_1\ne 0\}}
{\bf 1}_{\{j_6=j_1\}}
{\bf 1}_{\{i_2=i_3\ne 0\}}
{\bf 1}_{\{j_2=j_3\}}
\zeta_{j_4}^{(i_4)}
\zeta_{j_5}^{(i_5)}
+
{\bf 1}_{\{i_6=i_2\ne 0\}}
{\bf 1}_{\{j_6=j_2\}}
{\bf 1}_{\{i_3=i_5\ne 0\}}
{\bf 1}_{\{j_3=j_5\}}
\zeta_{j_1}^{(i_1)}
\zeta_{j_4}^{(i_4)}+
$$
$$
+
{\bf 1}_{\{i_6=i_2\ne 0\}}
{\bf 1}_{\{j_6=j_2\}}
{\bf 1}_{\{i_4=i_5\ne 0\}}
{\bf 1}_{\{j_4=j_5\}}
\zeta_{j_1}^{(i_1)}
\zeta_{j_3}^{(i_3)}
+
{\bf 1}_{\{i_6=i_2\ne 0\}}
{\bf 1}_{\{j_6=j_2\}}
{\bf 1}_{\{i_3=i_4\ne 0\}}
{\bf 1}_{\{j_3=j_4\}}
\zeta_{j_1}^{(i_1)}
\zeta_{j_5}^{(i_5)}+
$$
$$
+
{\bf 1}_{\{i_6=i_2\ne 0\}}
{\bf 1}_{\{j_6=j_2\}}
{\bf 1}_{\{i_1=i_5\ne 0\}}
{\bf 1}_{\{j_1=j_5\}}
\zeta_{j_3}^{(i_3)}
\zeta_{j_4}^{(i_4)}
+
{\bf 1}_{\{i_6=i_2\ne 0\}}
{\bf 1}_{\{j_6=j_2\}}
{\bf 1}_{\{i_1=i_4\ne 0\}}
{\bf 1}_{\{j_1=j_4\}}
\zeta_{j_3}^{(i_3)}
\zeta_{j_5}^{(i_5)}+
$$
$$
+
{\bf 1}_{\{i_6=i_2\ne 0\}}
{\bf 1}_{\{j_6=j_2\}}
{\bf 1}_{\{i_1=i_3\ne 0\}}
{\bf 1}_{\{j_1=j_3\}}
\zeta_{j_4}^{(i_4)}
\zeta_{j_5}^{(i_5)}
+
{\bf 1}_{\{i_6=i_3\ne 0\}}
{\bf 1}_{\{j_6=j_3\}}
{\bf 1}_{\{i_2=i_5\ne 0\}}
{\bf 1}_{\{j_2=j_5\}}
\zeta_{j_1}^{(i_1)}
\zeta_{j_4}^{(i_4)}+
$$
$$
+
{\bf 1}_{\{i_6=i_3\ne 0\}}
{\bf 1}_{\{j_6=j_3\}}
{\bf 1}_{\{i_4=i_5\ne 0\}}
{\bf 1}_{\{j_4=j_5\}}
\zeta_{j_1}^{(i_1)}
\zeta_{j_2}^{(i_2)}
+
{\bf 1}_{\{i_6=i_3\ne 0\}}
{\bf 1}_{\{j_6=j_3\}}
{\bf 1}_{\{i_2=i_4\ne 0\}}
{\bf 1}_{\{j_2=j_4\}}
\zeta_{j_1}^{(i_1)}
\zeta_{j_5}^{(i_5)}+
$$
$$
+
{\bf 1}_{\{i_6=i_3\ne 0\}}
{\bf 1}_{\{j_6=j_3\}}
{\bf 1}_{\{i_1=i_5\ne 0\}}
{\bf 1}_{\{j_1=j_5\}}
\zeta_{j_2}^{(i_2)}
\zeta_{j_4}^{(i_4)}
+
{\bf 1}_{\{i_6=i_3\ne 0\}}
{\bf 1}_{\{j_6=j_3\}}
{\bf 1}_{\{i_1=i_4\ne 0\}}
{\bf 1}_{\{j_1=j_4\}}
\zeta_{j_2}^{(i_2)}
\zeta_{j_5}^{(i_5)}+
$$
$$
+
{\bf 1}_{\{i_6=i_3\ne 0\}}
{\bf 1}_{\{j_6=j_3\}}
{\bf 1}_{\{i_1=i_2\ne 0\}}
{\bf 1}_{\{j_1=j_2\}}
\zeta_{j_4}^{(i_4)}
\zeta_{j_5}^{(i_5)}
+
{\bf 1}_{\{i_6=i_4\ne 0\}}
{\bf 1}_{\{j_6=j_4\}}
{\bf 1}_{\{i_3=i_5\ne 0\}}
{\bf 1}_{\{j_3=j_5\}}
\zeta_{j_1}^{(i_1)}
\zeta_{j_2}^{(i_2)}+
$$
$$
+
{\bf 1}_{\{i_6=i_4\ne 0\}}
{\bf 1}_{\{j_6=j_4\}}
{\bf 1}_{\{i_2=i_5\ne 0\}}
{\bf 1}_{\{j_2=j_5\}}
\zeta_{j_1}^{(i_1)}
\zeta_{j_3}^{(i_3)}
+
{\bf 1}_{\{i_6=i_4\ne 0\}}
{\bf 1}_{\{j_6=j_4\}}
{\bf 1}_{\{i_2=i_3\ne 0\}}
{\bf 1}_{\{j_2=j_3\}}
\zeta_{j_1}^{(i_1)}
\zeta_{j_5}^{(i_5)}+
$$
$$
+
{\bf 1}_{\{i_6=i_4\ne 0\}}
{\bf 1}_{\{j_6=j_4\}}
{\bf 1}_{\{i_1=i_5\ne 0\}}
{\bf 1}_{\{j_1=j_5\}}
\zeta_{j_2}^{(i_2)}
\zeta_{j_3}^{(i_3)}
+
{\bf 1}_{\{i_6=i_4\ne 0\}}
{\bf 1}_{\{j_6=j_4\}}
{\bf 1}_{\{i_1=i_3\ne 0\}}
{\bf 1}_{\{j_1=j_3\}}
\zeta_{j_2}^{(i_2)}
\zeta_{j_5}^{(i_5)}+
$$
$$
+
{\bf 1}_{\{i_6=i_4\ne 0\}}
{\bf 1}_{\{j_6=j_4\}}
{\bf 1}_{\{i_1=i_2\ne 0\}}
{\bf 1}_{\{j_1=j_2\}}
\zeta_{j_3}^{(i_3)}
\zeta_{j_5}^{(i_5)}
+
{\bf 1}_{\{i_6=i_5\ne 0\}}
{\bf 1}_{\{j_6=j_5\}}
{\bf 1}_{\{i_3=i_4\ne 0\}}
{\bf 1}_{\{j_3=j_4\}}
\zeta_{j_1}^{(i_1)}
\zeta_{j_2}^{(i_2)}+
$$
$$
+
{\bf 1}_{\{i_6=i_5\ne 0\}}
{\bf 1}_{\{j_6=j_5\}}
{\bf 1}_{\{i_2=i_4\ne 0\}}
{\bf 1}_{\{j_2=j_4\}}
\zeta_{j_1}^{(i_1)}
\zeta_{j_3}^{(i_3)}
+
{\bf 1}_{\{i_6=i_5\ne 0\}}
{\bf 1}_{\{j_6=j_5\}}
{\bf 1}_{\{i_2=i_3\ne 0\}}
{\bf 1}_{\{j_2=j_3\}}
\zeta_{j_1}^{(i_1)}
\zeta_{j_4}^{(i_4)}+
$$
$$
+
{\bf 1}_{\{i_6=i_5\ne 0\}}
{\bf 1}_{\{j_6=j_5\}}
{\bf 1}_{\{i_1=i_4\ne 0\}}
{\bf 1}_{\{j_1=j_4\}}
\zeta_{j_2}^{(i_2)}
\zeta_{j_3}^{(i_3)}
+
{\bf 1}_{\{i_6=i_5\ne 0\}}
{\bf 1}_{\{j_6=j_5\}}
{\bf 1}_{\{i_1=i_3\ne 0\}}
{\bf 1}_{\{j_1=j_3\}}
\zeta_{j_2}^{(i_2)}
\zeta_{j_4}^{(i_4)}+
$$
$$
+
{\bf 1}_{\{i_6=i_5\ne 0\}}
{\bf 1}_{\{j_6=j_5\}}
{\bf 1}_{\{i_1=i_2\ne 0\}}
{\bf 1}_{\{j_1=j_2\}}
\zeta_{j_3}^{(i_3)}
\zeta_{j_4}^{(i_4)}-
$$
$$
-
{\bf 1}_{\{i_6=i_1\ne 0\}}
{\bf 1}_{\{j_6=j_1\}}
{\bf 1}_{\{i_2=i_5\ne 0\}}
{\bf 1}_{\{j_2=j_5\}}
{\bf 1}_{\{i_3=i_4\ne 0\}}
{\bf 1}_{\{j_3=j_4\}}-
$$
$$
-
{\bf 1}_{\{i_6=i_1\ne 0\}}
{\bf 1}_{\{j_6=j_1\}}
{\bf 1}_{\{i_2=i_4\ne 0\}}
{\bf 1}_{\{j_2=j_4\}}
{\bf 1}_{\{i_3=i_5\ne 0\}}
{\bf 1}_{\{j_3=j_5\}}-
$$
$$
-
{\bf 1}_{\{i_6=i_1\ne 0\}}
{\bf 1}_{\{j_6=j_1\}}
{\bf 1}_{\{i_2=i_3\ne 0\}}
{\bf 1}_{\{j_2=j_3\}}
{\bf 1}_{\{i_4=i_5\ne 0\}}
{\bf 1}_{\{j_4=j_5\}}-
$$
$$
-
{\bf 1}_{\{i_6=i_2\ne 0\}}
{\bf 1}_{\{j_6=j_2\}}
{\bf 1}_{\{i_1=i_5\ne 0\}}
{\bf 1}_{\{j_1=j_5\}}
{\bf 1}_{\{i_3=i_4\ne 0\}}
{\bf 1}_{\{j_3=j_4\}}-
$$
$$
-
{\bf 1}_{\{i_6=i_2\ne 0\}}
{\bf 1}_{\{j_6=j_2\}}
{\bf 1}_{\{i_1=i_4\ne 0\}}
{\bf 1}_{\{j_1=j_4\}}
{\bf 1}_{\{i_3=i_5\ne 0\}}
{\bf 1}_{\{j_3=j_5\}}-
$$
$$
-
{\bf 1}_{\{i_6=i_2\ne 0\}}
{\bf 1}_{\{j_6=j_2\}}
{\bf 1}_{\{i_1=i_3\ne 0\}}
{\bf 1}_{\{j_1=j_3\}}
{\bf 1}_{\{i_4=i_5\ne 0\}}
{\bf 1}_{\{j_4=j_5\}}-
$$
$$
-
{\bf 1}_{\{i_6=i_3\ne 0\}}
{\bf 1}_{\{j_6=j_3\}}
{\bf 1}_{\{i_1=i_5\ne 0\}}
{\bf 1}_{\{j_1=j_5\}}
{\bf 1}_{\{i_2=i_4\ne 0\}}
{\bf 1}_{\{j_2=j_4\}}-
$$
$$
-
{\bf 1}_{\{i_6=i_3\ne 0\}}
{\bf 1}_{\{j_6=j_3\}}
{\bf 1}_{\{i_1=i_4\ne 0\}}
{\bf 1}_{\{j_1=j_4\}}
{\bf 1}_{\{i_2=i_5\ne 0\}}
{\bf 1}_{\{j_2=j_5\}}-
$$
$$
-
{\bf 1}_{\{i_3=i_6\ne 0\}}
{\bf 1}_{\{j_3=j_6\}}
{\bf 1}_{\{i_1=i_2\ne 0\}}
{\bf 1}_{\{j_1=j_2\}}
{\bf 1}_{\{i_4=i_5\ne 0\}}
{\bf 1}_{\{j_4=j_5\}}-
$$
$$
-
{\bf 1}_{\{i_6=i_4\ne 0\}}
{\bf 1}_{\{j_6=j_4\}}
{\bf 1}_{\{i_1=i_5\ne 0\}}
{\bf 1}_{\{j_1=j_5\}}
{\bf 1}_{\{i_2=i_3\ne 0\}}
{\bf 1}_{\{j_2=j_3\}}-
$$
$$
-
{\bf 1}_{\{i_6=i_4\ne 0\}}
{\bf 1}_{\{j_6=j_4\}}
{\bf 1}_{\{i_1=i_3\ne 0\}}
{\bf 1}_{\{j_1=j_3\}}
{\bf 1}_{\{i_2=i_5\ne 0\}}
{\bf 1}_{\{j_2=j_5\}}-
$$
$$
-
{\bf 1}_{\{i_6=i_4\ne 0\}}
{\bf 1}_{\{j_6=j_4\}}
{\bf 1}_{\{i_1=i_2\ne 0\}}
{\bf 1}_{\{j_1=j_2\}}
{\bf 1}_{\{i_3=i_5\ne 0\}}
{\bf 1}_{\{j_3=j_5\}}-
$$
$$
-
{\bf 1}_{\{i_6=i_5\ne 0\}}
{\bf 1}_{\{j_6=j_5\}}
{\bf 1}_{\{i_1=i_4\ne 0\}}
{\bf 1}_{\{j_1=j_4\}}
{\bf 1}_{\{i_2=i_3\ne 0\}}
{\bf 1}_{\{j_2=j_3\}}-
$$
$$
-
{\bf 1}_{\{i_6=i_5\ne 0\}}
{\bf 1}_{\{j_6=j_5\}}
{\bf 1}_{\{i_1=i_2\ne 0\}}
{\bf 1}_{\{j_1=j_2\}}
{\bf 1}_{\{i_3=i_4\ne 0\}}
{\bf 1}_{\{j_3=j_4\}}-
$$
\begin{equation}
\label{a6}
\Biggl.-
{\bf 1}_{\{i_6=i_5\ne 0\}}
{\bf 1}_{\{j_6=j_5\}}
{\bf 1}_{\{i_1=i_3\ne 0\}}
{\bf 1}_{\{j_1=j_3\}}
{\bf 1}_{\{i_2=i_4\ne 0\}}
{\bf 1}_{\{j_2=j_4\}}\Biggr),
\end{equation}

\vspace{5mm}
\noindent
where ${\bf 1}_A$ is the indicator of the set $A$.

Thus, we obtain the following useful possibilities and advantages
of the method of generalized multiple Fourier series.

1. There is the explicit formula (see (\ref{ppppa})) for calculation 
of expansion coefficients 
of the iterated Ito stochastic integral (\ref{ito}) with any
fixed multiplicity $k$. 

2. We have new possibilities for exact calculation of the mean-square 
error of approximation 
of the iterated Ito stochastic integral (\ref{ito})
(see \cite{17}, \cite{19}, \cite{20}-\cite{20xx2}, \cite{26}).

3. Since the used
multiple Fourier series is a generalized in the sense
that it is built using various complete orthonormal
systems of functions in the space $L_2([t, T])$, then we 
have new possibilities 
for approximation --- we can
use not only trigonometric functions as in \cite{KlPl2}-\cite{Mi3}
but Legendre polynomials.

4. As it turned out (see \cite{3}-\cite{19}, \cite{20}-\cite{31}), 
it is more convenient to work 
with Legendre polynomials for construction the approximations 
of iterated Ito stochastic integrals (\ref{ito}). 
Approximations based on the Legendre polynomials essentially simpler 
than their analogues based on the trigonometric functions
(see \cite{3}-\cite{19}, \cite{20}-\cite{31}).
Another advantages of the application of Legendre polynomials 
in the framework of the mentioned problem are considered
in \cite{20xx}-\cite{20xx2} (Sect.~5.3), \cite{29}, \cite{30}.

5. The approach based on the Karhunen--Loeve expansion
of the Brownian bridge process (also see \cite{rr})
leads to 
iterated application of the operation of limit
transition (the operation of limit transition 
is implemented only once in Theorem 1 and Theorem 2 (see below))
starting from the 
second multiplicity (in the general case) 
and third multiplicity (for the case
$\psi_1(\tau), \psi_2(\tau), \psi_3(\tau)\equiv 1;$ 
$i_1, i_2, i_3=0,1,\ldots,m$)
of iterated stochastic integrals.
Multiple series (the operation of limit transition 
is implemented only once) are more convenient 
for approximation than the iterated ones
(iterated application of the operation of limit
transition), 
since partial sums of multiple series converge for any possible case of  
convergence to infinity of their upper limits of summation 
(let us denote them as $p_1,\ldots, p_k$). 
For example,
when $p_1=\ldots=p_k=p\to\infty$. 
For iterated series, the condition $p_1=\ldots=p_k=p\to\infty$ obviously 
does not guarantee the convergence of this series.
However, in 
\cite{KlPl2}
(Sect.~5.8, pp.~202--204), \cite{KPS} (pp.~82-84),
\cite{KPW} (pp.~438-439),  
\cite{Zapad-9} (pp.~263-264) the authors use 
(without rigorous proof)
the condition $p_1=p_2=p_3=p\to\infty$
within the frames of the mentioned approach
based on the Karhunen--Loeve expansion of the Brownian bridge
process \cite{Mi2} together with the Wong--Zakai approximation
\cite{W-Z-1}-\cite{Watanabe} (see discussion
in Sect.~9 of this paper for detail).

Note that the correctness of formulas (\ref{a1})--(\ref{a6}) 
can be 
verified 
by the fact that if 
$i_1=\ldots=i_6=i=1,\ldots,m$
and $\psi_1(\tau),\ldots,\psi_6(\tau)\equiv \psi(\tau)$,
then we can derive from (\ref{a1})--(\ref{a6}) the well known
equalities, which be fulfilled w.~p.~1 
\cite{8}-\cite{16}, \cite{19}, \cite{20}-\cite{20xx2}

$$
J[\psi^{(1)}]_{T,t}
=\frac{1}{1!}\delta_{T,t},
$$

\vspace{1mm}
$$
J[\psi^{(2)}]_{T,t}
=\frac{1}{2!}\left(\delta^2_{T,t}-\Delta_{T,t}\right),\
$$

\vspace{1mm}
$$
J[\psi^{(3)}]_{T,t}
=\frac{1}{3!}\left(\delta_{T,t}^3-3\delta_{T,t}\Delta_{T,t}\right),
$$

\vspace{1mm}
$$
J[\psi^{(4)}]_{T,t}
=\frac{1}{4!}\left(\delta^4_{T,t}-6\delta_{T,t}^2\Delta_{T,t}
+3\Delta^2_{T,t}\right),\
$$

\vspace{1mm}
$$
J[\psi^{(5)}]_{T,t}
=\frac{1}{5!}\left(\delta^5_{T,t}-10\delta_{T,t}^3\Delta_{T,t}
+15\delta_{T,t}\Delta^2_{T,t}\right),
$$

\vspace{1mm}
$$
J[\psi^{(6)}]_{T,t}
=\frac{1}{6!}\left(\delta^6_{T,t}-15\delta_{T,t}^4\Delta_{T,t}
+45\delta_{T,t}^2\Delta^2_{T,t}-15\Delta_{T,t}^3\right),
$$

\vspace{4mm}
\noindent
where 
$$
\delta_{T,t}=\int\limits_t^T\psi(\tau)d{\bf f}_{\tau}^{(i)},\ \ \
\Delta_{T,t}=\int\limits_t^T\psi^2(\tau)d\tau.
$$

\vspace{4mm}

The above relations can be independently  
obtained using the Ito formula and Hermite polynomials.

\vspace{5mm}

\section{Generalization of Theorem 1 to the Case of an Arbitrary 
Complete Ortho\-nor\-mal System of Functions in the Space $L_2([t, T])$
and $\psi_1(\tau),\ldots,\psi_k(\tau)\in L_2([t, T])$}

\vspace{5mm}

For further consideration, let us 
consider the generalization of formulas (\ref{a1})--(\ref{a6})                 
for the case of an arbitrary multiplicity $k$ $(k\in\mathbb{N})$ of 
the iterated Ito stochastic integral $J[\psi^{(k)}]_{T,t}$ defined by (\ref{ito}).
In order to do this, let us
introduce some notations. 
Consider the unordered
set $\{1, 2, \ldots, k\}$ 
and separate it into two parts:
the first part consists of $r$ unordered 
pairs (sequence order of these pairs is also unimportant) and the 
second one consists of the 
remaining $k-2r$ numbers.
So, we have

\vspace{1mm}
\begin{equation}
\label{leto5007}
(\{
\underbrace{\{g_1, g_2\}, \ldots, 
\{g_{2r-1}, g_{2r}\}}_{\small{\hbox{part 1}}}
\},
\{\underbrace{q_1, \ldots, q_{k-2r}}_{\small{\hbox{part 2}}}
\}),
\end{equation}

\vspace{5mm}
\noindent
where 

\vspace{-2mm}
$$
\{g_1, g_2, \ldots, 
g_{2r-1}, g_{2r}, q_1, \ldots, q_{k-2r}\}=\{1, 2, \ldots, k\},
$$

\vspace{4mm}
\noindent
braces   
mean an unordered 
set, and pa\-ren\-the\-ses mean an ordered set.

We will say that (\ref{leto5007}) is a partition 
and consider the sum with respect to all possible
partitions

\vspace{1mm}
\begin{equation}
\label{leto5008}
\sum_{\stackrel{(\{\{g_1, g_2\}, \ldots, 
\{g_{2r-1}, g_{2r}\}\}, \{q_1, \ldots, q_{k-2r}\})}
{{}_{\{g_1, g_2, \ldots, 
g_{2r-1}, g_{2r}, q_1, \ldots, q_{k-2r}\}=\{1, 2, \ldots, k\}}}}
a_{g_1 g_2, \ldots, 
g_{2r-1} g_{2r}, q_1 \ldots q_{k-2r}}.
\end{equation}

\vspace{5mm}

Below there are several examples of sums in the form (\ref{leto5008})

\vspace{2mm}
$$
\sum_{\stackrel{(\{g_1, g_2\})}{{}_{\{g_1, g_2\}=\{1, 2\}}}}
a_{g_1 g_2}=a_{12},
$$

\vspace{4mm}
$$
\sum_{\stackrel{(\{\{g_1, g_2\}, \{g_3, g_4\}\})}
{{}_{\{g_1, g_2, g_3, g_4\}=\{1, 2, 3, 4\}}}}
a_{g_1 g_2 g_3 g_4}=a_{1234} + a_{1324} + a_{2314},
$$

\vspace{4mm}
$$
\sum_{\stackrel{(\{g_1, g_2\}, \{q_1, q_{2}\})}
{{}_{\{g_1, g_2, q_1, q_{2}\}=\{1, 2, 3, 4\}}}}
a_{g_1 g_2, q_1 q_{2}}=
$$

$$
=a_{12,34}+a_{13,24}+a_{14,23}
+a_{23,14}+a_{24,13}+a_{34,12},
$$

\vspace{5mm}
$$
\sum_{\stackrel{(\{g_1, g_2\}, \{q_1, q_{2}, q_3\})}
{{}_{\{g_1, g_2, q_1, q_{2}, q_3\}=\{1, 2, 3, 4, 5\}}}}
a_{g_1 g_2, q_1 q_{2}q_3}
=
$$

$$
=a_{12,345}+a_{13,245}+a_{14,235}
+a_{15,234}+a_{23,145}+a_{24,135}+
$$
$$
+a_{25,134}+a_{34,125}+a_{35,124}+a_{45,123},
$$

\vspace{6mm}
$$
\sum_{\stackrel{(\{\{g_1, g_2\}, \{g_3, g_{4}\}\}, \{q_1\})}
{{}_{\{g_1, g_2, g_3, g_{4}, q_1\}=\{1, 2, 3, 4, 5\}}}}
a_{g_1 g_2, g_3 g_{4},q_1}
=
$$

$$
=
a_{12,34,5}+a_{13,24,5}+a_{14,23,5}+
a_{12,35,4}+a_{13,25,4}+a_{15,23,4}+
$$
$$
+a_{12,54,3}+a_{15,24,3}+a_{14,25,3}+a_{15,34,2}+a_{13,54,2}+a_{14,53,2}+
$$
$$
+
a_{52,34,1}+a_{53,24,1}+a_{54,23,1}.
$$

\vspace{9mm}

Now we can write (\ref{tyyy}) as

\vspace{1mm}

$$
J[\psi^{(k)}]_{T,t}=
\hbox{\vtop{\offinterlineskip\halign{
\hfil#\hfil\cr
{\rm l.i.m.}\cr
$\stackrel{}{{}_{p_1,\ldots,p_k\to \infty}}$\cr
}} }
\sum\limits_{j_1=0}^{p_1}\ldots
\sum\limits_{j_k=0}^{p_k}
C_{j_k\ldots j_1}\Biggl(
\prod_{l=1}^k\zeta_{j_l}^{(i_l)}+\sum\limits_{r=1}^{[k/2]}
(-1)^r \times
\Biggr.
$$

\vspace{3mm}
\begin{equation}
\label{leto6000hh}
\times
\sum_{\stackrel{(\{\{g_1, g_2\}, \ldots, 
\{g_{2r-1}, g_{2r}\}\}, \{q_1, \ldots, q_{k-2r}\})}
{{}_{\{g_1, g_2, \ldots, 
g_{2r-1}, g_{2r}, q_1, \ldots, q_{k-2r}\}=\{1, 2, \ldots, k\}}}}
\prod\limits_{s=1}^r
{\bf 1}_{\{i_{g_{{}_{2s-1}}}=~i_{g_{{}_{2s}}}\ne 0\}}
\Biggl.{\bf 1}_{\{j_{g_{{}_{2s-1}}}=~j_{g_{{}_{2s}}}\}}
\prod_{l=1}^{k-2r}\zeta_{j_{q_l}}^{(i_{q_l})}\Biggr),
\end{equation}

\vspace{5mm}
\noindent
where $[x]$ is an integer part of a real number $x;$
another notations are the same as in Theorem {\bf 1}.

\vspace{2mm}

In particular, from (\ref{leto6000hh}) for $k=5$ we obtain

\vspace{3mm}

$$
J[\psi^{(5)}]_{T,t}=
\hbox{\vtop{\offinterlineskip\halign{
\hfil#\hfil\cr
{\rm l.i.m.}\cr
$\stackrel{}{{}_{p_1,\ldots,p_5\to \infty}}$\cr
}} }\sum_{j_1=0}^{p_1}\ldots\sum_{j_5=0}^{p_5}
C_{j_5\ldots j_1}\Biggl(
\prod_{l=1}^5\zeta_{j_l}^{(i_l)}-\Biggr.
$$

\vspace{2mm}
$$
-
\sum\limits_{\stackrel{(\{g_1, g_2\}, \{q_1, q_{2}, q_3\})}
{{}_{\{g_1, g_2, q_{1}, q_{2}, q_3\}=\{1, 2, 3, 4, 5\}}}}
{\bf 1}_{\{i_{g_{{}_{1}}}=~i_{g_{{}_{2}}}\ne 0\}}
{\bf 1}_{\{j_{g_{{}_{1}}}=~j_{g_{{}_{2}}}\}}
\prod_{l=1}^{3}\zeta_{j_{q_l}}^{(i_{q_l})}+
$$

\vspace{2mm}
$$
+
\sum_{\stackrel{(\{\{g_1, g_2\}, 
\{g_{3}, g_{4}\}\}, \{q_1\})}
{{}_{\{g_1, g_2, g_{3}, g_{4}, q_1\}=\{1, 2, 3, 4, 5\}}}}
{\bf 1}_{\{i_{g_{{}_{1}}}=~i_{g_{{}_{2}}}\ne 0\}}
{\bf 1}_{\{j_{g_{{}_{1}}}=~j_{g_{{}_{2}}}\}}
\Biggl.{\bf 1}_{\{i_{g_{{}_{3}}}=~i_{g_{{}_{4}}}\ne 0\}}
{\bf 1}_{\{j_{g_{{}_{3}}}=~j_{g_{{}_{4}}}\}}
\zeta_{j_{q_1}}^{(i_{q_1})}\Biggr).
$$

\vspace{7mm}
\noindent
The last equality obviously agrees with
(\ref{a5}).

Let us consider the generalization of Theorem 1 for the case
of an arbitrary complete orthonormal system 
of functions in the space $L_2([t,T])$ 
and $\psi_1(\tau),\ldots,\psi_k(\tau)\in L_2([t, T]).$

\vspace{2mm}

{\bf Theorem~2}\ \cite{20xx} (Sect.~1.11), \cite{26a} (Sect.~15).
{\it Suppose that
$\psi_1(\tau),\ldots,\psi_k(\tau)\in L_2([t, T])$ and
$\{\phi_j(x)\}_{j=0}^{\infty}$ is an arbitrary complete orthonormal system  
of functions in the space $L_2([t,T]).$
Then the following expansion

\vspace{1mm}
$$
J[\psi^{(k)}]_{T,t}=
\hbox{\vtop{\offinterlineskip\halign{
\hfil#\hfil\cr
{\rm l.i.m.}\cr
$\stackrel{}{{}_{p_1,\ldots,p_k\to \infty}}$\cr
}} }
\sum\limits_{j_1=0}^{p_1}\ldots
\sum\limits_{j_k=0}^{p_k}
C_{j_k\ldots j_1}\Biggl(
\prod_{l=1}^k\zeta_{j_l}^{(i_l)}+\sum\limits_{r=1}^{[k/2]}
(-1)^r \times
\Biggr.
$$

\vspace{2mm}
\begin{equation}
\label{leto6000}
\times
\sum_{\stackrel{(\{\{g_1, g_2\}, \ldots, 
\{g_{2r-1}, g_{2r}\}\}, \{q_1, \ldots, q_{k-2r}\})}
{{}_{\{g_1, g_2, \ldots, 
g_{2r-1}, g_{2r}, q_1, \ldots, q_{k-2r}\}=\{1, 2, \ldots, k\}}}}
\prod\limits_{s=1}^r
{\bf 1}_{\{i_{g_{{}_{2s-1}}}=~i_{g_{{}_{2s}}}\ne 0\}}
\Biggl.{\bf 1}_{\{j_{g_{{}_{2s-1}}}=~j_{g_{{}_{2s}}}\}}
\prod_{l=1}^{k-2r}\zeta_{j_{q_l}}^{(i_{q_l})}\Biggr)
\end{equation}

\vspace{6mm}
\noindent
con\-verg\-ing in the mean-square sense is valid,
where $[x]$ is an integer part of a real number $x;$
another notations are the same as in Theorem~{\rm 1}.}

\vspace{2mm}

It should be noted that an analogue of Theorem 2 was considered 
in \cite{Rybakov1000}. 
Note that we use another notations 
\cite{20xx} (Sect.~1.11), \cite{26a} (Sect.~15)
in comparison with \cite{Rybakov1000}.
Moreover, the proof of an analogue of Theorem 2
from \cite{Rybakov1000} is somewhat different from the proof given in 
\cite{20xx} (Sect.~1.11), \cite{26a} (Sect.~15).

As it turned out, the adaptation of the method of generalized multiple 
Fourier series (Theorems 1, 2) to the iterated Stratonovich stochastic 
integrals (\ref{str}) leads simpler expansions than
(\ref{a1})--(\ref{a6}).
The article is devoted to deriving the analogues
of Theorems 1, 2 for triple Stratonovich
stochastic integrals from the so-called Taylor--Stratonovich
expansion \cite{KlPl2}. In this work, we 
use triple Fourier--Legendre series
as well as triple trigonometric Fourier series
for construction of expansions of 
the iterated Stratonovich stochastic 
integrals (\ref{str}). At that, we
consider the general case of series summation (Sect.~4--6).

The rest of the article is organized as follows.
In Sect.~4, we formulate and prove Theorem 3
on expansion of iterated Stratonovich stochastic integrals (\ref{str}) 
of third multiplicity
with 
constant weight functions using triple Fourier--Legendre series. 
Sect.~5 is devoted to the generalization of Theorem 3 for the 
case of binomial weight functions.
In Sect.~6, we obtain an analogue
of Theorem 3 using triple trigonometric Fourier series.
Sect.~7 is devoted to modifications of Theorems 3--5.
In Sect.~8, we consider some recent results on expansions of
iterated Stratonovich stochastic integrals of multiplicities 3 to 6.
Sect.~9 is devoted to the discussion of main results of this article
from point of view of the Wong--Zakai approximation
\cite{W-Z-1}-\cite{Watanabe}.

\vspace{5mm}

\section{Expansion of Iterated Stratonovich Stochastic Integrals of 
Multiplicity 3.
The Case of Legendre Polynomials}

\vspace{5mm}

{\bf Theorem 3}\ \cite{12}-\cite{16}, \cite{19}, \cite{20}-\cite{20xx2}. 
{\it Suppose that
$\{\phi_j(x)\}_{j=0}^{\infty}$ is a complete orthonormal
system of Legendre polynomials
in the space $L_2([t, T])$.
Then, for the iterated Stratonovich stochastic integral of 
third multiplicity

\vspace{1mm}
$$
{\int\limits_t^{*}}^T
{\int\limits_t^{*}}^{t_3}
{\int\limits_t^{*}}^{t_2}
d{\bf f}_{t_1}^{(i_1)}
d{\bf f}_{t_2}^{(i_2)}d{\bf f}_{t_3}^{(i_3)}\ \ \ (i_1, i_2, i_3=1,\ldots,m)
$$

\vspace{4mm}
\noindent
the following 
expansion 

\vspace{1mm}
\begin{equation}
\label{feto19001}
{\int\limits_t^{*}}^T
{\int\limits_t^{*}}^{t_3}
{\int\limits_t^{*}}^{t_2}
d{\bf f}_{t_1}^{(i_1)}
d{\bf f}_{t_2}^{(i_2)}d{\bf f}_{t_3}^{(i_3)}\ 
=
\hbox{\vtop{\offinterlineskip\halign{
\hfil#\hfil\cr
{\rm l.i.m.}\cr
$\stackrel{}{{}_{p_1,p_2,p_3\to \infty}}$\cr
}} }\sum_{j_1=0}^{p_1}\sum_{j_2=0}^{p_2}\sum_{j_3=0}^{p_3}
C_{j_3 j_2 j_1}\zeta_{j_1}^{(i_1)}\zeta_{j_2}^{(i_2)}\zeta_{j_3}^{(i_3)}
\end{equation}

\vspace{4mm}
\noindent
converging in the mean-square sense is valid, where

\vspace{-2mm}
$$
C_{j_3 j_2 j_1}=\int\limits_t^T
\phi_{j_3}(s)\int\limits_t^s
\phi_{j_2}(s_1)
\int\limits_t^{s_1}
\phi_{j_1}(s_2)ds_2ds_1ds.
$$
}

\vspace{3mm}

{\bf Proof.} If we prove w.~p.~1 the following equalities

\begin{equation}
\label{ogo12}
\hbox{\vtop{\offinterlineskip\halign{
\hfil#\hfil\cr
{\rm l.i.m.}\cr
$\stackrel{}{{}_{p_1, p_3\to \infty}}$\cr
}} }
\sum\limits_{j_1=0}^{p_1}\sum\limits_{j_3=0}^{p_3}
C_{j_3 j_1 j_1}\zeta_{j_3}^{(i_3)}
=
\frac{1}{4}(T-t)^{3/2}\left(
\zeta_0^{(i_3)}+\frac{1}{\sqrt{3}}\zeta_1^{(i_3)}\right),
\end{equation}

\vspace{3mm}
\begin{equation}
\label{ogo13}
\hbox{\vtop{\offinterlineskip\halign{
\hfil#\hfil\cr
{\rm l.i.m.}\cr
$\stackrel{}{{}_{p_1, p_3\to \infty}}$\cr
}} }
\sum\limits_{j_1=0}^{p_1}\sum\limits_{j_3=0}^{p_3}
C_{j_3 j_3 j_1}\zeta_{j_1}^{(i_1)}
=
\frac{1}{4}(T-t)^{3/2}\left(
\zeta_0^{(i_1)}-\frac{1}{\sqrt{3}}\zeta_1^{(i_1)}\right),
\end{equation}

\vspace{2mm}
\begin{equation}
\label{ogo13a}
\hbox{\vtop{\offinterlineskip\halign{
\hfil#\hfil\cr
{\rm l.i.m.}\cr
$\stackrel{}{{}_{p_1, p_3\to \infty}}$\cr
}} }
\sum\limits_{j_1=0}^{p_1}\sum\limits_{j_3=0}^{p_3}
C_{j_1 j_3 j_1}\zeta_{j_3}^{(i_2)}
=0,
\end{equation}

\vspace{5mm}
\noindent
then in accordance with Theorems 1, 2 (see (\ref{a3})), 
formulas (\ref{ogo12})--(\ref{ogo13a}), 
standard 
relations between iterated Ito and Stratonovich 
stochastic integrals as well as in accordance with 
the formulas (they also follow
from Theorems 1, 2)

\vspace{-1mm}
$$
\frac{1}{2}\int\limits_t^T\int\limits_t^{\tau}dsd{\bf f}_{\tau}^{(i_3)}=
\frac{1}{4}(T-t)^{3/2}\left(
\zeta_0^{(i_3)}+\frac{1}{\sqrt{3}}\zeta_1^{(i_3)}\right)\ \ \ \hbox{w.\ p.\ 1},
$$

$$
\frac{1}{2}\int\limits_t^T\int\limits_t^{\tau}d{\bf f}_{s}^{(i_1)}d\tau=
\frac{1}{4}(T-t)^{3/2}\left(
\zeta_0^{(i_1)}-\frac{1}{\sqrt{3}}\zeta_1^{(i_1)}\right)\ \ \ \hbox{w.\ p.\ 1}
$$

\vspace{2mm}
\noindent
we will have

$$
\int\limits_t^T\int\limits_t^{t_3}\int\limits_t^{t_2}
d{\bf f}_{t_1}^{(i_1)}d{\bf f}_{t_2}^{(i_2)}d{\bf f}_{t_3}^{(i_3)}=
\hbox{\vtop{\offinterlineskip\halign{
\hfil#\hfil\cr
{\rm l.i.m.}\cr
$\stackrel{}{{}_{p_1,p_2,p_3\to \infty}}$\cr
}} }\sum_{j_1=0}^{p_1}\sum_{j_2=0}^{p_2}\sum_{j_3=0}^{p_3}
C_{j_3 j_2 j_1}\zeta_{j_1}^{(i_1)}\zeta_{j_2}^{(i_2)}\zeta_{j_3}^{(i_3)}
-
$$

\vspace{2mm}
$$
-
{\bf 1}_{\{i_1=i_2\}}
\frac{1}{2}\int\limits_t^T\int\limits_t^{\tau}dsd{\bf f}_{\tau}^{(i_3)}-
{\bf 1}_{\{i_2=i_3\}}
\frac{1}{2}\int\limits_t^T\int\limits_t^{\tau}d{\bf f}_{s}^{(i_1)}d\tau.
$$

\vspace{5mm}

It means that the expansion (\ref{feto19001}) will be proved.

First let us prove that

\vspace{-1mm}
\begin{equation}
\label{ogo3}
\sum\limits_{j_1=0}^{\infty}C_{0 j_1 j_1}=\frac{1}{4}(T-t)^{3/2},
\end{equation}

\begin{equation}
\label{ogo4}
\sum\limits_{j_1=0}^{\infty}C_{1 j_1 j_1}=
\frac{1}{4\sqrt{3}}(T-t)^{3/2}.
\end{equation}

\vspace{3mm}

We have
$$
C_{000}=\frac{(T-t)^{3/2}}{6},
$$

\vspace{2mm}
$$
C_{0 j_1 j_1}=\int\limits_t^T\phi_0(s)\int\limits_t^s\phi_{j_1}(s_1)
\int\limits_t^{s_1}\phi_{j_1}(s_2)ds_2ds_1ds
=
$$

\begin{equation}
\label{ogo6}
=\frac{1}{2}\int\limits_t^T\phi_0(s)\left(
\int\limits_t^s\phi_{j_1}(s_1)ds_1\right)^2ds,\ \ \ j_1\ge 1.
\end{equation}

\vspace{3mm}

Here $\phi_j(s)$ looks as follows

\vspace{-1mm}
\begin{equation}
\label{ogo7}
\phi_j(s)=\sqrt{\frac{2j+1}{T-t}}P_j\left(\left(
s-\frac{T+t}{2}\right)\frac{2}{T-t}\right),\ \ \ j\ge 0,
\end{equation}

\vspace{3mm}
\noindent
where $P_j(x)$ is the Legendre polynomial.

Let us substitute (\ref{ogo7}) into (\ref{ogo6}) and
calculate $C_{0 j_1 j_1}$\  $(j_1\ge 1)$

\vspace{2mm}
$$
C_{0 j_1 j_1}=\frac{2j_1+1}{2(T-t)^{3/2}}
\int\limits_t^T
\left(\int\limits_{-1}^{z(s)}
P_{j_1}(y)\frac{T-t}{2}dy\right)^2ds=
$$

\vspace{2mm}
$$
=\frac{(2j_1+1)\sqrt{T-t}}{8}
\int\limits_t^T
\left(\int\limits_{-1}^{z(s)}
\frac{1}{2j_1+1}\left(P_{j_1+1}^{'}(y)-P_{j_1-1}^{'}(y)\right)dy
\right)^2ds=
$$

\vspace{2mm}
\begin{equation}
\label{ogo8}
=\frac{\sqrt{T-t}}{8(2j_1+1)}
\int\limits_t^T\left(P_{j_1+1}(z(s))-P_{j_1-1}(z(s))\right)^2ds,
\end{equation}

\vspace{5mm}
\noindent
where here and further

\vspace{-2mm}
$$
z(s)=\left(s-\frac{T+t}{2}\right)\frac{2}{T-t},
$$ 

\vspace{4mm}
\noindent
and
we used the following well-known properties of the  Legendre polynomials

$$
P_j(y)=\frac{1}{2j+1}\left(P_{j+1}^{'}(y)-P_{j-1}^{'}(y)\right),\ \ \ 
P_j(-1)=(-1)^j,\ \ \ j\ge 1.
$$

\vspace{4mm}
\noindent

Also we denote 

\vspace{-2mm}
$$
\frac{dP_j}{dy}(y)\stackrel{{\rm def}}{=}P_j^{'}(y).
$$

\vspace{4mm}

From (\ref{ogo8}) using the property of orthogonality of the Legendre 
polynomials, we get the following relation

\vspace{-2mm}
$$
C_{0 j_1 j_1}=\frac{(T-t)^{3/2}}{16(2j_1+1)}
\int\limits_{-1}^1\left(P_{j_1+1}^2(y)+P_{j_1-1}^2(y)\right)dy=
$$

\vspace{1mm}
$$
=
\frac{(T-t)^{3/2}}{8(2j_1+1)}
\left(\frac{1}{2j_1+3}+\frac{1}{2j_1-1}\right),
$$

\vspace{5mm}
\noindent
where we used the property

\vspace{-1mm}
$$
\int\limits_{-1}^1 P_j^2(y)dy=\frac{2}{2j+1},\ \ \ j\ge 0.
$$

\vspace{2mm}

Then

$$
\sum\limits_{j_1=0}^{\infty}C_{0 j_1 j_1}=
\frac{(T-t)^{3/2}}{6}+
\frac{(T-t)^{3/2}}{8}
\left(
\sum_{j_1=1}^{\infty}\frac{1}{(2j_1+1)(2j_1+3)}+
\sum_{j_1=1}^{\infty}\frac{1}{4j_1^2-1}\right)=
$$

\vspace{2mm}
$$
=\frac{(T-t)^{3/2}}{6}+\frac{(T-t)^{3/2}}{8}
\left(\sum_{j_1=1}^{\infty}\frac{1}{4j_1^2-1}-\frac{1}{3}
+\sum_{j_1=1}^{\infty}\frac{1}{4j_1^2-1}\right)=
$$

\vspace{2mm}
$$
=\frac{(T-t)^{3/2}}{6}+\frac{(T-t)^{3/2}}{8}
\left(\frac{1}{2}-\frac{1}{3}+\frac{1}{2}\right)=
\frac{(T-t)^{3/2}}{4}.
$$

\vspace{7mm}

The relation (\ref{ogo3}) is proved.

Let us check the correctness of (\ref{ogo4}). 
Let us represent $C_{1 j_1 j_1}$ in the form

$$
C_{1 j_1 j_1}=\frac{1}{2}\int\limits_t^T
\phi_1(s)\left(\int\limits_t^s\phi_{j_1}(s_1)ds_1\right)^2 ds=
$$

\vspace{1mm}
$$
=\frac{(T-t)^{3/2}(2j_1+1)\sqrt{3}}{16}
\int\limits_{-1}^{1}
P_1(y)\left(\int\limits_{-1}^y P_{j_1}(y_1)dy_1\right)^2 dy,\ \ \ j_1\ge 1.
$$

\vspace{4mm}

Since the functions

\vspace{-3mm}
$$
\left(\int\limits_{-1}^y P_{j_1}(y_1)dy_1\right)^2,\ \ \ j_1\ge 1
$$ 

\vspace{4mm}
\noindent
are even, then, correspondently the functions

\vspace{-1mm}
$$
P_1(y)\left(\int\limits_{-1}^y P_{j_1}(y_1)dy_1\right)^2 dy,\ \ \ j_1\ge 1
$$

\vspace{4mm}
\noindent
are uneven. 

It means that $C_{1 j_1 j_1}=0$ $(j_1\ge 1).$ From the other hand

\vspace{-1mm}
$$
C_{100}=\frac{\sqrt{3}(T-t)^{3/2}}{16}
\int\limits_{-1}^1 y(y+1)^2 dy=\frac{(T-t)^{3/2}}{4\sqrt{3}}.
$$

\vspace{2mm}

Then 

\vspace{-2mm}
$$
\sum\limits_{j_1=0}^{\infty}C_{1 j_1 j_1}=C_{100}+
\sum\limits_{j_1=1}^{\infty}C_{1 j_1 j_1}=
\frac{(T-t)^{3/2}}{4\sqrt{3}}.
$$

\vspace{3mm}

The relation (\ref{ogo4}) is proved. 

Let us prove the equality (\ref{ogo12}). Using (\ref{ogo4}), we get

$$
\sum\limits_{j_1=0}^{p_1}\sum\limits_{j_3=0}^{p_3}
C_{j_3 j_1 j_1}\zeta_{j_3}^{(i_3)}=
\sum\limits_{j_1=0}^{p_1}C_{0 j_1 j_1}\zeta_{0}^{(i_3)}+
\frac{(T-t)^{3/2}}{4\sqrt{3}}\zeta_{1}^{(i_3)}+
\sum\limits_{j_1=0}^{p_1}\sum\limits_{j_3=2}^{p_3}
C_{j_3 j_1 j_1}\zeta_{j_3}^{(i_3)}=
$$

\vspace{2mm}
\begin{equation}
\label{ogo15}
=\sum\limits_{j_1=0}^{p_1}C_{0 j_1 j_1}\zeta_{0}^{(i_3)}+
\frac{(T-t)^{3/2}}{4\sqrt{3}}\zeta_{1}^{(i_3)}+
\sum\limits_{j_1=0}^{p_1}\ \ \sum\limits_{j_3=2, j_3 - {\rm even}}^{2j_1+2}
C_{j_3 j_1 j_1}\zeta_{j_3}^{(i_3)}.
\end{equation}

\vspace{4mm}

Since

\vspace{-2mm}
$$
C_{j_3j_1j_1}=\frac{(T-t)^{3/2}(2j_1+1)\sqrt{2j_3+1}}{16}
\int\limits_{-1}^{1}
P_{j_3}(y)\left(\int\limits_{-1}^y P_{j_1}(y_1)dy_1\right)^2 dy
$$

\vspace{4mm}
\noindent
and degree of the polynomial

\vspace{-2mm}
$$
\left(\int\limits_{-1}^y P_{j_1}(y_1)dy_1\right)^2
$$ 

\vspace{4mm}
\noindent
equals
to $2j_1+2$, then 
$C_{j_3j_1j_1}=0$ for $j_3>2j_1+2.$ It explains
the circumstance that we put
$2j_1+2$ instead of $p_3$ on the right-hand side 
of the formula (\ref{ogo15}).

Moreover, the function 

\vspace{-2mm}
$$
\left(\int\limits_{-1}^y P_{j_1}(y_1)dy_1\right)^2
$$

\vspace{3mm}
\noindent
is even. It means that the function

\vspace{-2mm}
$$
P_{j_3}(y)\left(\int\limits_{-1}^y P_{j_1}(y_1)dy_1\right)^2
$$

\vspace{3mm}
\noindent
is uneven
for uneven
$j_3.$ It means that $C_{j_3 j_1j_1}=0$ for 
uneven
$j_3.$
That is why we 
summarize using even
$j_3$ on the right-hand side
of the formula (\ref{ogo15}).

Then we have

\vspace{-1mm}
$$
\sum\limits_{j_1=0}^{p_1}\ \sum\limits_{j_3=2, j_3 - {\rm even}}^{2j_1+2}
C_{j_3 j_1 j_1}\zeta_{j_3}^{(i_3)}=
\sum\limits_{j_3=2, j_3 - {\rm even}}^{2p_1+2}\ \
\sum\limits_{j_1=(j_3-2)/2}^{p_1}
C_{j_3 j_1 j_1}\zeta_{j_3}^{(i_3)}=
$$

\vspace{2mm}
\begin{equation}
\label{ogo16}
=
\sum\limits_{j_3=2, j_3 - {\rm even}}^{2p_1+2}\ \
\sum\limits_{j_1=0}^{p_1}
C_{j_3 j_1 j_1}\zeta_{j_3}^{(i_3)}.
\end{equation}

\vspace{4mm}

We replaced $(j_3-2)/2$ by zero on the right-hand side
of the formula (\ref{ogo16}), since $C_{j_3j_1j_1}=0$ for 
$0\le j_1< (j_3-2)/2.$

Let us substitute (\ref{ogo16}) into (\ref{ogo15})

$$
\sum\limits_{j_1=0}^{p_1}\sum\limits_{j_3=0}^{p_3}
C_{j_3 j_1 j_1}\zeta_{j_3}^{(i_3)}=
\sum\limits_{j_1=0}^{p_1}C_{0 j_1 j_1}\zeta_{0}^{(i_3)}+
\frac{(T-t)^{3/2}}{4\sqrt{3}}\zeta_{1}^{(i_3)}+
$$

\vspace{2mm}
\begin{equation}
\label{ogo17}
+
\sum\limits_{j_3=2, j_3 - {\rm even}}^{2p_1+2}\ \
\sum\limits_{j_1=0}^{p_1}
C_{j_3 j_1 j_1}\zeta_{j_3}^{(i_3)}.
\end{equation}

\vspace{4mm}

It is easy to see that the right-hand side
of the formula (\ref{ogo17}) does not depend on $p_3.$ 

If we prove that

\vspace{-1mm}
\begin{equation}
\label{ogo18}
\hbox{\vtop{\offinterlineskip\halign{
\hfil#\hfil\cr
{\rm lim}\cr
$\stackrel{}{{}_{p_1\to \infty}}$\cr
}} }
{\sf M}\left\{\left(
\sum\limits_{j_1=0}^{p_1}\sum\limits_{j_3=0}^{p_3}
C_{j_3 j_1 j_1}\zeta_{j_3}^{(i_3)}-
\frac{1}{4}(T-t)^{3/2}\left(
\zeta_0^{(i_3)}+\frac{1}{\sqrt{3}}\zeta_1^{(i_3)}\right)\right)^2\right\}=0,
\end{equation}

\vspace{4mm}
\noindent
then the relaion (\ref{ogo12}) will be proved.

Using (\ref{ogo17}) and (\ref{ogo3}), we can write the left-hand side 
of (\ref{ogo18})
in the following form

\vspace{2mm}
$$
\hbox{\vtop{\offinterlineskip\halign{
\hfil#\hfil\cr
{\rm lim}\cr
$\stackrel{}{{}_{p_1\to \infty}}$\cr
}} }
{\sf M}\left\{\left(
\left(\sum\limits_{j_1=0}^{p_1}C_{0j_1j_1}-
\frac{(T-t)^{3/2}}{4}\right)\zeta_0^{(i_3)}+
\sum\limits_{j_3=2, j_3 - {\rm even}}^{2p_1+2}\ \
\sum\limits_{j_1=0}^{p_1}
C_{j_3 j_1 j_1}\zeta_{j_3}^{(i_3)}\right)^2\right\}=
$$

\vspace{2mm}
$$
=\hbox{\vtop{\offinterlineskip\halign{
\hfil#\hfil\cr
{\rm lim}\cr
$\stackrel{}{{}_{p_1\to \infty}}$\cr
}} }\left(\sum\limits_{j_1=0}^{p_1}C_{0j_1j_1}-
\frac{(T-t)^{3/2}}{4}\right)^2+
\hbox{\vtop{\offinterlineskip\halign{
\hfil#\hfil\cr
{\rm lim}\cr
$\stackrel{}{{}_{p_1\to \infty}}$\cr
}} }
\sum\limits_{j_3=2, j_3 - {\rm even}}^{2p_1+2}
\left(\sum\limits_{j_1=0}^{p_1}
C_{j_3 j_1 j_1}\right)^2=
$$

\vspace{2mm}
\begin{equation}
\label{ogo19}
=\hbox{\vtop{\offinterlineskip\halign{
\hfil#\hfil\cr
{\rm lim}\cr
$\stackrel{}{{}_{p_1\to \infty}}$\cr
}} }
\sum\limits_{j_3=2, j_3 - {\rm even}}^{2p_1+2}
\left(\sum\limits_{j_1=0}^{p_1}
C_{j_3 j_1 j_1}\right)^2.
\end{equation}

\vspace{7mm}

If we prove that

\vspace{-3mm}
\begin{equation}
\label{ogo20}
\hbox{\vtop{\offinterlineskip\halign{
\hfil#\hfil\cr
{\rm lim}\cr
$\stackrel{}{{}_{p_1\to \infty}}$\cr
}} }
\sum\limits_{j_3=2, j_3 - {\rm even}}^{2p_1+2}
\left(\sum\limits_{j_1=0}^{p_1}
C_{j_3 j_1 j_1}\right)^2=0,
\end{equation}

\vspace{5mm}
\noindent
then the relation (\ref{ogo12}) will be proved.

We have

\vspace{-1mm}
$$
\sum\limits_{j_3=2, j_3 - {\rm even}}^{2p_1+2}
\left(\sum\limits_{j_1=0}^{p_1}
C_{j_3 j_1 j_1}\right)^2=
$$

\vspace{2mm}
$$
=
\frac{1}{4}
\sum\limits_{j_3=2, j_3 - {\rm even}}^{2p_1+2}
\left(\int\limits_t^T\phi_{j_3}(s)\sum\limits_{j_1=0}^{p_1}
\left(\int\limits_t^s\phi_{j_1}(s_1)ds_1\right)^2ds\right)^2=
$$

\vspace{2mm}
$$
=\frac{1}{4}
\sum\limits_{j_3=2, j_3 - {\rm even}}^{2p_1+2}
\left(\int\limits_t^T\phi_{j_3}(s)\left((s-t)-\sum\limits_{j_1=p_1+1}^{\infty}
\left(\int\limits_t^s\phi_{j_1}(s_1)ds_1\right)^2\right)ds\right)^2=
$$

\vspace{2mm}
$$
=\frac{1}{4}
\sum\limits_{j_3=2, j_3 - {\rm even}}^{2p_1+2}
\left(\int\limits_t^T\phi_{j_3}(s)\sum\limits_{j_1=p_1+1}^{\infty}
\left(\int\limits_t^s\phi_{j_1}(s_1)ds_1\right)^2 ds\right)^2\le
$$

\vspace{2mm}
\begin{equation}
\label{ogo21}
\le\frac{1}{4}
\sum\limits_{j_3=2, j_3 - {\rm even}}^{2p_1+2}
\left(\int\limits_t^T|\phi_{j_3}(s)| \sum\limits_{j_1=p_1+1}^{\infty}
\left(\int\limits_t^s\phi_{j_1}(s_1)ds_1\right)^2 ds\right)^2.
\end{equation}

\vspace{6mm}

Obtaining (\ref{ogo21}), we used 
the Parseval equality in the form

\vspace{-2mm}
\begin{equation}
\label{ogo10}
\sum_{j_1=0}^{\infty}\left(\int\limits_t^s\phi_{j_1}(s_1)ds_1\right)^2=
\int\limits_t^T \left({\bf 1}_{\{s_1<s\}}\right)^2ds_1=s-t
\end{equation}

\vspace{3mm}
\noindent
and a property of othogonality of the Legendre polynomials

\vspace{-1mm}
\begin{equation}
\label{ogo11}
\int\limits_t^T\phi_{j_3}(s)(s-t)ds=0,\ \ \ j_3\ge 2.
\end{equation}

\vspace{2mm}

Then we have for $j_1\in\mathbb{N}$

\vspace{-1mm}
$$
\left(\int\limits_t^s\phi_{j_1}(s_1)ds_1\right)^2=
\frac{(T-t)(2j_1+1)}{4}
\left(\int\limits_{-1}^{z(s)}
P_{j_1}(y)dy\right)^2=
$$

$$
=\frac{T-t}{4(2j_1+1)}
\left(\int\limits_{-1}^{z(s)}
\left(P_{j_1+1}^{'}(y)-P_{j_1-1}^{'}(y)\right)dy
\right)^2=
$$

\vspace{3mm}
$$
=\frac{T-t}{4(2j_1+1)}
\left(P_{j_1+1}\left(z(s)\right)-
P_{j_1-1}\left(z(s)\right)\right)^2
\le
$$

\vspace{3mm}
\begin{equation}
\label{ogo22}
\le
\frac{T-t}{2(2j_1+1)}
\left(P_{j_1+1}^2\left(z(s)\right)+
P_{j_1-1}^2\left(z(s)\right)\right).
\end{equation}

\vspace{7mm}

For the Legendre polynomials the following well-known 
estimate is correct

\vspace{1mm}
\begin{equation}
\label{ogo23}
|P_n(y)|<\frac{K}{\sqrt{n+1}(1-y^2)^{1/4}},\ \ \ 
y\in (-1, 1),\ \ \ n\in \mathbb{N},
\end{equation}

\vspace{5mm}
\noindent
where constant $K$ does not depend on $y$ and $n.$

The estimate (\ref{ogo23}) can be written for the 
function $\phi_n(s)$ in 
the following form

\vspace{1mm}
$$
|\phi_n(s)|< \sqrt{\frac{2n+1}{n+1}}\frac{K}{\sqrt{T-t}}
\frac{1}
{\left(1-z^2(s)\right)^{1/4}}
<
$$

\vspace{2mm}
\begin{equation}
\label{ogo24}
<\frac{K_1}{\sqrt{T-t}}
\frac{1}
{\left(1-z^2(s)\right)^{1/4}},
\end{equation}

\vspace{4mm}
\noindent
where
$n\in\mathbb{N},$ $K_1=K\sqrt{2},$\  $s\in (t, T).$

Let us estimate the right-hand side of (\ref{ogo22}) using the estimate
(\ref{ogo23})

$$
\left(\int\limits_t^s\phi_{j_1}(s_1)ds_1\right)^2 <
\frac{T-t}{2(2j_1+1)}\left(\frac{K^2}{j_1+2}+\frac{K^2}{j_1}\right)
\frac{1}
{(1-(z(s))^2)^{1/2}} <
$$

\begin{equation}
\label{ogo25}
<
\frac{(T-t)K^2}{2j_1^2}
\frac{1}
{(1-(z(s))^2)^{1/2}},
\end{equation}

\vspace{3mm}
\noindent
where $j_1\in\mathbb{N},$ $s\in(t, T).$

Substituting the estimate (\ref{ogo25}) into the relation (\ref{ogo21})
and using in (\ref{ogo21}) the estimate (\ref{ogo24})
for $|\phi_{j_3}(s)|$, we obtain

\vspace{-2mm}
$$
\sum\limits_{j_3=2, j_3 - {\rm even}}^{2p_1+2}
\left(\sum\limits_{j_1=0}^{p_1}
C_{j_3 j_1 j_1}\right)^2<
$$

\vspace{2mm}
$$
<
\frac{(T-t)K^4 K_1^2}{16}
\sum\limits_{j_3=2, j_3 - {\rm even}}^{2p_1+2}
\left(\int\limits_t^T
\frac{ds}
{\left(1-\left(z(s)\right)^2
\right)^{3/4}}\sum\limits_{j_1=p_1+1}^{\infty}\frac{1}{j_1^2}
\right)^2=
$$

\vspace{2mm}
\begin{equation}
\label{ogo26}
=\frac{(T-t)^3K^4 K_1^2(p_1+1)}{64}
\left(\int\limits_{-1}^1
\frac{dy}
{\left(1-y^2\right)^{3/4}}\right)^2\left(
\sum\limits_{j_1=p_1+1}^{\infty}\frac{1}{j_1^2}
\right)^2.
\end{equation}

\vspace{3mm}

Since
\begin{equation}
\label{ogo27}
\int\limits_{-1}^1
\frac{dy}
{\left(1-y^2\right)^{3/4}}<\infty
\end{equation}

\vspace{2mm}
\noindent
and
\begin{equation}
\label{ogo28}
\sum\limits_{j_1=p_1+1}^{\infty}\frac{1}{j_1^2}
\le \int\limits_{p_1}^{\infty}\frac{dx}{x^2}=\frac{1}{p_1},
\end{equation}

\vspace{2mm}
\noindent
then from (\ref{ogo26}) we obtain

\vspace{-1mm}
\begin{equation}
\label{ogo29}
\sum\limits_{j_3=2, j_3 - {\rm even}}^{2p_1+2}
\left(\sum\limits_{j_1=0}^{p_1}
C_{j_3 j_1 j_1}\right)^2<\frac{C(T-t)^3 (p_1+1)}{p_1^2}\ \  \to\  0\ \ \ 
\hbox{if}\ \
p_1\to \infty,
\end{equation}

\vspace{3mm}
\noindent
where constant $C$ does not depend on $p_1$ and $T-t.$
From (\ref{ogo29}) it follows (\ref{ogo20}), and the relation
(\ref{ogo20})
implies the formula (\ref{ogo12}).

Let us prove the equaity (\ref{ogo13}). 
First let us prove that

\vspace{-1mm}
\begin{equation}
\label{ogo30}
\sum\limits_{j_3=0}^{\infty}C_{j_3 j_3 0}=\frac{1}{4}(T-t)^{3/2},
\end{equation}

\begin{equation}
\label{ogo31}
\sum\limits_{j_3=0}^{\infty}C_{j_3 j_3 j_1}=
-\frac{1}{4\sqrt{3}}(T-t)^{3/2}.
\end{equation}

\vspace{2mm}

We have

\vspace{-2mm}
$$
\sum_{j_3=0}^{\infty}C_{j_3 j_3 0}=C_{000}+
\sum_{j_3=1}^{\infty}C_{j_3 j_3 0},
$$

\vspace{1mm}
$$
C_{000}=\frac{(T-t)^{3/2}}{6},
$$

\vspace{3mm}
$$
C_{j_3 j_3 0}=\frac{(T-t)^{3/2}}{16(2j_3+1)}
\int\limits_{-1}^1\left(P_{j_3+1}^2(y)+P_{j_3-1}^2(y)\right)dy=
$$

\vspace{2mm}
$$
=
\frac{(T-t)^{3/2}}{8(2j_3+1)}
\left(\frac{1}{2j_3+3}+\frac{1}{2j_3-1}\right),\ \ \ j_3\ge 1.
$$

\vspace{6mm}

Then

\vspace{1mm}
$$
\sum\limits_{j_3=0}^{\infty}C_{j_3 j_3 0}=
\frac{(T-t)^{3/2}}{6}+
\frac{(T-t)^{3/2}}{8}
\left(
\sum_{j_3=1}^{\infty}\frac{1}{(2j_3+1)(2j_3+3)}+
\sum_{j_3=1}^{\infty}\frac{1}{4j_3^2-1}\right)=
$$

\vspace{2mm}
$$
=\frac{(T-t)^{3/2}}{6}+\frac{(T-t)^{3/2}}{8}
\left(\sum_{j_3=1}^{\infty}\frac{1}{4j_3^2-1}-\frac{1}{3}
+\sum_{j_3=1}^{\infty}\frac{1}{4j_3^2-1}\right)=
$$

\vspace{3mm}
$$
=\frac{(T-t)^{3/2}}{6}+\frac{(T-t)^{3/2}}{8}
\left(\frac{1}{2}-\frac{1}{3}+\frac{1}{2}\right)=
\frac{(T-t)^{3/2}}{4}.
$$

\vspace{8mm}

The relation (\ref{ogo30}) is proved.
Let us check the equality (\ref{ogo31}). We have

\vspace{2mm}
$$
C_{j_3 j_3 j_1}=\int\limits_t^T
\phi_{j_3}(s)\int\limits_t^s
\phi_{j_3}(s_1)\int\limits_t^{s_1}
\phi_{j_1}(s_2)ds_2ds_1ds=
$$

$$
=
\int\limits_t^T\phi_{j_1}(s_2)ds_2
\int\limits_{s_2}^T
\phi_{j_3}(s_1)ds_1\int\limits_{s_1}^T
\phi_{j_3}(s)ds=
$$

$$
=\frac{1}{2}\int\limits_t^T
\phi_{j_1}(s_2)\left(\int\limits_{s_2}^T\phi_{j_3}(s_1)ds_1\right)^2 ds_2=
$$

\begin{equation}
\label{ogo33}
=\frac{(T-t)^{3/2}(2j_3+1)\sqrt{2j_1+1}}{16}
\int\limits_{-1}^{1}
P_{j_1}(y)\left(\int\limits_{y}^1 P_{j_3}(y_1)dy_1\right)^2 dy,\ \ \ j_3\ge 1.
\end{equation}

\vspace{7mm}

Since the functions
$$
\left(\int\limits_{y}^1 P_{j_3}(y_1)dy_1\right)^2,\ \ \ j_3\ge 1
$$ 

\vspace{3mm}
\noindent
are even, then the functions

\vspace{-2mm}
$$
P_1(y)\left(\int\limits_{y}^1 P_{j_3}(y_1)dy_1\right)^2 dy,\ \ \ j_3\ge 1
$$

\vspace{3mm}
\noindent
are uneven. It means that $C_{j_3 j_3 1}=0$ $(j_3\ge 1).$

Moreover,
$$
C_{001}=\frac{\sqrt{3}(T-t)^{3/2}}{16}
\int\limits_{-1}^1 y(1-y)^2 dy=-\frac{(T-t)^{3/2}}{4\sqrt{3}}.
$$

\vspace{3mm}

Then 
$$
\sum\limits_{j_3=0}^{\infty}C_{j_3 j_3 1}=C_{001}+
\sum\limits_{j_3=1}^{\infty}C_{j_3 j_3 1}=
-\frac{(T-t)^{3/2}}{4\sqrt{3}}.
$$

\vspace{4mm}

The relation (\ref{ogo31}) is proved.

Using the obtained results, we have

\vspace{3mm}
$$
\sum\limits_{j_1=0}^{p_1}\sum\limits_{j_3=0}^{p_3}
C_{j_3 j_3 j_1}\zeta_{j_1}^{(i_1)}=
\sum\limits_{j_3=0}^{p_3}C_{j_3 j_3 0}\zeta_{0}^{(i_1)}-
\frac{(T-t)^{3/2}}{4\sqrt{3}}\zeta_{1}^{(i_1)}+
\sum\limits_{j_3=0}^{p_3}\sum\limits_{j_1=2}^{p_1}
C_{j_3 j_3 j_1}\zeta_{j_1}^{(i_1)}=
$$

\vspace{3mm}
\begin{equation}
\label{ogoo5}
=\sum\limits_{j_3=0}^{p_3}C_{j_3 j_3 0}\zeta_{0}^{(i_1)}-
\frac{(T-t)^{3/2}}{4\sqrt{3}}\zeta_{1}^{(i_1)}+
\sum\limits_{j_3=0}^{p_3}\ \ \sum\limits_{j_1=2, j_1 - {\rm even}}^{2j_3+2}
C_{j_3 j_3 j_1}\zeta_{j_1}^{(i_1)}.
\end{equation}

\vspace{5mm}

Since 
$$
C_{j_3j_3j_1}=
\frac{(T-t)^{3/2}(2j_3+1)\sqrt{2j_1+1}}{16}
\int\limits_{-1}^{1}
P_{j_1}(y)\left(\int\limits_{y}^1 P_{j_3}(y_1)dy_1\right)^2 dy,\ \ \ j_3\ge 1,
$$

\vspace{2mm}
\noindent
and degree of the polynomial

\vspace{-2mm}
$$
\left(\int\limits_{y}^1 P_{j_3}(y_1)dy_1\right)^2
$$ 

\vspace{2mm}
\noindent
equals to 
$2j_3+2$, 
then
$C_{j_3j_3j_1}=0$ for $j_1>2j_3+2.$ It explains the circumstance
that we put $2j_3+2$ instead of $p_1$ on the right-hand side 
of the formula (\ref{ogoo5}).

Moreover, the function 

\vspace{-2mm}
$$
\left(\int\limits_{y}^1 P_{j_3}(y_1)dy_1\right)^2
$$

\vspace{3mm}
\noindent
is even.
It means that the 
function  

\vspace{-2mm}
$$
P_{j_1}(y)\left(\int\limits_{y}^1 P_{j_3}(y_1)dy_1\right)^2
$$ 

\vspace{3mm}
\noindent
is uneven for uneven $j_1.$
It means that $C_{j_3 j_3j_1}=0$ for uneven  $j_1.$
It explains the summation with respect to
even $j_1$ on the right-hand side of (\ref{ogoo5}).

Then we have

$$
\sum\limits_{j_3=0}^{p_3}\ \ \sum\limits_{j_1=2, j_1 - {\rm even}}^{2j_3+2}
C_{j_3 j_3 j_1}\zeta_{j_1}^{(i_1)}=
\sum\limits_{j_1=2, j_1 - {\rm even}}^{2p_3+2}\ \ 
\sum\limits_{j_3=(j_1-2)/2}^{p_3}
C_{j_3 j_3 j_1}\zeta_{j_1}^{(i_1)}
=
$$

\vspace{2mm}
\begin{equation}
\label{ogoo11}
=\sum\limits_{j_1=2, j_1 - {\rm even}}^{2p_3+2}\ \ 
\sum\limits_{j_3=0}^{p_3}
C_{j_3 j_3 j_1}\zeta_{j_1}^{(i_1)}.
\end{equation}

\vspace{5mm}

We replaced $(j_1-2)/2$ by zero on the right-hand side
of (\ref{ogoo11}), since $C_{j_3j_3j_1}=0$ for
$0\le j_3< (j_1-2)/2.$

Let us substitute (\ref{ogoo11}) into (\ref{ogoo5})

$$
\sum\limits_{j_1=0}^{p_1}\sum\limits_{j_3=0}^{p_3}
C_{j_3 j_3 j_1}\zeta_{j_1}^{(i_1)}=
\sum\limits_{j_3=0}^{p_3}C_{j_3 j_3 0}\zeta_{0}^{(i_1)}-
\frac{(T-t)^{3/2}}{4\sqrt{3}}\zeta_{1}^{(i_1)}+
$$

\vspace{1mm}

\begin{equation}
\label{ogoo12}
+
\sum\limits_{j_1=2, j_1 - {\rm even}}^{2p_3+2}\ \ 
\sum\limits_{j_3=0}^{p_3}
C_{j_3 j_3 j_1}\zeta_{j_1}^{(i_1)}.
\end{equation}

\vspace{4mm}

It is easy to see that the right-hand side of the formula 
(\ref{ogoo12}) does not depend on $p_1.$

If we prove that

\begin{equation}
\label{ogoo13}
\hbox{\vtop{\offinterlineskip\halign{
\hfil#\hfil\cr
{\rm lim}\cr
$\stackrel{}{{}_{p_3\to \infty}}$\cr
}} }
{\sf M}\left\{\left(
\sum\limits_{j_1=0}^{p_1}\sum\limits_{j_3=0}^{p_3}
C_{j_3 j_3 j_1}\zeta_{j_1}^{(i_1)}-
\frac{1}{4}(T-t)^{3/2}\left(
\zeta_0^{(i_1)}-\frac{1}{\sqrt{3}}\zeta_1^{(i_1)}\right)\right)^2\right\}=0,
\end{equation}

\vspace{5mm}
\noindent
then (\ref{ogo13}) will be proved.

Using (\ref{ogoo12}) and (\ref{ogo30}), (\ref{ogo31}), we can write 
the left-hand side of the formula (\ref{ogoo13}) in the 
following form

\vspace{2mm}
$$
\hbox{\vtop{\offinterlineskip\halign{
\hfil#\hfil\cr
{\rm lim}\cr
$\stackrel{}{{}_{p_3\to \infty}}$\cr
}} }
{\sf M}\left\{\left(
\left(\sum\limits_{j_3=0}^{p_3}C_{j_3j_3 0}-
\frac{(T-t)^{3/2}}{4}\right)\zeta_0^{(i_1)}+
\sum\limits_{j_1=2, j_1 - {\rm even}}^{2p_3+2}\ \
\sum\limits_{j_3=0}^{p_3}
C_{j_3 j_3 j_1}\zeta_{j_1}^{(i_1)}\right)^2\right\}=
$$

\vspace{1mm}
$$
=\hbox{\vtop{\offinterlineskip\halign{
\hfil#\hfil\cr
{\rm lim}\cr
$\stackrel{}{{}_{p_3\to \infty}}$\cr
}} }\left(\sum\limits_{j_3=0}^{p_1}C_{j_3j_3 0}-
\frac{(T-t)^{3/2}}{4}\right)^2+
\hbox{\vtop{\offinterlineskip\halign{
\hfil#\hfil\cr
{\rm lim}\cr
$\stackrel{}{{}_{p_3\to \infty}}$\cr
}} }
\sum\limits_{j_1=2, j_1 - {\rm even}}^{2p_3+2}
\left(\sum\limits_{j_3=0}^{p_3}
C_{j_3 j_3 j_1}\right)^2=
$$

\vspace{2mm}
$$
=\hbox{\vtop{\offinterlineskip\halign{
\hfil#\hfil\cr
{\rm lim}\cr
$\stackrel{}{{}_{p_3\to \infty}}$\cr
}} }
\sum\limits_{j_1=2, j_1 - {\rm even}}^{2p_3+2}
\left(\sum\limits_{j_3=0}^{p_3}
C_{j_3 j_3 j_1}\right)^2.
$$

\vspace{4mm}

If we prove that

\vspace{-2mm}
\begin{equation}
\label{ogoo15}
\hbox{\vtop{\offinterlineskip\halign{
\hfil#\hfil\cr
{\rm lim}\cr
$\stackrel{}{{}_{p_3\to \infty}}$\cr
}} }
\sum\limits_{j_1=2, j_1 - {\rm even}}^{2p_3+2}
\left(\sum\limits_{j_3=0}^{p_3}
C_{j_3 j_3 j_1}\right)^2=0,
\end{equation}

\vspace{5mm}
\noindent
then the relation (\ref{ogo13}) will be proved.

From (\ref{ogo33}) we obtain

\vspace{-1mm}
$$
\sum\limits_{j_1=2, j_1 - {\rm even}}^{2p_3+2}
\left(\sum\limits_{j_3=0}^{p_3}
C_{j_3 j_3 j_1}\right)^2=
$$

\vspace{2mm}
$$
=
\frac{1}{4}
\sum\limits_{j_1=2, j_1 - {\rm even}}^{2p_3+2}
\left(\int\limits_t^T\phi_{j_1}(s_2)\sum\limits_{j_3=0}^{p_3}
\left(\int\limits_{s_2}^T\phi_{j_3}(s_1)ds_1\right)^2ds_2\right)^2=
$$

\vspace{2mm}
$$
=\frac{1}{4}
\sum\limits_{j_1=2, j_1 - {\rm even}}^{2p_3+2}
\left(\int\limits_t^T\phi_{j_1}(s_2)\left((T-s_2)-
\sum\limits_{j_3=p_3+1}^{\infty}
\left(\int\limits_{s_2}^T\phi_{j_3}(s_1)ds_1\right)^2\right)ds_2\right)^2=
$$

\vspace{2mm}
$$
=\frac{1}{4}
\sum\limits_{j_1=2, j_1 - {\rm even}}^{2p_3+2}
\left(\int\limits_t^T\phi_{j_1}(s_2)\sum\limits_{j_3=p_3+1}^{\infty}
\left(\int\limits_{s_2}^T\phi_{j_3}(s_1)ds_1\right)^2 ds_2\right)^2\le
$$

\vspace{2mm}
\begin{equation}
\label{ogoo21}
\le\frac{1}{4}
\sum\limits_{j_1=2, j_1 - {\rm even}}^{2p_3+2}
\left(\int\limits_t^T|\phi_{j_1}(s_2)|\sum\limits_{j_3=p_3+1}^{\infty}
\left(\int\limits_{s_2}^T\phi_{j_3}(s_1)ds_1\right)^2 ds_2\right)^2.
\end{equation}

\vspace{8mm}

In order to get (\ref{ogoo21}) we used 
the Parseval equality in the form

\vspace{-1mm}
\begin{equation}
\label{ogo10ee}
\sum_{j_1=0}^{\infty}\left(\int\limits_s^T\phi_{j_1}(s_1)ds_1\right)^2=
\int\limits_t^T \left({\bf 1}_{\{s<s_1\}}\right)^2ds_1=T-s
\end{equation}

\vspace{3mm}
\noindent
and a property of othogonality of the Legendre polynomials

\vspace{-2mm}
\begin{equation}
\label{ogo11e}
\int\limits_t^T\phi_{j_3}(s)(T-s)ds=0,\ \ \ j_3\ge 2.
\end{equation}

\vspace{3mm}

Then we have for $j_3\in\mathbb{N}$

$$
\left(\int\limits_{s_2}^T\phi_{j_3}(s_1)ds_1\right)^2=
\frac{(T-t)}{4(2j_3+1)}
\left(P_{j_3+1}\left(
z(s_2)\right)-
P_{j_3-1}\left(
z(s_2)\right)\right)^2\le
$$

\vspace{2mm}
$$
\le
\frac{T-t}{2(2j_3+1)}
\left(P_{j_3+1}^2\left(
z(s_2)\right)+
P_{j_3-1}^2\left(
z(s_2)\right)\right)
<
$$

\vspace{2mm}
$$
<\frac{T-t}{2(2j_3+1)}\left(\frac{K^2}{j_3+2}+\frac{K^2}{j_3}\right)
\frac{1}
{(1-(z(s_2))^2)^{1/2}} <
$$

\vspace{2mm}
\begin{equation}
\label{ogoo25}
< \frac{(T-t)K^2}{2j_3^2}
\frac{1}
{(1-(z(s_2))^2)^{1/2}},\ \ \ s\in(t, T).
\end{equation}

\vspace{6mm}

In order to get (\ref{ogoo25}) we used the estimate
(\ref{ogo23}). 

Substituting the estimate (\ref{ogoo25}) into the relation (\ref{ogoo21})
and using in (\ref{ogoo21}) the estimate (\ref{ogo24})
for $|\phi_{j_1}(s_2)|$, we obtain

$$
\sum\limits_{j_1=2, j_1 - {\rm even}}^{2p_3+2}
\left(\sum\limits_{j_3=0}^{p_3}
C_{j_3 j_3 j_1}\right)^2<
$$

\vspace{2mm}
$$
<
\frac{(T-t)K^4 K_1^2}{16}
\sum\limits_{j_1=2, j_1 - {\rm even}}^{2p_3+2}
\left(\int\limits_t^T
\frac{ds_2}
{(1-z^2(s_2))^{3/4}}\sum\limits_{j_3=p_3+1}^{\infty}\frac{1}{j_3^2}
\right)^2=
$$

\vspace{2mm}
\begin{equation}
\label{ogoo26}
=\frac{(T-t)^3K^4 K_1^2(p_3+1)}{64}
\left(\int\limits_{-1}^1
\frac{dy}
{\left(1-y^2\right)^{3/4}}\right)^2\left(
\sum\limits_{j_3=p_3+1}^{\infty}\frac{1}{j_3^2}
\right)^2.
\end{equation}

\vspace{7mm}

Using (\ref{ogo27}) and (\ref{ogo28}) in
(\ref{ogoo26}), we get

\begin{equation}
\label{ogoo29}
\sum\limits_{j_1=2, j_1 - {\rm even}}^{2p_3+2}
\left(\sum\limits_{j_3=0}^{p_3}
C_{j_3 j_3 j_1}\right)^2<\frac{C(T-t)^3 (p_3+1)}{p_3^2}\ \ \to\ 0\ \ \
\hbox{with}\ \ p_3\to \infty,
\end{equation}

\vspace{3mm}
\noindent
where constant $C$ does not depend on $p_3$ and $T-t.$

From (\ref{ogoo29}) it follows (\ref{ogoo15}), and the relation
(\ref{ogoo15})
implies the formula (\ref{ogo13}). The relation (\ref{ogo13}) is proved.

Let us prove the equality (\ref{ogo13a}).
Since $\psi_1(\tau),$ $\psi_2(\tau),$ $\psi_3(\tau)\equiv 1,$
then the following relation 
for the Fourier coefficients is correct

$$
C_{j_1 j_1 j_3}+C_{j_1 j_3 j_1}+C_{j_3 j_1 j_1}=\frac{1}{2}
C_{j_1}^2 C_{j_3},
$$ 

\vspace{4mm}
\noindent
where $C_j=0$ for $j\ge 1$ and $C_0=\sqrt{T-t}.$
Then w.~p.~1

\vspace{1mm}
\begin{equation}
\label{sodom31}
\hbox{\vtop{\offinterlineskip\halign{
\hfil#\hfil\cr
{\rm l.i.m.}\cr
$\stackrel{}{{}_{p_1, p_3\to \infty}}$\cr
}} }
\sum\limits_{j_1=0}^{p_1}\sum\limits_{j_3=0}^{p_3}
C_{j_1 j_3 j_1}\zeta_{j_3}^{(i_2)}=
\hbox{\vtop{\offinterlineskip\halign{
\hfil#\hfil\cr
{\rm l.i.m.}\cr
$\stackrel{}{{}_{p_1, p_3\to \infty}}$\cr
}} }
\sum\limits_{j_1=0}^{p_1}\sum\limits_{j_3=0}^{p_3}
\left(\frac{1}{2}C_{j_1}^2 C_{j_3}-C_{j_1 j_1 j_3}-C_{j_3 j_1 j_1}
\right)\zeta_{j_3}^{(i_2)}.
\end{equation}

\vspace{5mm}

Therefore, considering 
(\ref{ogo12}) and (\ref{ogo13}), we can write 
w.~p.~1 

\vspace{1mm}
$$
\hbox{\vtop{\offinterlineskip\halign{
\hfil#\hfil\cr
{\rm l.i.m.}\cr
$\stackrel{}{{}_{p_1, p_3\to \infty}}$\cr
}} }
\sum\limits_{j_1=0}^{p_1}\sum\limits_{j_3=0}^{p_3}
C_{j_1 j_3 j_1}\zeta_{j_3}^{(i_2)}=
$$

\vspace{2mm}

$$
=
\frac{1}{2}C_0^3\zeta_0^{(i_2)}-
\hbox{\vtop{\offinterlineskip\halign{
\hfil#\hfil\cr
{\rm l.i.m.}\cr
$\stackrel{}{{}_{p_1, p_3\to \infty}}$\cr
}} }
\sum\limits_{j_1=0}^{p_1}\sum\limits_{j_3=0}^{p_3}
C_{j_1 j_1 j_3}\zeta_{j_3}^{(i_2)}-
\hbox{\vtop{\offinterlineskip\halign{
\hfil#\hfil\cr
{\rm l.i.m.}\cr
$\stackrel{}{{}_{p_1, p_3\to \infty}}$\cr
}} }
\sum\limits_{j_1=0}^{p_1}\sum\limits_{j_3=0}^{p_3}
C_{j_3 j_1 j_1}\zeta_{j_3}^{(i_2)}
=
$$

\vspace{2mm}
$$
=\frac{1}{2}(T-t)^{3/2}
\zeta_0^{(i_2)}
-\frac{1}{4}(T-t)^{3/2}\left(
\zeta_0^{(i_2)}-\frac{1}{\sqrt{3}}\zeta_1^{(i_2)}\right)
-
$$

\vspace{2mm}
\begin{equation}
\label{sodom3}
-\frac{1}{4}(T-t)^{3/2}\left(
\zeta_0^{(i_2)}+\frac{1}{\sqrt{3}}\zeta_1^{(i_2)}\right)=0.
\end{equation}

\vspace{7mm}

The relation (\ref{ogo13a}) is proved. Theorem 3 is proved.

It is easy to see that the formula (\ref{feto19001})  
can be proved for the case $i_1=i_2=i_3$  
using the Ito formula

\vspace{1mm}
$$
\int\limits_t^{*T}\int\limits_t^{*t_3}\int\limits_t^{*t_2}
d{\bf f}_{t_1}^{(i_1)}d{\bf f}_{t_2}^{(i_1)}d{\bf f}_{t_3}^{(i_1)}=
\frac{1}{6}\left(\int\limits_t^T d{\bf f}_{s}^{(i_1)}\right)^3=
\frac{1}{6}\left(C_0\zeta_{0}^{(i_1)}\right)^3=
$$

\vspace{3mm}
$$
=
C_{000}\zeta_{0}^{(i_1)}\zeta_{0}^{(i_1)}\zeta_{0}^{(i_1)},
$$

\vspace{4mm}
\noindent
where the equality is fulfilled w.~p.~1.

\vspace{5mm}

\section{Generalization of Theorem 3}

\vspace{5mm}

Let us consider the following generalization of Theorem 3.

{\bf Theorem 4}\ \cite{12}-\cite{16}, \cite{19}, \cite{20}-\cite{20xx2}. {\it Suppose that
$\{\phi_j(x)\}_{j=0}^{\infty}$ is a complete orthonormal
system of Legendre polynomials
in the space $L_2([t, T])$.
Then, for the iterated Stratonovich stochastic integral of 
third multiplicity

$$
I_{{l_1l_2l_3}_{T,t}}^{*(i_1i_2i_3)}={\int\limits_t^{*}}^T(t-t_3)^{l_3}
{\int\limits_t^{*}}^{t_3}(t-t_2)^{l_2}
{\int\limits_t^{*}}^{t_2}(t-t_1)^{l_1}
d{\bf f}_{t_1}^{(i_1)}
d{\bf f}_{t_2}^{(i_2)}d{\bf f}_{t_3}^{(i_3)}\ \ \ (i_1, i_2, i_3=1,\ldots,m)
$$

\vspace{4mm}
\noindent
the following 
expansion 

\vspace{1mm}
\begin{equation}
\label{feto1900}
I_{{l_1l_2l_3}_{T,t}}^{*(i_1i_2i_3)}=
\hbox{\vtop{\offinterlineskip\halign{
\hfil#\hfil\cr
{\rm l.i.m.}\cr
$\stackrel{}{{}_{p_1,p_2,p_3\to \infty}}$\cr
}} }\sum_{j_1=0}^{p_1}\sum_{j_2=0}^{p_2}\sum_{j_3=0}^{p_3}
C_{j_3 j_2 j_1}\zeta_{j_1}^{(i_1)}\zeta_{j_2}^{(i_2)}\zeta_{j_3}^{(i_3)}
\end{equation}

\vspace{5mm}
\noindent
converging in the mean-square sense is valid 
for each of the following cases

\vspace{2mm}
\noindent
{\rm 1}.\ $i_1\ne i_2,\ i_2\ne i_3,\ i_1\ne i_3$\ and
$l_1, l_2, l_3=0, 1, 2,\ldots $\\
{\rm 2}.\ $i_1=i_2\ne i_3$ and $l_1=l_2\ne l_3$\ and
$l_1, l_2, l_3=0, 1, 2,\ldots $\\
{\rm 3}.\ $i_1\ne i_2=i_3$ and $l_1\ne l_2=l_3$\ and
$l_1, l_2, l_3=0, 1, 2,\ldots $\\
{\rm 4}.\ $i_1, i_2, i_3=1,\ldots,m;$ $l_1=l_2=l_3=l$\ and $l=0, 1, 
2,\ldots,$\\

\vspace{1mm}
\noindent
where

\vspace{-2mm}
$$
C_{j_3 j_2 j_1}=\int\limits_t^T(t-s)^{l_3}\phi_{j_3}(s)
\int\limits_t^s(t-s_1)^{l_2}\phi_{j_2}(s_1)
\int\limits_t^{s_1}(t-s_2)^{l_1}\phi_{j_1}(s_2)ds_2ds_1ds.
$$
}

\vspace{3mm}

{\bf Proof.} Case 1 directly follows from (\ref{a3}).

Let us consider Case 2 ($i_1=i_2\ne i_3$,\ $l_1=l_2=l\ne l_3$\ and
$l_1, l_3=0, 1, 2,\ldots$). So, we prove 
the following expansion 

\vspace{1mm}
\begin{equation}
\label{ogo101}
I_{{l_1 l_1 l_3}_{T,t}}^{*(i_1i_1i_3)}=
\hbox{\vtop{\offinterlineskip\halign{
\hfil#\hfil\cr
{\rm l.i.m.}\cr
$\stackrel{}{{}_{p_1,p_2,p_3\to \infty}}$\cr
}} }\sum_{j_1=0}^{p_1}\sum_{j_2=0}^{p_2}\sum_{j_3=0}^{p_3}
C_{j_3 j_2 j_1}\zeta_{j_1}^{(i_1)}\zeta_{j_2}^{(i_1)}\zeta_{j_3}^{(i_3)}\ \ \
(i_1, i_2, i_3=1,\ldots,m),
\end{equation}

\vspace{5mm}
\noindent
where 
$l, l_3=0, 1, 2,\ldots$, and

\vspace{-1mm}
\begin{equation}
\label{ogo199}
C_{j_3 j_2 j_1}=\int\limits_t^T
\phi_{j_3}(s)(t-s)^{l_3}\int\limits_t^s(t-s_1)^{l}
\phi_{j_2}(s_1)
\int\limits_t^{s_1}(t-s_2)^l
\phi_{j_1}(s_2)ds_2ds_1ds.
\end{equation}

\vspace{4mm}

If we prove w.~p.~1 the formula

\begin{equation}
\label{ogo200}
\hbox{\vtop{\offinterlineskip\halign{
\hfil#\hfil\cr
{\rm l.i.m.}\cr
$\stackrel{}{{}_{p_1, p_3\to \infty}}$\cr
}} }
\sum\limits_{j_1=0}^{p_1}\sum\limits_{j_3=0}^{p_3}
C_{j_3 j_1 j_1}\zeta_{j_3}^{(i_3)}=
\frac{1}{2}\int\limits_t^T(t-s)^{l_3}
\int\limits_t^s(t-s_1)^{2l}ds_1d{\bf f}_s^{(i_3)},
\end{equation}

\vspace{4mm}
\noindent
where 
coefficients $C_{j_3 j_1 j_1}$ has the form (\ref{ogo199}), 
then using Theorems 1, 2 and 
standard relations between iterated
Ito and Stratonovich stochastic integrals we obtain 
the expansion (\ref{ogo101}).

Using Theorems 
1 and 2, we can write 

\vspace{-2mm}
$$
\frac{1}{2}\int\limits_t^T(t-s)^{l_3}
\int\limits_t^s(t-s_1)^{2l}ds_1d{\bf f}_s^{(i_3)}=
\frac{1}{2}\sum\limits_{j_3=0}^{2l+l_3+1}
\tilde C_{j_3}\zeta_{j_3}^{(i_3)}\ \ \ \hbox{w.\ p.\ 1},
$$

\vspace{3mm}
\noindent
where
$$
\tilde C_{j_3}=
\int\limits_t^T
\phi_{j_3}(s)(t-s)^{l_3}\int\limits_t^s(t-s_1)^{2l}ds_1ds.
$$

\vspace{2mm}

Then
$$
\sum\limits_{j_3=0}^{p_3}\sum\limits_{j_1=0}^{p_1}
C_{j_3 j_1 j_1}\zeta_{j_3}^{(i_3)}-
\frac{1}{2}\sum\limits_{j_3=0}^{2l+l_3+1}
\tilde C_{j_3}\zeta_{j_3}^{(i_3)}=
$$

\vspace{2mm}
$$
=\sum\limits_{j_3=0}^{2l+l_3+1}
\left(\sum\limits_{j_1=0}^{p_1}
C_{j_3 j_1 j_1}-\frac{1}{2}\tilde C_{j_3}\right)
\zeta_{j_3}^{(i_3)}+
\sum\limits_{j_3=2l+l_3+2}^{p_3}
\sum\limits_{j_1=0}^{p_1}
C_{j_3 j_1 j_1}\zeta_{j_3}^{(i_3)}.
$$

\vspace{5mm}

Therefore,

\vspace{-1mm}
$$
\hbox{\vtop{\offinterlineskip\halign{
\hfil#\hfil\cr
{\rm lim}\cr
$\stackrel{}{{}_{p_1,p_3\to \infty}}$\cr
}} }
{\sf M}\left\{\left(
\sum\limits_{j_3=0}^{p_3}\sum\limits_{j_1=0}^{p_1}C_{j_3j_1 j_1}
\zeta_{j_3}^{(i_3)}-
\frac{1}{2}\int\limits_t^T(t-s)^{l_3}
\int\limits_t^s(t-s_1)^{2l}ds_1d{\bf f}_s^{(i_3)}\right)^2\right\}=
$$

\vspace{2mm}
\begin{equation}
\label{ogo210}
=\hbox{\vtop{\offinterlineskip\halign{
\hfil#\hfil\cr
{\rm lim}\cr
$\stackrel{}{{}_{p_1\to \infty}}$\cr
}} }\sum\limits_{j_3=0}^{2l+l_3+1}
\left(\sum\limits_{j_1=0}^{p_1}C_{j_3j_1 j_1}-
\frac{1}{2}\tilde C_{j_3}\right)^2+
\hbox{\vtop{\offinterlineskip\halign{
\hfil#\hfil\cr
{\rm lim}\cr
$\stackrel{}{{}_{p_1,p_3\to \infty}}$\cr
}} }{\sf M}\left\{\left(
\sum\limits_{j_3=2l+l_3+2}^{p_3}
\sum\limits_{j_1=0}^{p_1}
C_{j_3 j_1 j_1}\zeta_{j_3}^{(i_3)}\right)^2\right\}.
\end{equation}

\vspace{5mm}

Let us prove that

\vspace{-3mm}
\begin{equation}
\label{ogo211}
\hbox{\vtop{\offinterlineskip\halign{
\hfil#\hfil\cr
{\rm lim}\cr
$\stackrel{}{{}_{p_1\to \infty}}$\cr
}} }
\left(\sum\limits_{j_1=0}^{p_1}C_{j_3j_1 j_1}-
\frac{1}{2}\tilde C_{j_3}\right)^2=0.
\end{equation}

\vspace{3mm}

We have

\vspace{-1mm}
$$
\left(\sum\limits_{j_1=0}^{p_1}C_{j_3j_1 j_1}-
\frac{1}{2}\tilde C_{j_3}\right)^2=
$$

\vspace{1mm}
$$
=\left(\frac{1}{2}\sum\limits_{j_1=0}^{p_1}
\int\limits_t^T\phi_{j_3}(s)(t-s)^{l_3}
\left(\int\limits_t^s\phi_{j_1}(s_1)(t-s_1)^{l}ds_1\right)^2ds-
\frac{1}{2}
\int\limits_t^T
\phi_{j_3}(s)(t-s)^{l_3}\int\limits_t^s(t-s_1)^{2l}ds_1ds\right)^2
$$

\vspace{1mm}
$$
=\frac{1}{4}\left(
\int\limits_t^T\phi_{j_3}(s)(t-s)^{l_3}\left(
\sum\limits_{j_1=0}^{p_1}
\left(\int\limits_t^s\phi_{j_1}(s_1)(t-s_1)^{l}ds_1\right)^2
-\int\limits_t^s(t-s_1)^{2l}ds_1\right)ds\right)^2=
$$

\vspace{1mm}
$$
=\frac{1}{4}\left(
\int\limits_t^T\phi_{j_3}(s)(t-s)^{l_3}\left(
\int\limits_t^s(t-s_1)^{2l}ds_1-\sum\limits_{j_1=p_1+1}^{\infty}
\left(\int\limits_t^s\phi_{j_1}(s_1)(t-s_1)^{l}ds_1\right)^2
-\right.\right.
$$

$$
-\left.\left.
\int\limits_t^s(t-s_1)^{2l}ds_1\right)ds\right)^2=
$$

\vspace{1mm}
\begin{equation}
\label{ogo300}
=\frac{1}{4}\left(
\int\limits_t^T\phi_{j_3}(s)(t-s)^{l_3}
\sum\limits_{j_1=p_1+1}^{\infty}
\left(\int\limits_t^s\phi_{j_1}(s_1)(t-s_1)^{l}ds_1\right)^2
ds\right)^2.
\end{equation}

\vspace{9mm}

In order to get (\ref{ogo300}) we used the Parseval equality, 
which looks as follows

\vspace{-1mm}
\begin{equation}
\label{ogo301}
\sum_{j_1=0}^{\infty}\left(\int\limits_t^s\phi_{j_1}(s_1)
(t-s_1)^lds_1\right)^2=
\int\limits_t^T K^2(s,s_1)ds_1,
\end{equation}

\vspace{3mm}
\noindent
where

\vspace{-2mm}
$$
K(s,s_1)=(t-s_1)^l{\bf 1}_{\{s_1<s\}},\ \ \ s, s_1\in [t, T].
$$

\vspace{5mm}

Taking into account the nondecreasing
of the functional sequence

\vspace{-1mm}
$$
u_n(s)=\sum_{j_1=0}^{n}\left(\int\limits_t^s\phi_{j_1}(s_1)
(t-s_1)^lds_1\right)^2,
$$

\vspace{3mm}
\noindent
continuity of its members and continuity of the limit function

\vspace{-1mm}
$$
u(s)=\int\limits_t^s(t-s_1)^{2l}ds_1
$$ 

\vspace{2mm}
\noindent
at the interval $[t, T]$
in accordance with the Dini Theorem we
have uniform
convergence of the functional sequences $u_n(s)$ to the limit function
$u(s)$ at the interval $[t, T]$.

From (\ref{ogo300}) using the inequality 
of Cauchy--Bunyakovsky, we obtain

$$
\left(\sum\limits_{j_1=0}^{p_1}C_{j_3j_1 j_1}-
\frac{1}{2}\tilde C_{j_3}\right)^2\le
$$

$$
\le
\frac{1}{4}
\int\limits_t^T\phi_{j_3}^2(s)(t-s)^{2l_3}ds
\int\limits_t^T\left(\sum\limits_{j_1=p_1+1}^{\infty}
\left(\int\limits_t^s\phi_{j_1}(s_1)(t-s_1)^{l}ds_1\right)^2
\right)^2 ds\le
$$

\begin{equation}
\label{ogo302}
\le\frac{1}{4}\varepsilon^2 (T-t)^{2l_3}\int\limits_t^T\phi_{j_3}^2(s)ds
(T-t)=\frac{1}{4}(T-t)^{2l_3+1}\varepsilon^2
\end{equation}

\vspace{4mm}
\noindent
when $p_1>N(\varepsilon),$ where $N(\varepsilon)$
exists for any $\varepsilon>0.$
From (\ref{ogo302}) it follows 
(\ref{ogo211}).

Further, 
\begin{equation}
\label{ogo303}
\sum\limits_{j_1=0}^{p_1}
\sum\limits_{j_3=2l+l_3+2}^{p_3}
C_{j_3 j_1 j_1}\zeta_{j_3}^{(i_3)}=
\sum\limits_{j_1=0}^{p_1}
\sum\limits_{j_3=2l+l_3+2}^{2(j_1+l+1)+l_3}
C_{j_3 j_1 j_1}\zeta_{j_3}^{(i_3)}.
\end{equation}

\vspace{3mm}

We put  $2(j_1+l+1)+l_3$ instead of $p_3$, since
$C_{j_3j_1j_1}=0$ for $j_3>2(j_1+l+1)+l_3.$ This conclusion
follows from the relation

\vspace{-1mm}
$$
C_{j_3j_1j_1}=
\frac{1}{2}
\int\limits_t^T\phi_{j_3}(s)(t-s)^{l_3}
\left(
\int\limits_t^s\phi_{j_1}(s_1)(t-s_1)^{l}ds_1\right)^2ds=
$$

$$
=
\frac{1}{2}\int\limits_t^T\phi_{j_3}(s)Q_{2(j_1+l+1)+l_3}(s)ds,
$$

\vspace{3mm}
\noindent
where $Q_{2(j_1+l+1)+l_3}(s)$ is a polynomial of the degree
$2(j_1+l+1)+l_3.$

It is easy to see that

\vspace{-1mm}
\begin{equation}
\label{ogo304}
\sum\limits_{j_1=0}^{p_1}
\sum\limits_{j_3=2l+l_3+2}^{2(j_1+l+1)+l_3}
C_{j_3 j_1 j_1}\zeta_{j_3}^{(i_3)}=
\sum\limits_{j_3=2l+l_3+2}^{2(p_1+l+1)+l_3}
\sum\limits_{j_1=0}^{p_1}
C_{j_3 j_1 j_1}\zeta_{j_3}^{(i_3)}.
\end{equation}

\vspace{2mm}

Note that we included some zero coefficients $C_{j_3 j_1 j_1}$ 
into the sum $\sum\limits_{j_1=0}^{p_1}$.
From (\ref{ogo303}) and (\ref{ogo304}) we have

\vspace{-2mm}
$$
{\sf M}\left\{\left(\sum\limits_{j_1=0}^{p_1}
\sum\limits_{j_3=2l+l_3+2}^{p_3}
C_{j_3 j_1 j_1}\zeta_{j_3}^{(i_3)}\right)^2\right\}=
$$

\vspace{2mm}
$$
=
{\sf M}\left\{\left(
\sum\limits_{j_3=2l+l_3+2}^{2(p_1+l+1)+l_3}
\sum\limits_{j_1=0}^{p_1}
C_{j_3 j_1 j_1}\zeta_{j_3}^{(i_3)}\right)^2\right\}
=\sum\limits_{j_3=2l+l_3+2}^{2(p_1+l+1)+l_3}
\left(\sum\limits_{j_1=0}^{p_1}
C_{j_3 j_1 j_1}\right)^2=
$$

\vspace{2mm}
$$
=
\sum\limits_{j_3=2l+l_3+2}^{2(p_1+l+1)+l_3}
\left(\frac{1}{2}\sum\limits_{j_1=0}^{p_1}
\int\limits_t^T\phi_{j_3}(s)(t-s)^{l_3}
\left(\int\limits_t^s\phi_{j_1}(s_1)(t-s_1)^{l}ds_1\right)^2ds\right)^2=
$$

\vspace{1mm}
$$
=\frac{1}{4}\sum\limits_{j_3=2l+l_3+2}^{2(p_1+l+1)+l_3}
\left(
\int\limits_t^T\phi_{j_3}(s)(t-s)^{l_3}
\sum\limits_{j_1=0}^{p_1}
\left(\int\limits_t^s\phi_{j_1}(s_1)(t-s_1)^{l}ds_1\right)^2
ds\right)^2=
$$

\vspace{1mm}
$$
=\frac{1}{4}\sum\limits_{j_3=2l+l_3+2}^{2(p_1+l+1)+l_3}\left(
\int\limits_t^T\phi_{j_3}(s)(t-s)^{l_3}\left(
\int\limits_t^s(t-s_1)^{2l}ds_1-\sum\limits_{j_1=p_1+1}^{\infty}
\left(\int\limits_t^s\phi_{j_1}(s_1)(t-s_1)^{l}ds_1\right)^2
\right)ds\right)^2
$$

\vspace{1mm}
\begin{equation}
\label{ogo310}
=\frac{1}{4}\sum\limits_{j_3=2l+l_3+2}^{2(p_1+l+1)+l_3}\left(
\int\limits_t^T\phi_{j_3}(s)(t-s)^{l_3}
\sum\limits_{j_1=p_1+1}^{\infty}
\left(\int\limits_t^s\phi_{j_1}(s_1)(t-s_1)^{l}ds_1\right)^2
ds\right)^2.
\end{equation}

\vspace{6mm}

In 
order to 
get
(\ref{ogo310}) we used the Parseval equality 
(\ref{ogo301}) and the following 
relation

\vspace{-2mm}
$$
\int\limits_t^T\phi_{j_3}(s)Q_{2l+1+l_3}(s)ds=0;\ j_3>2l+1+l_3,
$$

\vspace{3mm}
\noindent
where $Q_{2l+1+l_3}(s)$ is a polynomial of degree
$2l+1+l_3.$

Further, we have for $j_1\in\mathbb{N}$

$$
\left(\int\limits_t^s\phi_{j_1}(s_1)(t-s_1)^lds_1\right)^2=
\frac{(T-t)^{2l+1}(2j_1+1)}{2^{2l+2}}
\left(\int\limits_{-1}^{z(s)}
P_{j_1}(y)(1+y)^ldy\right)^2=
$$

$$
=\frac{(T-t)^{2l+1}}{2^{2l+2}(2j_1+1)}\left(
\left(1+z(s)\right)^l
R_{j_1}(s)
-l\int\limits_{-1}^{z(s)}
\left(P_{j_1+1}(y)-P_{j_1-1}(y)\right)\left(1+y\right)^{l-1}dy\right)^2
\le
$$

$$
\le\frac{(T-t)^{2l+1}2}{2^{2l+2}(2j_1+1)}\left(
\left(\frac{2(s-t)}{T-t}\right)^{2l}
R_{j_1}^2(s)
+l^2
\left(
\int\limits_{-1}^{z(s)}
\left(P_{j_1+1}(y)-P_{j_1-1}(y)\right)\left(1+y\right)^{l-1}dy\right)^2
\right)\le
$$

\vspace{1mm}
$$
\le\frac{(T-t)^{2l+1}}{2^{2l+1}(2j_1+1)}\left(
2^{2l+1}
Z_{j_1}(s)+
l^2
\int\limits_{-1}^{z(s)}
(1+y)^{2l-2}dy
\int\limits_{-1}^{z(s)}
\left(P_{j_1+1}(y)-P_{j_1-1}(y)\right)^2dy
\right)\le
$$

\vspace{1mm}
$$
\le\frac{(T-t)^{2l+1}}{2^{2l+1}(2j_1+1)}\left(
2^{2l+1}
Z_{j_1}(s)
+\frac{2l^2}{2l-1}\left(\frac{2(s-t)}{T-t}\right)^{2l-1}
\int\limits_{-1}^{z(s)}
\left(P_{j_1+1}^2(y)+P_{j_1-1}^2(y)\right)dy
\right)\le
$$

\vspace{1mm}
\begin{equation}
\label{ogo400}
\le\frac{(T-t)^{2l+1}}{2(2j_1+1)}\left(
2
Z_{j_1}(s)
+\frac{l^2}{2l-1}
\int\limits_{-1}^{z(s)}
\left(P_{j_1+1}^2(y)+P_{j_1-1}^2(y)\right)dy
\right),
\end{equation}

\vspace{6mm}
\noindent
where
$$
R_{j_1}(s)=P_{j_1+1}(z(s))-P_{j_1-1}(z(s)),
$$

\vspace{-2mm}
$$
Z_{j_1}(s)=P_{j_1+1}^2(z(s))+P_{j_1-1}^2(z(s)).
$$

\vspace{7mm}

Let us estimate the right-hand side 
of (\ref{ogo400}) using (\ref{ogo23})

$$
\left(\int\limits_t^s\phi_{j_1}(s_1)(t-s_1)^lds_1\right)^2 <
$$

$$
<
\frac{(T-t)^{2l+1}}{2(2j_1+1)}\left(\frac{K^2}{j_1+2}+\frac{K^2}{j_1}
\right)\left(\frac{2}
{(1-\left(z(s))^2
\right)^{1/2}}
+\frac{l^2}{2l-1}
\int\limits_{-1}^{z(s)}
\frac{dy}{\left(1-y^2\right)^{1/2}}\right)<
$$

\begin{equation}
\label{ogo401}
<\frac{(T-t)^{2l+1}K^2}{2j_1^2}\left(
\frac{2}
{(1-\left(z(s))^2
\right)^{1/2}}+
\frac{l^2\pi}{2l-1}\right),\ \ \ s\in(t, T).
\end{equation}

\vspace{5mm}

From (\ref{ogo310}) and (\ref{ogo401}) we obtain

$$
{\sf M}\left\{\left(\sum\limits_{j_1=0}^{p_1}
\sum\limits_{j_3=2l+l_3+2}^{p_3}
C_{j_3 j_1 j_1}\zeta_{j_3}^{(i_3)}\right)^2\right\}\le
$$

$$
\le
\frac{1}{4}\sum\limits_{j_3=2l+l_3+2}^{2(p_1+l+1)+l_3}\left(
\int\limits_t^T|\phi_{j_3}(s)|(t-s)^{l_3}
\sum\limits_{j_1=p_1+1}^{\infty}
\left(\int\limits_t^s\phi_{j_1}(s_1)(t-s_1)^{l}ds_1\right)^2
ds\right)^2\le
$$

$$
\le
\frac{1}{4}(T-t)^{2l_3}\sum\limits_{j_3=2l+l_3+2}^{2(p_1+l+1)+l_3}\left(
\int\limits_t^T|\phi_{j_3}(s)|
\sum\limits_{j_1=p_1+1}^{\infty}
\left(\int\limits_t^s\phi_{j_1}(s_1)(t-s_1)^{l}ds_1\right)^2
ds\right)^2<
$$

\vspace{6mm}
$$
<
\frac{(T-t)^{4l+2l_3+1}K^4 K_1^2}{16}\times
$$

$$
\times
\sum\limits_{j_3=2l+l_3+2}^{2(p_1+l+1)+l_3}
\left(\left(\int\limits_t^T
\frac{2ds}
{\left(1-\left(z(s)\right)^2
\right)^{3/4}}
+\frac{l^2\pi}{2l-1}
\int\limits_t^T
\frac{ds}
{\left(1-\left(z(s)\right)^2
\right)^{1/4}}\right)
\sum\limits_{j_1=p_1+1}^{\infty}\frac{1}{j_1^2}
\right)^2\le 
$$

\vspace{4mm}
$$
\le
\frac{(T-t)^{4l+2l_3+3}K^4 K_1^2}{64}\cdot\frac{2p_1+1}{p_1^2}
\left(\int\limits_{-1}^1
\frac{2dy}
{(1-y^2)^{3/4}}
+\frac{l^2\pi}{2l-1}
\int\limits_{-1}^1
\frac{dy}
{(1-y^2)^{1/4}}\right)^2\le
$$

\vspace{4mm}
\begin{equation}
\label{ogo500}
\le (T-t)^{4l+2l_3+3}C\frac{2p_1+1}{p_1^2}\to 0\ \ \ 
\hbox{when}\ \ p_1\ \to\ \infty,
\end{equation}

\vspace{6mm}
\noindent
where constant $C$ does not depend on $p_1$ and $T-t$.

From (\ref{ogo210}), (\ref{ogo211}), and (\ref{ogo500}) it follows 
(\ref{ogo200}), and the relation (\ref{ogo200}) implies
the formula
(\ref{ogo101}).

Let us consider Case 3 ($i_2=i_3\ne i_1$,\ $l_2=l_3=l\ne l_1$,\ and
$l_1, l_3=0, 1, 2,\ldots$). So, we prove 
the following expansion 

\begin{equation}
\label{ogo101ee}
I_{{l_1 l_3 l_3}_{T,t}}^{*(i_1i_3i_3)}=
\hbox{\vtop{\offinterlineskip\halign{
\hfil#\hfil\cr
{\rm l.i.m.}\cr
$\stackrel{}{{}_{p_1,p_2,p_3\to \infty}}$\cr
}} }\sum_{j_1=0}^{p_1}\sum_{j_2=0}^{p_2}\sum_{j_3=0}^{p_3}
C_{j_3 j_2 j_1}\zeta_{j_1}^{(i_1)}\zeta_{j_2}^{(i_3)}\zeta_{j_3}^{(i_3)}\ \ \
(i_1, i_2, i_3=1,\ldots,m),
\end{equation}

\vspace{3mm}
\noindent
where $l, l_1=0, 1, 2,\ldots$, and

\vspace{-1mm}
\begin{equation}
\label{ogo1991}
C_{j_3 j_2 j_1}=\int\limits_t^T
\phi_{j_3}(s)(t-s)^{l}\int\limits_t^s(t-s_1)^{l}
\phi_{j_2}(s_1)
\int\limits_t^{s_1}(t-s_2)^{l_1}
\phi_{j_1}(s_2)ds_2ds_1ds.
\end{equation}

\vspace{3mm}

If we prove w.~p.~1 the formula

\vspace{-1mm}
\begin{equation}
\label{ogo2000}
\hbox{\vtop{\offinterlineskip\halign{
\hfil#\hfil\cr
{\rm l.i.m.}\cr
$\stackrel{}{{}_{p_1, p_3\to \infty}}$\cr
}} }
\sum\limits_{j_1=0}^{p_1}\sum\limits_{j_3=0}^{p_3}
C_{j_3 j_3 j_1}\zeta_{j_1}^{(i_1)}=
\frac{1}{2}\int\limits_t^T(t-s)^{2l}
\int\limits_t^s(t-s_1)^{l_1}d{\bf f}_{s_1}^{(i_1)}ds,
\end{equation}

\vspace{2mm}
\noindent
where  
coefficients $C_{j_3 j_3 j_1}$ has the form (\ref{ogo1991}), 
then using Theorems 1, 2 
and standard relations between iterated Ito and 
Stratonovich stochastic integrals we obtain the expansion (\ref{ogo101ee}).

Using Theorems 1, 2 and the Ito formula we 
can write 

\vspace{-1mm}
$$
\frac{1}{2}\int\limits_t^T(t-s)^{2l}
\int\limits_t^s(t-s_1)^{l_1}d{\bf f}_{s_1}^{(i_1)}ds=
\frac{1}{2}\int\limits_t^T(t-s_1)^{l_1}
\int\limits_{s_1}^T(t-s)^{2l}dsd{\bf f}_{s_1}^{(i_1)}=
$$

$$
=
\frac{1}{2}\sum\limits_{j_1=0}^{2l+l_1+1}
\tilde C_{j_1}\zeta_{j_1}^{(i_1)}\ \ \ \hbox{w.\ p.\ 1},
$$

\vspace{3mm}
\noindent
where
$$
\tilde C_{j_1}=
\int\limits_t^T
\phi_{j_1}(s_1)(t-s_1)^{l_1}\int\limits_{s_1}^T(t-s)^{2l}dsds_1.
$$

\vspace{3mm}

Then
$$
\sum\limits_{j_1=0}^{p_1}\sum\limits_{j_3=0}^{p_3}
C_{j_3 j_3 j_1}\zeta_{j_1}^{(i_1)}-
\frac{1}{2}\sum\limits_{j_1=0}^{2l+l_1+1}
\tilde C_{j_1}\zeta_{j_1}^{(i_1)}=
$$

\vspace{2mm}
$$
=
\sum\limits_{j_1=0}^{2l+l_1+1}
\left(\sum\limits_{j_3=0}^{p_3}
C_{j_3 j_3 j_1}-\frac{1}{2}\tilde C_{j_1}\right)
\zeta_{j_1}^{(i_1)}+
\sum\limits_{j_1=2l+l_1+2}^{p_1}
\sum\limits_{j_3=0}^{p_3}
C_{j_3 j_3 j_1}\zeta_{j_1}^{(i_1)}.
$$

\vspace{6mm}

Therefore,

\vspace{-2mm}
$$
\hbox{\vtop{\offinterlineskip\halign{
\hfil#\hfil\cr
{\rm lim}\cr
$\stackrel{}{{}_{p_1,p_3\to \infty}}$\cr
}} }
{\sf M}\left\{\left(
\sum\limits_{j_1=0}^{p_1}\sum\limits_{j_3=0}^{p_3}C_{j_3j_3 j_1}
\zeta_{j_1}^{(i_1)}-
\frac{1}{2}\int\limits_t^T(t-s)^{2l}
\int\limits_t^s(t-s_1)^{l_1}d{\bf f}_{s_1}^{(i_1)}ds\right)^2\right\}=
$$

\vspace{2mm}
\begin{equation}
\label{ogo2100}
=\hbox{\vtop{\offinterlineskip\halign{
\hfil#\hfil\cr
{\rm lim}\cr
$\stackrel{}{{}_{p_3\to \infty}}$\cr
}} }\sum\limits_{j_1=0}^{2l+l_1+1}
\left(\sum\limits_{j_3=0}^{p_3}C_{j_3j_3 j_1}-
\frac{1}{2}\tilde C_{j_1}\right)^2+
\hbox{\vtop{\offinterlineskip\halign{
\hfil#\hfil\cr
{\rm lim}\cr
$\stackrel{}{{}_{p_1,p_3\to \infty}}$\cr
}} }{\sf M}\left\{\left(
\sum\limits_{j_1=2l+l_1+2}^{p_1}
\sum\limits_{j_3=0}^{p_3}
C_{j_3 j_3 j_1}\zeta_{j_1}^{(i_1)}\right)^2\right\}.
\end{equation}

\vspace{7mm}

Let us prove that

\vspace{-3mm}
\begin{equation}
\label{ogo2110}
\hbox{\vtop{\offinterlineskip\halign{
\hfil#\hfil\cr
{\rm lim}\cr
$\stackrel{}{{}_{p_3\to \infty}}$\cr
}} }
\left(\sum\limits_{j_3=0}^{p_3}C_{j_3j_3 j_1}-
\frac{1}{2}\tilde C_{j_1}\right)^2=0.
\end{equation}

\vspace{2mm}

We have

\vspace{-2mm}
$$
\left(\sum\limits_{j_3=0}^{p_3}C_{j_3j_3 j_1}-
\frac{1}{2}\tilde C_{j_1}\right)^2=
$$

$$
=
\left(\sum\limits_{j_3=0}^{p_3}
\int\limits_t^T\phi_{j_1}(s_2)(t-s_2)^{l_1}ds_2
\int\limits_{s_2}^T\phi_{j_3}(s_1)(t-s_1)^{l}ds_1
\int\limits_{s_1}^T\phi_{j_3}(s)(t-s)^{l}ds-\right.
$$

$$
-\left.
\frac{1}{2}
\int\limits_t^T
\phi_{j_1}(s_1)(t-s_1)^{l_1}\int\limits_{s_1}^T(t-s)^{2l}dsds_1\right)^2=
$$

$$
=\left(\frac{1}{2}\sum\limits_{j_3=0}^{p_3}
\int\limits_t^T\phi_{j_1}(s_2)(t-s_2)^{l_1}
\left(\int\limits_{s_2}^T\phi_{j_3}(s_1)(t-s_1)^{l}ds_1\right)^2ds_2-
\right.
$$

$$
\left.-
\frac{1}{2}
\int\limits_t^T
\phi_{j_1}(s_1)(t-s_1)^{l_1}\int\limits_{s_1}^T(t-s)^{2l}dsds_1\right)^2=
$$

$$
=\frac{1}{4}\left(
\int\limits_t^T\phi_{j_1}(s_1)(t-s_1)^{l_1}\left(
\sum\limits_{j_3=0}^{p_3}
\left(\int\limits_{s_1}^T\phi_{j_3}(s)(t-s)^{l}ds\right)^2
-\int\limits_{s_1}^T(t-s)^{2l}ds\right)ds_1\right)^2=
$$

$$
=\frac{1}{4}\left(
\int\limits_t^T\phi_{j_1}(s_1)(t-s_1)^{l_1}\left(
\int\limits_{s_1}^T(t-s)^{2l}ds-
\sum\limits_{j_3=p_3+1}^{\infty}
\left(\int\limits_{s_1}^T\phi_{j_3}(s)(t-s)^{l}ds\right)^2
-\int\limits_{s_1}^T(t-s)^{2l}ds\right)ds_1\right)^2
$$

\begin{equation}
\label{ogo3000}
=\frac{1}{4}\left(
\int\limits_t^T\phi_{j_1}(s_1)(t-s_1)^{l_1}
\sum\limits_{j_3=p_3+1}^{\infty}
\left(\int\limits_{s_1}^T\phi_{j_3}(s)(t-s)^{l}ds\right)^2
ds_1\right)^2.
\end{equation}

\vspace{6mm}

In order to 
get (\ref{ogo3000}) we used the Parseval equality, which 
looks as follows

\vspace{-2mm}
\begin{equation}
\label{ogo3010}
\sum_{j_3=0}^{\infty}\left(\int\limits_{s_1}^T\phi_{j_3}(s)
(t-s)^lds\right)^2=
\int\limits_t^T K^2(s,s_1)ds,
\end{equation}

\vspace{3mm}
\noindent
where
$$
K(s,s_1)=(t-s)^l
{\bf 1}_{\{s_1<s\}},\ \ \ s, s_1\in [t, T].
$$

\vspace{4mm}

Taking into account nondecreasing
of the functional sequence

\vspace{-1mm}
$$
u_n(s_1)=\sum_{j_3=0}^{n}\left(\int\limits_{s_1}^T\phi_{j_3}(s)
(t-s)^lds\right)^2,
$$

\vspace{3mm}
\noindent
continuity of its members and continuity of the limit
function 

\vspace{-1mm}
$$
u(s_1)=\int\limits_{s_1}^T(t-s)^{2l}ds
$$ 

\vspace{2mm}
\noindent
at the interval $[t, T]$ according to the Dini Theorem we have 
uniform 
convergence
of the functional sequence $u_n(s_1)$ to the limit 
function $u(s_1)$ at the interval $[t, T]$.

From (\ref{ogo3000}) using the inequality of Cauchy--Bunyakovsky, we obtain

$$
\left(\sum\limits_{j_3=0}^{p_3}C_{j_3j_3 j_1}-
\frac{1}{2}\tilde C_{j_1}\right)^2\le
$$

$$
\le
\frac{1}{4}
\int\limits_t^T\phi_{j_1}^2(s_1)(t-s_1)^{2l_1}ds_1
\int\limits_t^T\left(\sum\limits_{j_3=p_3+1}^{\infty}
\left(\int\limits_{s_1}^T\phi_{j_3}(s)(t-s)^{l}ds\right)^2
\right)^2 ds_1\le
$$

\begin{equation}
\label{ogo3020}
\le\frac{1}{4}\varepsilon^2 (T-t)^{2l_1}\int\limits_t^T\phi_{j_1}^2(s_1)ds_1
(T-t)=\frac{1}{4}(T-t)^{2l_1+1}\varepsilon^2
\end{equation}

\vspace{3mm}
\noindent 
when $p_3>N(\varepsilon),$ where $N(\varepsilon)$
exists for any $\varepsilon>0.$

From (\ref{ogo3020}) it follows 
(\ref{ogo2110}).

We have

\vspace{-2mm}
\begin{equation}
\label{ogo3030}
\sum\limits_{j_3=0}^{p_3}
\sum\limits_{j_1=2l+l_1+2}^{p_1}
C_{j_3 j_3 j_1}\zeta_{j_1}^{(i_1)}
=
\sum\limits_{j_3=0}^{p_3}
\sum\limits_{j_1=2l+l_1+2}^{2(j_3+l+1)+l_1}
C_{j_3 j_3 j_1}\zeta_{j_1}^{(i_1)}.
\end{equation}

\vspace{4mm}

We put  $2(j_3+l+1)+l_1$ instead of $p_1$, since
$C_{j_3j_3j_1}=0$ when $j_1>2(j_3+l+1)+l_1.$ 
This conclusion follows from the relation

\vspace{-2mm}
$$
C_{j_3j_3j_1}=
\frac{1}{2}
\int\limits_t^T\phi_{j_1}(s_2)(t-s_2)^{l_1}
\left(
\int\limits_{s_2}^T\phi_{j_3}(s_1)(t-s_1)^{l}ds_1\right)^2ds_2=
$$

$$
=
\frac{1}{2}\int\limits_t^T\phi_{j_1}(s_2)Q_{2(j_3+l+1)+l_1}(s_2)ds_2,
$$

\vspace{3mm}
\noindent
where $Q_{2(j_3+l+1)+l_1}(s)$ is a polynomial of degree
$2(j_3+l+1)+l_1.$

It is easy to see that

\vspace{-2mm}
\begin{equation}
\label{ogo3040}
\sum\limits_{j_3=0}^{p_3}
\sum\limits_{j_1=2l+l_1+2}^{2(j_3+l+1)+l_1}
C_{j_3 j_3 j_1}\zeta_{j_1}^{(i_1)}=
\sum\limits_{j_1=2l+l_1+2}^{2(p_3+l+1)+l_1}
\sum\limits_{j_3=0}^{p_3}
C_{j_3 j_3 j_1}\zeta_{j_1}^{(i_1)}.
\end{equation}

\vspace{3mm}

Note that we 
included
some zero coefficients $C_{j_3 j_3 j_1}$ 
into the sum 
$\sum\limits_{j_3=0}^{p_3}$. 

From (\ref{ogo3030}) and (\ref{ogo3040}) we have

$$
{\sf M}\left\{\left(\sum\limits_{j_3=0}^{p_3}
\sum\limits_{j_1=2l+l_1+2}^{p_1}
C_{j_3 j_3 j_1}\zeta_{j_1}^{(i_1)}\right)^2\right\}=
$$

\vspace{2mm}
$$
=
{\sf M}\left\{\left(
\sum\limits_{j_1=2l+l_1+2}^{2(p_3+l+1)+l_1}
\sum\limits_{j_3=0}^{p_3}
C_{j_3 j_3 j_1}\zeta_{j_1}^{(i_1)}\right)^2\right\}
=\sum\limits_{j_1=2l+l_1+2}^{2(p_3+l+1)+l_1}
\left(\sum\limits_{j_3=0}^{p_3}
C_{j_3 j_3 j_1}\right)^2=
$$

\vspace{2mm}
$$
=\sum\limits_{j_1=2l+l_1+2}^{2(p_3+l+1)+l_1}
\left(\frac{1}{2}\sum\limits_{j_3=0}^{p_3}
\int\limits_t^T\phi_{j_1}(s_2)(t-s_2)^{l_1}
\left(\int\limits_{s_2}^T\phi_{j_3}(s_1)(t-s_1)^{l}ds_1\right)^2ds_2\right)^2=
$$

\vspace{2mm}
$$
=\frac{1}{4}\sum\limits_{j_1=2l+l_1+2}^{2(p_3+l+1)+l_1}
\left(
\int\limits_t^T\phi_{j_1}(s_2)(t-s_2)^{l_1}
\sum\limits_{j_3=0}^{p_3}
\left(\int\limits_{s_2}^T\phi_{j_3}(s_1)(t-s_1)^{l}ds_1\right)^2
ds_2\right)^2=
$$

\vspace{2mm}
$$
=\frac{1}{4}\sum\limits_{j_1=2l+l_1+2}^{2(p_3+l+1)+l_1}\left(
\int\limits_t^T\phi_{j_1}(s_2)(t-s_2)^{l_1}\left(
\int\limits_{s_2}^T(t-s_1)^{2l}ds_1-
\sum\limits_{j_3=p_3+1}^{\infty}
\left(\int\limits_{s_2}^T\phi_{j_3}(s_1)(t-s_1)^{l}ds_1\right)^2
\right)ds_2\right)^2
$$

\vspace{2mm}
\begin{equation}
\label{ogo3100}
=\frac{1}{4}\sum\limits_{j_1=2l+l_1+2}^{2(p_3+l+1)+l_1}\left(
\int\limits_t^T\phi_{j_1}(s_2)(t-s_2)^{l_1}
\sum\limits_{j_3=p_3+1}^{\infty}
\left(\int\limits_{s_2}^T\phi_{j_3}(s_1)(t-s_1)^{l}ds_1\right)^2
ds_2\right)^2.
\end{equation}

\vspace{6mm}

In order to 
get
(\ref{ogo3100}) we used the Parseval equality 
(\ref{ogo3010}) and the following relation

\vspace{-2mm} 
$$
\int\limits_t^T\phi_{j_1}(s)Q_{2l+1+l_1}(s)ds=0,\ \ \ j_1>2l+1+l_1,
$$

\vspace{3mm}
\noindent
where $Q_{2l+1+l_1}(s)$ is a polynomial of degree
$2l+1+l_1.$

Further, we have for $j_3\in\mathbb{N}$

$$
\left(\int\limits_{s_2}^T\phi_{j_3}(s_1)(t-s_1)^lds_1\right)^2=
\frac{(T-t)^{2l+1}(2j_3+1)}{2^{2l+2}}
\left(\int\limits_{z(s_2)}^1
P_{j_3}(y)(1+y)^ldy\right)^2=
$$

\vspace{2mm}
$$
=\frac{(T-t)^{2l+1}}{2^{2l+2}(2j_3+1)}\left(
\left(1+z(s_2)\right)^l
Q_{j_3}(s_2)-
l\int\limits_{z(s_2)}^1
\left(P_{j_3+1}(y)-P_{j_3-1}(y)\right)\left(1+y\right)^{l-1}dy\right)^2\le
$$

\vspace{2mm}
$$
\le\frac{(T-t)^{2l+1}2}{2^{2l+2}(2j_3+1)}\left(
\left(\frac{2(s_2-t)}{T-t}\right)^{2l}
Q_{j_3}^2(s_2))+
l^2
\left(
\int\limits_{z(s_2)}^1
\left(P_{j_3+1}(y)-P_{j_3-1}(y)\right)\left(1+y\right)^{l-1}dy\right)^2
\right)\le
$$

\vspace{2mm}
$$
\le\frac{(T-t)^{2l+1}}{2^{2l+1}(2j_3+1)}\left(
2^{2l+1}
H_{j_3}(s_2)+
l^2
\int\limits_{z(s_2)}^1
(1+y)^{2l-2}dy
\int\limits_{z(s_2)}^1
\left(P_{j_3+1}(y)-P_{j_3-1}(y)\right)^2dy
\right)\le
$$

\vspace{2mm}
$$
\le\frac{(T-t)^{2l+1}}{2^{2l+1}(2j_3+1)}\left(
2^{2l+1}
H_{j_3}(s_2)
+\frac{2^{2l}l^2}{2l-1}\left(1-\left(\frac{(s_2-t)}{T-t}\right)^{2l-1}\right)
\int\limits_{z(s_2)}^1
\left(P_{j_3+1}^2(y)+P_{j_3-1}^2(y)\right)dy
\right)\le
$$

\vspace{2mm}
\begin{equation}
\label{ogo4000}
\le\frac{(T-t)^{2l+1}}{2(2j_3+1)}\Biggl(
2
H_{j_3}(s_2)
\Biggl.
+\frac{l^2}{2l-1}
\int\limits_{z(s_2)}^1
\left(P_{j_3+1}^2(y)+P_{j_3-1}^2(y)\right)dy
\Biggr),
\end{equation}

\vspace{6mm}
\noindent
where
$$
Q_{j_3}(s_2)=P_{j_3-1}(z(s_2))-P_{j_3+1}(z(s_2)),
$$

\vspace{-2mm}
$$
H_{j_3}(s_2)=P_{j_3-1}^2(z(s_2))+P_{j_3+1}^2(z(s_2)).
$$

\vspace{7mm}

Let us estimate the right-hand side
of (\ref{ogo4000}) using
(\ref{ogo23})

\vspace{1mm}
$$
\left(\int\limits_{s_2}^T\phi_{j_3}(s_1)(t-s_1)^lds_1\right)^2 <
$$

\vspace{1mm}
$$
<
\frac{(T-t)^{2l+1}}{2(2j_3+1)}\left(\frac{K^2}{j_3+2}+\frac{K^2}{j_3}
\right)\left(\frac{2}
{(1-\left(z(s_2))^2
\right)^{1/2}}
+\frac{l^2}{2l-1}
\int\limits_{z(s_2)}^1
\frac{dy}{\left(1-y^2\right)^{1/2}}\right)<
$$

\vspace{2mm}
\begin{equation}
\label{ogo4010}
<\frac{(T-t)^{2l+1}K^2}{2j_3^2}\left(
\frac{2}
{(1-\left(z(s_2))^2
\right)^{1/2}}+
\frac{l^2\pi}{2l-1}\right),\ \ \ s\in(t, T).
\end{equation}

\vspace{5mm}

From (\ref{ogo3100}) and (\ref{ogo4010}) we obtain

\vspace{1mm}
$$
{\rm M}\left\{\left(\sum\limits_{j_3=0}^{p_3}
\sum\limits_{j_1=2l+l_1+2}^{p_1}
C_{j_3 j_3 j_1}\zeta_{j_1}^{(i_1)}\right)^2\right\}\le
$$

\vspace{2mm}
$$
\le
\frac{1}{4}\sum\limits_{j_1=2l+l_1+2}^{2(p_3+l+1)+l_1}\left(
\int\limits_t^T|\phi_{j_1}(s_2)|(t-s_2)^{l_1}
\sum\limits_{j_3=p_3+1}^{\infty}
\left(\int\limits_{s_2}^T\phi_{j_3}(s_1)(t-s_1)^{l}ds_1\right)^2
ds_2\right)^2\le
$$

\vspace{2mm}
$$
\le
\frac{1}{4}(T-t)^{2l_1}\sum\limits_{j_1=2l+l_1+2}^{2(p_3+l+1)+l_1}\left(
\int\limits_t^T|\phi_{j_1}(s_2)|
\sum\limits_{j_3=p_3+1}^{\infty}
\left(\int\limits_{s_2}^T\phi_{j_3}(s_1)(t-s_1)^{l}ds_1\right)^2
ds_2\right)^2<
$$

\vspace{6mm}
$$
<
\frac{(T-t)^{4l+2l_1+1}K^4 K_1^2}{16}\times
$$

\vspace{1mm}
$$
\times
\sum\limits_{j_1=2l+l_1+2}^{2(p_3+l+1)+l_1}
\left(\left(\int\limits_t^T
\frac{2ds_2}
{(1-\left(z(s_2))^2
\right)^{3/4}}\right.\right.
\left.\left.
+\frac{l^2\pi}{2l-1}
\int\limits_t^T
\frac{ds_2}
{(1-\left(z(s_2))^2
\right)^{1/4}}\right)
\sum\limits_{j_3=p_3+1}^{\infty}\frac{1}{j_3^2}
\right)^2\le 
$$

\vspace{2mm}
$$
\le
\frac{(T-t)^{4l+2l_1+3}K^4 K_1^2}{64}\cdot\frac{2p_3+1}{p_3^2}
\left(\int\limits_{-1}^1
\frac{2dy}
{(1-y^2)^{3/4}}
+\frac{l^2\pi}{2l-1}
\int\limits_{-1}^1
\frac{dy}
{(1-y^2)^{1/4}}\right)^2\le
$$

\vspace{2mm}
\begin{equation}
\label{ogo5000}
\le (T-t)^{4l+2l_1+3}C\frac{2p_3+1}{p_3^2}\to 0\ \ \ 
\hbox{when}\ \ p_3\ \to\ \infty,
\end{equation}

\vspace{7mm}
\noindent
where constant $C$ does not depend on $p_3$ and $T-t$.

From (\ref{ogo2100}), (\ref{ogo2110}), and (\ref{ogo5000}) it follows 
(\ref{ogo2000}), and the relation (\ref{ogo2000}) implies the expansion
(\ref{ogo101ee}).

Let us consider Case 4 ($l_1=l_2=l_3=l=0, 1, 2,\ldots$ and 
$i_1, i_2, i_3=1,\ldots,m$). So, we will prove the following expansion
for iterated Stratonovich stochastic integral of third multiplicity

\vspace{1mm}
\begin{equation}
\label{ogo10100}
I_{{l l l}_{T,t}}^{*(i_1i_2i_3)}=
\hbox{\vtop{\offinterlineskip\halign{
\hfil#\hfil\cr
{\rm l.i.m.}\cr
$\stackrel{}{{}_{p_1,p_2,p_3\to \infty}}$\cr
}} }\sum_{j_1=0}^{p_1}\sum_{j_2=0}^{p_2}\sum_{j_3=0}^{p_3}
C_{j_3 j_2 j_1}\zeta_{j_1}^{(i_1)}\zeta_{j_2}^{(i_2)}\zeta_{j_3}^{(i_3)}\ \ \
(i_1, i_2, i_3=1,\ldots,m),
\end{equation}

\vspace{4mm}
\noindent
where the series converges in the mean-square sense,\ \ 
$l=0, 1, 2,\ldots$, and

\begin{equation}
\label{ogo19900}
C_{j_3 j_2 j_1}=\int\limits_t^T
\phi_{j_3}(s)(t-s)^{l}\int\limits_t^s(t-s_1)^{l}
\phi_{j_2}(s_1)
\int\limits_t^{s_1}(t-s_2)^{l}
\phi_{j_1}(s_2)ds_2ds_1ds.
\end{equation}

\vspace{3mm}

If we prove w.~p.~1 the following formula

\vspace{-1mm}
\begin{equation}
\label{ogo20000}
\hbox{\vtop{\offinterlineskip\halign{
\hfil#\hfil\cr
{\rm l.i.m.}\cr
$\stackrel{}{{}_{p_1,p_3\to \infty}}$\cr
}} }\sum_{j_1=0}^{p_1}\sum_{j_3=0}^{p_3}
C_{j_1 j_3 j_1}\zeta_{j_3}^{(i_2)}=0,
\end{equation}

\vspace{4mm}
\noindent
where 
coefficients $C_{j_3 j_2 j_1}$ 
have
the form (\ref{ogo19900}), then
using Theorems 1, 2, relations (\ref{ogo200}), (\ref{ogo2000})
when $l_1=l_3=l$
and standard relations
between 
iterated Ito and Stratonovich stochastic integrals
we will have the expansion (\ref{ogo10100}).

Since 
$\psi_1(s),$ $\psi_2(s),$ $\psi_3(s)\equiv (t-s)^l$, then
the following equality for the Fourier coefficients takes place

\vspace{-1mm}
$$
C_{j_1 j_1 j_3}+C_{j_1 j_3 j_1}+C_{j_3 j_1 j_1}=\frac{1}{2}
C_{j_1}^2 C_{j_3},
$$ 

\vspace{4mm}
\noindent
where $C_{j_3 j_2 j_1}$ has the form (\ref{ogo19900}) and

\vspace{-2mm}
$$
C_{j_1}=\int\limits_t^T
\phi_{j_1}(s)(t-s)^{l}ds.
$$

\vspace{3mm}

Then w.~p.~1

\begin{equation}
\label{sodom310}
\hbox{\vtop{\offinterlineskip\halign{
\hfil#\hfil\cr
{\rm l.i.m.}\cr
$\stackrel{}{{}_{p_1,p_3\to \infty}}$\cr
}} }\sum_{j_1=0}^{p_1}\sum_{j_3=0}^{p_3}
C_{j_1 j_3 j_1}\zeta_{j_3}^{(i_2)}=
\hbox{\vtop{\offinterlineskip\halign{
\hfil#\hfil\cr
{\rm l.i.m.}\cr
$\stackrel{}{{}_{p_1,p_3\to \infty}}$\cr
}} }\sum_{j_1=0}^{p_1}\sum_{j_3=0}^{p_3}
\left(\frac{1}{2}C_{j_1}^2 C_{j_3}-C_{j_1 j_1 j_3}-C_{j_3 j_1 j_1}
\right)\zeta_{j_3}^{(i_2)}.
\end{equation}

\vspace{5mm}
 
Taking into account (\ref{ogo200}) and (\ref{ogo2000}) 
when $l_3=l_1=l$ and the Ito formula,
we have w.~p.~1

\vspace{1mm}
$$
\hbox{\vtop{\offinterlineskip\halign{
\hfil#\hfil\cr
{\rm l.i.m.}\cr
$\stackrel{}{{}_{p_1,p_3\to \infty}}$\cr
}} }\sum_{j_1=0}^{p_1}\sum_{j_3=0}^{p_3}
C_{j_1 j_3 j_1}\zeta_{j_3}^{(i_2)}=
$$

\vspace{1mm}
$$
=
\frac{1}{2}\sum\limits_{j_1=0}^{l}
C_{j_1}^2\sum\limits_{j_3=0}^{l}C_{j_3}\zeta_{j_3}^{(i_2)}-
\hbox{\vtop{\offinterlineskip\halign{
\hfil#\hfil\cr
{\rm l.i.m.}\cr
$\stackrel{}{{}_{p_1,p_3\to \infty}}$\cr
}} }\sum_{j_1=0}^{p_1}\sum_{j_3=0}^{p_3}
C_{j_1 j_1 j_3}\zeta_{j_3}^{(i_2)}-
\hbox{\vtop{\offinterlineskip\halign{
\hfil#\hfil\cr
{\rm l.i.m.}\cr
$\stackrel{}{{}_{p_1,p_3\to \infty}}$\cr
}} }\sum_{j_1=0}^{p_1}\sum_{j_3=0}^{p_3}
C_{j_3 j_1 j_1}\zeta_{j_3}^{(i_2)}
=
$$

\vspace{2mm}
$$
=\frac{1}{2}\sum\limits_{j_1=0}^{l}
C_{j_1}^2\int\limits_t^T(t-s)^ld{\bf f}_s^{(i_2)}
-\frac{1}{2}\int\limits_t^T(t-s)^{l}
\int\limits_t^s(t-s_1)^{2l}ds_1d{\bf f}_s^{(i_2)}-
$$

$$
-\frac{1}{2}\int\limits_t^T(t-s)^{2l}
\int\limits_t^s(t-s_1)^{l}d{\bf f}_{s_1}^{(i_2)}ds=
$$

\vspace{2mm}
$$
=\frac{1}{2}\sum\limits_{j_1=0}^{l}
C_{j_1}^2\int\limits_t^T(t-s)^ld{\bf f}_s^{(i_2)}
+\frac{1}{2(2l+1)}\int\limits_t^T(t-s)^{3l+1}d{\bf f}_s^{(i_2)}-
$$

$$
-\frac{1}{2}\int\limits_t^T(t-s_1)^{l}
\int\limits_{s_1}^T(t-s)^{2l}dsd{\bf f}_{s_1}^{(i_2)}=
$$

\vspace{2mm}
$$
=\frac{1}{2}\sum\limits_{j_1=0}^{l}
C_{j_1}^2\int\limits_t^T(t-s)^ld{\bf f}_s^{(i_2)}
+\frac{1}{2(2l+1)}\int\limits_t^T(t-s)^{3l+1}d{\bf f}_s^{(i_2)}-
$$

$$
-\frac{1}{2(2l+1)}\left((T-t)^{2l+1}\int\limits_t^T(t-s)^{l}
d{\bf f}_{s}^{(i_2)}+\int\limits_t^T(t-s)^{3l+1}
d{\bf f}_{s}^{(i_2)}\right)=
$$

\vspace{4mm}
$$
=\frac{1}{2}\sum\limits_{j_1=0}^{l}
C_{j_1}^2\int\limits_t^T(t-s)^ld{\bf f}_s^{(i_2)}-
\frac{(T-t)^{2l+1}}{2(2l+1)}\int\limits_t^T(t-s)^{l}
d{\bf f}_{s}^{(i_2)}=
$$

\vspace{2mm}
$$
=
\frac{1}{2}\left(\sum\limits_{j_1=0}^{l}
C_{j_1}^2-\int\limits_t^T(t-s)^{2l}ds\right)
\int\limits_t^T(t-s)^ld{\bf f}_s^{(i_2)}=0.
$$

\vspace{7mm}

Here, the Parseval equality looks as follows

\vspace{-1mm}
$$
\sum\limits_{j_1=0}^{\infty}
C_{j_1}^2=
\sum\limits_{j_1=0}^{l}
C_{j_1}^2=\int\limits_t^T(t-s)^{2l}ds=\frac{(T-t)^{2l+1}}{2l+1}
$$

\vspace{4mm}
\noindent
and      
$$
\int\limits_t^{T}(t-s)^{l}
d{\bf f}_{s}^{(i_2)}=
\sum\limits_{j_3=0}^{l}C_{j_3}\zeta_{j_3}^{(i_2)}\ \ \ \hbox{w.\ p.\ 1}.
$$

\vspace{4mm}

The expansion (\ref{ogo10100}) is proved. Theorem 4 is proved.

It is easy to see that using the Ito formula,
we obtain for the case $i_1=i_2=i_3$ 

\vspace{1mm}
$$
\int\limits_t^{*T}(t-s)^{l}\int\limits_t^{*s}
(t-s_1)^l\int\limits_t^{*s_1}(t-s_2)^{l}
d{\bf f}_{s_2}^{(i_1)}d{\bf f}_{s_1}^{(i_1)}d{\bf f}_{s}^{(i_1)}=
$$

\vspace{1mm}
$$
=\frac{1}{6}\Biggl(
\int\limits_t^{T}(t-s)^{l}
d{\bf f}_{s}^{(i_1)}\Biggr)^3
=\frac{1}{6}\left(
\sum\limits_{j_1=0}^{l}C_{j_1}\zeta_{j_1}^{(i_1)}\right)^3=
$$

\vspace{1mm}
\begin{equation}
\label{sodom400}
=
\sum\limits_{j_1, j_2, j_3=0}^{l}
C_{j_3 j_2 j_1}\zeta_{j_1}^{(i_1)}\zeta_{j_2}^{(i_1)}\zeta_{j_3}^{(i_1)}\ \ \
\hbox{w.\ p.\ 1}.
\end{equation}

\vspace{5mm}

\section{Expansion of Iterated Stratonovich Stochastic Integrals of 
Multiplicity 3. The Case of Trigonometric Functions}

\vspace{5mm}

In this section we will prove the following theorem.

\vspace{2mm}

{\bf Theorem 5}\ \cite{12}-\cite{16}, \cite{19}, \cite{20}-\cite{20xx2}.\
{\it Suppose that
$\{\phi_j(x)\}_{j=0}^{\infty}$ is a complete orthonormal
system of trigonometric functions
in the space $L_2([t, T])$.
Then, for the iterated Stratonovich stochastic integral 
of third multiplicity

$$
{\int\limits_t^{*}}^T
{\int\limits_t^{*}}^{t_3}
{\int\limits_t^{*}}^{t_2}
d{\bf f}_{t_1}^{(i_1)}
d{\bf f}_{t_2}^{(i_2)}d{\bf f}_{t_3}^{(i_3)}\ \ \ (i_1, i_2, i_3=1,\ldots,m)
$$

\vspace{4mm}
\noindent
the following 
expansion

\begin{equation}
\label{feto19001ee}
{\int\limits_t^{*}}^T
{\int\limits_t^{*}}^{t_3}
{\int\limits_t^{*}}^{t_2}
d{\bf f}_{t_1}^{(i_1)}
d{\bf f}_{t_2}^{(i_2)}d{\bf f}_{t_3}^{(i_3)}\ 
=
\hbox{\vtop{\offinterlineskip\halign{
\hfil#\hfil\cr
{\rm l.i.m.}\cr
$\stackrel{}{{}_{p_1,p_2,p_3\to \infty}}$\cr
}} }\sum_{j_1=0}^{p_1}\sum_{j_2=0}^{p_2}\sum_{j_3=0}^{p_3}
C_{j_3 j_2 j_1}\zeta_{j_1}^{(i_1)}\zeta_{j_2}^{(i_2)}\zeta_{j_3}^{(i_3)}
\end{equation}

\vspace{5mm}
\noindent
converging in the mean-square sense is valid, where

\vspace{-2mm}
$$
C_{j_3 j_2 j_1}=\int\limits_t^T
\phi_{j_3}(s)\int\limits_t^s
\phi_{j_2}(s_1)
\int\limits_t^{s_1}
\phi_{j_1}(s_2)ds_2ds_1ds.
$$
}

\vspace{4mm}

{\bf Proof.}\ If we prove w.~p.~1 the following formulas

\vspace{-1mm}
\begin{equation}
\label{ogo1299}
\hbox{\vtop{\offinterlineskip\halign{
\hfil#\hfil\cr
{\rm l.i.m.}\cr
$\stackrel{}{{}_{p_1, p_3\to \infty}}$\cr
}} }
\sum\limits_{j_1=0}^{p_1}\sum\limits_{j_3=0}^{p_3}
C_{j_3 j_1 j_1}\zeta_{j_3}^{(i_3)}
=
\frac{1}{2}\int\limits_t^T\int\limits_t^{\tau}dsd{\bf f}_{\tau}^{(i_3)},
\end{equation}

\vspace{3mm}
\begin{equation}
\label{ogo1399}
\hbox{\vtop{\offinterlineskip\halign{
\hfil#\hfil\cr
{\rm l.i.m.}\cr
$\stackrel{}{{}_{p_1, p_3\to \infty}}$\cr
}} }
\sum\limits_{j_1=0}^{p_1}\sum\limits_{j_3=0}^{p_3}
C_{j_3 j_3 j_1}\zeta_{j_1}^{(i_1)}
=
\frac{1}{2}\int\limits_t^T\int\limits_t^{\tau}d{\bf f}_{s}^{(i_1)}d\tau,
\end{equation}

\vspace{3mm}
\begin{equation}
\label{ogo13a99}
\hbox{\vtop{\offinterlineskip\halign{
\hfil#\hfil\cr
{\rm l.i.m.}\cr
$\stackrel{}{{}_{p_1, p_3\to \infty}}$\cr
}} }
\sum\limits_{j_1=0}^{p_1}\sum\limits_{j_3=0}^{p_3}
C_{j_1 j_3 j_1}\zeta_{j_3}^{(i_2)}
=0,
\end{equation}

\vspace{5mm}
\noindent
then from the equalities (\ref{ogo1299})--(\ref{ogo13a99}),
Theorems 1, 2, and
standard relations 
between 
iterated Ito and
Stratonovich stochastic integrals we will obtain
the expansion (\ref{feto19001ee}).

We have

\vspace{-1mm}
$$
S_{p_1,p_3}\stackrel{\sf def}{=}\sum\limits_{j_3=0}^{p_3}\sum\limits_{j_1=0}^{p_1}
C_{j_3 j_1 j_1}\zeta_{j_3}^{(i_3)}=
\frac{(T-t)^{3/2}}{6}\zeta_{0}^{(i_3)}+
$$

\vspace{3mm}
$$
+\sum\limits_{j_1=1}^{p_1}C_{0,2j_1,2j_1}
\zeta_0^{(i_3)}+\sum\limits_{j_1=1}^{p_1}C_{0,2j_1-1,2j_1-1}\zeta_0^{(i_3)}+
\sum\limits_{j_3=1}^{p_1}C_{2j_3,0,0}\zeta_{2j_3}^{(i_3)}+
$$

\vspace{3mm}
$$
+
\sum\limits_{j_3=1}^{p_3}\sum_{j_1=1}^{p_1}
C_{2j_3,2j_1,2j_1}\zeta_{2j_3}^{(i_3)}+
\sum\limits_{j_3=1}^{p_3}\sum_{j_1=1}^{p_1}
C_{2j_3,2j_1-1,2j_1-1}\zeta_{2j_3}^{(i_3)}+
\sum\limits_{j_3=1}^{p_3}
C_{2j_3-1,0,0}\zeta_{2j_3-1}^{(i_3)}+
$$

\vspace{3mm}
\begin{equation}
\label{ogo900}
+\sum\limits_{j_3=1}^{p_3}\sum_{j_1=1}^{p_1}
C_{2j_3-1,2j_1,2j_1}\zeta_{2j_3-1}^{(i_3)}+
\sum\limits_{j_3=1}^{p_3}\sum_{j_1=1}^{p_1}
C_{2j_3-1,2j_1-1,2j_1-1}\zeta_{2j_3-1}^{(i_3)},
\end{equation}

\vspace{7mm}
\noindent
where the summation is stopped, when $2j_1,$ $2j_1-1> p_1$
or $2j_3,$ $2j_3-1> p_3$ and

\vspace{1mm}
\begin{equation}
\label{ogo901}
C_{0,2l,2l}=\frac{(T-t)^{3/2}}{8\pi^2l^2},\ \ \
C_{0,2l-1,2l-1}=\frac{3(T-t)^{3/2}}{8\pi^2l^2},\ \ \
C_{2l,0,0}=\frac{\sqrt{2}(T-t)^{3/2}}{4\pi^2l^2},
\end{equation}

\vspace{3mm}
\begin{equation}
\label{ogo903}
C_{2r-1,2l,2l}=0,\ \ \
C_{2l-1,0,0}=-\frac{\sqrt{2}(T-t)^{3/2}}{4\pi l},\ \ \
C_{2r-1,2l-1,2l-1}=0,
\end{equation}

\vspace{3mm}
\begin{equation}
\label{ogo902}
C_{2r,2l,2l}=
\begin{cases}
-\sqrt{2}(T-t)^{3/2}/(16\pi^2l^2),\ &r=2l
\cr
\cr
0,\ & r\ne 2l\
\end{cases},
\end{equation}

\vspace{3mm}
\begin{equation}
\label{ogo902a}
C_{2r,2l-1,2l-1}=
\begin{cases}
\sqrt{2}(T-t)^{3/2}/(16\pi^2l^2),\ &r=2l
\cr
\cr
-\sqrt{2}(T-t)^{3/2}/(4\pi^2l^2),\ &r=l
\cr
\cr
0,\ & r\ne l,\ r\ne 2l
\end{cases}.
\end{equation}

\vspace{5mm}

Let us show that

\vspace{-1mm}
\begin{equation}
\label{agenty100}
\hbox{\vtop{\offinterlineskip\halign{
\hfil#\hfil\cr
{\rm l.i.m.}\cr
$\stackrel{}{{}_{p_1, p_3\to \infty}}$\cr
}} }S_{2p_1,2p_3}=
\hbox{\vtop{\offinterlineskip\halign{
\hfil#\hfil\cr
{\rm l.i.m.}\cr
$\stackrel{}{{}_{p_1, p_3\to \infty}}$\cr
}} }S_{2p_1,2p_3-1}=
\hbox{\vtop{\offinterlineskip\halign{
\hfil#\hfil\cr
{\rm l.i.m.}\cr
$\stackrel{}{{}_{p_1, p_3\to \infty}}$\cr
}} }S_{2p_1-1,2p_3-1}=
\hbox{\vtop{\offinterlineskip\halign{
\hfil#\hfil\cr
{\rm l.i.m.}\cr
$\stackrel{}{{}_{p_1, p_3\to \infty}}$\cr
}} }S_{2p_1-1,2p_3}.
\end{equation}

\vspace{4mm}

We have

\vspace{-2mm}
\begin{equation}
\label{agenty101}
S_{2p_1,2p_3}=S_{2p_1,2p_3-1}+
\sum\limits_{j_1=0}^{2p_1}
C_{2p_3, j_1, j_1}\zeta_{2p_3}^{(i_3)}.
\end{equation}

\vspace{5mm}

Using the relations (\ref{ogo901}), (\ref{ogo902}), and (\ref{ogo902a}), we obtain

\vspace{-1mm}
$$
\sum\limits_{j_1=0}^{2p_1}
C_{2p_3, j_1, j_1}=
C_{2p_3, 0, 0}+\sum\limits_{j_1=1}^{2p_1}
C_{2p_3, j_1, j_1}=
$$

\vspace{3mm}
$$
=C_{2p_3, 0, 0}+\sum\limits_{j_1=1}^{p_1}
\biggl(C_{2p_3, 2j_1-1, 2j_1-1}+C_{2p_3, 2j_1, 2j_1}\biggr)=
$$

\vspace{3mm}
\begin{equation}
\label{agenty102}
=\frac{\sqrt{2}(T-t)^{3/2}}{4\pi^2 p_3^2}\bigl(1-{\bf 1}_{\{p_1\ge p_3\}}\bigr).
\end{equation}

\vspace{5mm}

From (\ref{agenty101}), (\ref{agenty102}) we get

\vspace{-1mm}
\begin{equation}
\label{agenty103}
\hbox{\vtop{\offinterlineskip\halign{
\hfil#\hfil\cr
{\rm l.i.m.}\cr
$\stackrel{}{{}_{p_1, p_3\to \infty}}$\cr
}} }S_{2p_1,2p_3}=
\hbox{\vtop{\offinterlineskip\halign{
\hfil#\hfil\cr
{\rm l.i.m.}\cr
$\stackrel{}{{}_{p_1, p_3\to \infty}}$\cr
}} }S_{2p_1,2p_3-1}.
\end{equation}

\vspace{5mm}

Further, we have (see (\ref{ogo901})--(\ref{ogo902}))

\vspace{-1mm}
\begin{equation}
\label{agenty104}
S_{2p_1,2p_3-1}=S_{2p_1-1,2p_3-1}+
\sum\limits_{j_3=0}^{2p_3-1}
C_{j_3, 2p_1, 2p_1}\zeta_{j_3}^{(i_3)},
\end{equation}

\vspace{3mm}
$$
\sum\limits_{j_3=0}^{2p_3-1}
C_{j_3, 2p_1, 2p_1}\zeta_{j_3}^{(i_3)}=
C_{0, 2p_1, 2p_1}\zeta_{0}^{(i_3)}+
\sum\limits_{j_3=1}^{2p_3}
C_{j_3, 2p_1, 2p_1}\zeta_{j_3}^{(i_3)}-
C_{2p_3, 2p_1, 2p_1}\zeta_{2p_3}^{(i_3)}=
$$

\vspace{3mm}
$$
=
C_{0, 2p_1, 2p_1}\zeta_{0}^{(i_3)}+
\sum\limits_{j_3=1}^{p_3}\biggl(
C_{2j_3-1, 2p_1, 2p_1}\zeta_{2j_3-1}^{(i_3)}
+C_{2j_3, 2p_1, 2p_1}\zeta_{2j_3}^{(i_3)}\biggr)
-
C_{2p_3, 2p_1, 2p_1}\zeta_{2p_3}^{(i_3)}=
$$

\vspace{3mm}
\begin{equation}
\label{agenty105}
=\frac{(T-t)^{3/2}}{8\pi^2 p_1^2}\zeta_{0}^{(i_3)}+\frac{\sqrt{2}(T-t)^{3/2}}{16\pi^2 p_1^2}
\bigl({\bf 1}_{\{p_3=2p_1\}}-{\bf 1}_{\{p_3\ge 2p_1\}}\bigr)\zeta_{4p_1}^{(i_3)}.
\end{equation}

\vspace{7mm}

From (\ref{agenty104}), (\ref{agenty105}) we obtain

\vspace{-1mm}
\begin{equation}
\label{agenty106}
\hbox{\vtop{\offinterlineskip\halign{
\hfil#\hfil\cr
{\rm l.i.m.}\cr
$\stackrel{}{{}_{p_1, p_3\to \infty}}$\cr
}} }S_{2p_1,2p_3-1}=
\hbox{\vtop{\offinterlineskip\halign{
\hfil#\hfil\cr
{\rm l.i.m.}\cr
$\stackrel{}{{}_{p_1, p_3\to \infty}}$\cr
}} }S_{2p_1-1,2p_3-1}.
\end{equation}

\vspace{5mm}

Further, we have

\vspace{-2mm}
\begin{equation}
\label{agenty107}
S_{2p_1,2p_3}=S_{2p_1-1,2p_3}+
\sum\limits_{j_3=0}^{2p_3}
C_{j_3, 2p_1, 2p_1}\zeta_{j_3}^{(i_3)},
\end{equation}

\vspace{3mm}
$$
\sum\limits_{j_3=0}^{2p_3}
C_{j_3, 2p_1, 2p_1}\zeta_{j_3}^{(i_3)}=
C_{0, 2p_1, 2p_1}\zeta_{0}^{(i_3)}+
\sum\limits_{j_3=1}^{2p_3}
C_{j_3, 2p_1, 2p_1}\zeta_{j_3}^{(i_3)}=
$$

\vspace{3mm}
\begin{equation}
\label{agenty108}
=
C_{0, 2p_1, 2p_1}\zeta_{0}^{(i_3)}+
\sum\limits_{j_3=1}^{p_3}
\biggl(C_{2j_3-1, 2p_1, 2p_1}\zeta_{2j_3-1}^{(i_3)}+
C_{2j_3, 2p_1, 2p_1}\zeta_{2j_3}^{(i_3)}\biggr).
\end{equation}

\vspace{7mm}

From (\ref{agenty108}), (\ref{ogo901})--(\ref{ogo902}) we obtain

\vspace{-1mm}
\begin{equation}
\label{agenty109}
\sum\limits_{j_3=0}^{2p_3}
C_{j_3, 2p_1, 2p_1}\zeta_{j_3}^{(i_3)}=
\frac{(T-t)^{3/2}}{8\pi^2 p_1^2}
\zeta_{0}^{(i_3)}-
\frac{\sqrt{2}(T-t)^{3/2}}{16\pi^2 p_1^2}
{\bf 1}_{\{p_3\ge 2p_1\}}\zeta_{4p_1}^{(i_3)}.
\end{equation}

\vspace{5mm}

The relations (\ref{agenty107}), (\ref{agenty109}) mean that

\vspace{-1mm}
\begin{equation}
\label{agenty110}
\hbox{\vtop{\offinterlineskip\halign{
\hfil#\hfil\cr
{\rm l.i.m.}\cr
$\stackrel{}{{}_{p_1, p_3\to \infty}}$\cr
}} }S_{2p_1,2p_3}=
\hbox{\vtop{\offinterlineskip\halign{
\hfil#\hfil\cr
{\rm l.i.m.}\cr
$\stackrel{}{{}_{p_1, p_3\to \infty}}$\cr
}} }S_{2p_1-1,2p_3}.
\end{equation}

\vspace{5mm}

The equalities (\ref{agenty103}), (\ref{agenty106}), and (\ref{agenty110})
imply (\ref{agenty100}). This means that instead of (\ref{ogo1299}) it is enough
to prove the following equality

\vspace{-1mm}
\begin{equation}
\label{agenty111}
\hbox{\vtop{\offinterlineskip\halign{
\hfil#\hfil\cr
{\rm l.i.m.}\cr
$\stackrel{}{{}_{p_1, p_3\to \infty}}$\cr
}} }
\sum\limits_{j_1=0}^{2p_1}\sum\limits_{j_3=0}^{2p_3}
C_{j_3 j_1 j_1}\zeta_{j_3}^{(i_3)}
=
\frac{1}{2}\int\limits_t^T\int\limits_t^{\tau}dsd{\bf f}_{\tau}^{(i_3)}\ \ \ \hbox{w.\ p.\ 1}.
\end{equation}

\vspace{5mm}

We have

\vspace{-1mm}
$$
S_{2p_1,2p_3}=\sum\limits_{j_3=0}^{2p_3}\sum\limits_{j_1=0}^{2p_1}
C_{j_3 j_1 j_1}\zeta_{j_3}^{(i_3)}=
\frac{(T-t)^{3/2}}{6}\zeta_{0}^{(i_3)}+
$$

\vspace{3mm}
$$
+\sum\limits_{j_1=1}^{p_1}C_{0,2j_1,2j_1}
\zeta_0^{(i_3)}+\sum\limits_{j_1=1}^{p_1}C_{0,2j_1-1,2j_1-1}\zeta_0^{(i_3)}+
\sum\limits_{j_3=1}^{p_1}C_{2j_3,0,0}\zeta_{2j_3}^{(i_3)}+
$$

\vspace{3mm}
$$
+
\sum\limits_{j_3=1}^{p_3}\sum_{j_1=1}^{p_1}
C_{2j_3,2j_1,2j_1}\zeta_{2j_3}^{(i_3)}+
\sum\limits_{j_3=1}^{p_3}\sum_{j_1=1}^{p_1}
C_{2j_3,2j_1-1,2j_1-1}\zeta_{2j_3}^{(i_3)}+
\sum\limits_{j_3=1}^{p_3}
C_{2j_3-1,0,0}\zeta_{2j_3-1}^{(i_3)}+
$$

\vspace{3mm}
\begin{equation}
\label{agenty112}
+\sum\limits_{j_3=1}^{p_3}\sum_{j_1=1}^{p_1}
C_{2j_3-1,2j_1,2j_1}\zeta_{2j_3-1}^{(i_3)}+
\sum\limits_{j_3=1}^{p_3}\sum_{j_1=1}^{p_1}
C_{2j_3-1,2j_1-1,2j_1-1}\zeta_{2j_3-1}^{(i_3)}.
\end{equation}

\vspace{5mm}

After 
substituting 
(\ref{ogo901})--(\ref{ogo902a})
into (\ref{agenty112}), we obtain

\vspace{1mm}
$$
\sum\limits_{j_3=0}^{2p_3}\sum\limits_{j_1=0}^{2p_1}
C_{j_3 j_1 j_1}\zeta_{j_3}^{(i_3)}=(T-t)^{3/2}
\left(\frac{1}{6}\zeta_{0}^{(i_3)} + \frac{1}{2\pi^2}
\sum\limits_{j_1=1}^{p_1}\frac{1}{j_1^2}\ \zeta_{0}^{(i_3)}-\right.
$$

\vspace{2mm}
\begin{equation}
\label{agenty113}
\left.-
\frac{\sqrt{2}}{4\pi}\sum\limits_{j_3=1}^{p_3}\frac{1}{j_3}
\zeta_{2j_3-1}^{(i_3)}-
\frac{\sqrt{2}}{4\pi^2}\sum\limits_{j_3=1}^{\min\{p_1,p_3\}}\frac{1}{j_3^2}
\zeta_{2j_3}^{(i_3)}+
\frac{\sqrt{2}}{4\pi^2}\sum\limits_{j_3=1}^{p_3}\frac{1}{j_3^2}
\zeta_{2j_3}^{(i_3)}\right).
\end{equation}

\vspace{6mm}

From (\ref{agenty113}) we have w.~p.~1

$$
\hbox{\vtop{\offinterlineskip\halign{
\hfil#\hfil\cr
{\rm l.i.m.}\cr
$\stackrel{}{{}_{p_1, p_3\to \infty}}$\cr
}} }
\sum\limits_{j_3=0}^{2p_3}\sum\limits_{j_1=0}^{2p_1}
C_{j_3 j_1 j_1}\zeta_{j_3}^{(i_3)}=(T-t)^{3/2}
\left(\frac{1}{6}\zeta_{0}^{(i_3)} + \frac{1}{2\pi^2}
\sum\limits_{j_1=1}^{\infty}\frac{1}{j_1^2}\ \zeta_{0}^{(i_3)}-\right.
$$

\vspace{2mm}
\begin{equation}
\label{ogo905}
\left.-\hbox{\vtop{\offinterlineskip\halign{
\hfil#\hfil\cr
{\rm l.i.m.}\cr
$\stackrel{}{{}_{p_3\to \infty}}$\cr
}} }
\frac{\sqrt{2}}{4\pi}\sum\limits_{j_3=1}^{p_3}\frac{1}{j_3}
\zeta_{2j_3-1}^{(i_3)}\right).
\end{equation}

\vspace{6mm}

Using Theorems 1, 2 and the system of trigonometric functions, we get
w.~p.~1

$$
\frac{1}{2}\int\limits_t^T\int\limits_t^s d\tau d{\bf f}_{s}^{(i_3)}=
\frac{1}{2}\int\limits_t^T(s-t)d{\bf f}_{s}^{(i_3)}
=
$$

\vspace{2mm}
\begin{equation}
\label{ogo906}
=\frac{(T-t)^{3/2}}{4}\hbox{\vtop{\offinterlineskip\halign{
\hfil#\hfil\cr
{\rm l.i.m.}\cr
$\stackrel{}{{}_{p_3\to \infty}}$\cr
}} }
\Biggl(\zeta_0^{(i_3)}-\frac{\sqrt{2}}{\pi}\sum_{j_3=1}^{p_3}
\frac{1}{j_3}
\zeta_{2j_3-1}^{(i_3)}
\Biggr).
\end{equation}

\vspace{5mm}

From (\ref{ogo905}) and (\ref{ogo906}) it follows that

\vspace{1mm}
$$
\hbox{\vtop{\offinterlineskip\halign{
\hfil#\hfil\cr
{\rm l.i.m.}\cr
$\stackrel{}{{}_{p_1, p_3\to \infty}}$\cr
}} }
\sum\limits_{j_3=0}^{2p_3}\sum\limits_{j_1=0}^{2p_1}
C_{j_3 j_1 j_1}\zeta_{j_3}^{(i_3)}=
$$

\vspace{3mm}
$$
=(T-t)^{3/2}
\left(\frac{1}{6}\zeta_{0}^{(i_3)} + \frac{1}{12}
\zeta_{0}^{(i_3)}-\hbox{\vtop{\offinterlineskip\halign{
\hfil#\hfil\cr
{\rm l.i.m.}\cr
$\stackrel{}{{}_{p_3\to \infty}}$\cr
}} }
\frac{\sqrt{2}}{4\pi}\sum\limits_{j_3=1}^{p_3}\frac{1}{j_3}
\zeta_{2j_3-1}^{(i_3)}\right)=
$$

\vspace{3mm}
$$
=(T-t)^{3/2}
\left(\frac{1}{4}\zeta_{0}^{(i_3)}
-\hbox{\vtop{\offinterlineskip\halign{
\hfil#\hfil\cr
{\rm l.i.m.}\cr
$\stackrel{}{{}_{p_3\to \infty}}$\cr
}} }
\frac{\sqrt{2}}{4\pi}\sum\limits_{j_3=1}^{p_3}\frac{1}{j_3}
\zeta_{2j_3-1}^{(i_3)}\right)
=
$$

\vspace{3mm}
$$
=\frac{1}{2}\int\limits_t^T\int\limits_t^s d\tau d{\bf f}_{s}^{(i_3)},
$$

\vspace{4mm}
\noindent
where the equality is fulfilled w.~p.~1.

So, the relations (\ref{agenty111}) and (\ref{ogo1299}) are proved for the case of 
trigonometric 
system of functions.

Let us prove the relation (\ref{ogo1399}). We have

\vspace{-1mm}
$$
S'_{p_1,p_3}\stackrel{\sf def}{=}\sum\limits_{j_1=0}^{p_1}\sum\limits_{j_3=0}^{p_3}
C_{j_3 j_3 j_1}\zeta_{j_1}^{(i_1)}=
\frac{(T-t)^{3/2}}{6}\zeta_0^{(i_1)}+
$$

\vspace{3mm}
$$
+
\sum\limits_{j_3=1}^{p_3}C_{2j_3,2j_3,0}
\zeta_0^{(i_1)}+\sum\limits_{j_3=1}^{p_3}C_{2j_3-1,2j_3-1,0}\zeta_0^{(i_1)}+
\sum\limits_{j_1=1}^{p_1}\sum\limits_{j_3=1}^{p_3}
C_{2j_3,2j_3,2j_1-1}\zeta_{2j_1-1}^{(i_1)}+
$$

\vspace{3mm}
$$
+
\sum\limits_{j_1=1}^{p_1}\sum_{j_3=1}^{p_3}
C_{2j_3-1,2j_3-1,2j_1-1}\zeta_{2j_1-1}^{(i_1)}+
\sum\limits_{j_1=1}^{p_1}
C_{0,0,2j_1-1}\zeta_{2j_1-1}^{(i_1)}+
\sum\limits_{j_1=1}^{p_1}\sum\limits_{j_3=1}^{p_3}
C_{2j_3,2j_3,2j_1}\zeta_{2j_1}^{(i_1)}+
$$

\vspace{3mm}
\begin{equation}
\label{ogo9000}
+\sum\limits_{j_1=1}^{p_1}\sum_{j_3=1}^{p_3}
C_{2j_3-1,2j_3-1,2j_1}\zeta_{2j_1}^{(i_1)}+
\sum_{j_1=1}^{p_1}
C_{0,0,2j_1}\zeta_{2j_1}^{(i_1)},
\end{equation}

\vspace{7mm}
\noindent
where the summation is stopped, when
$2j_3,$ $2j_3-1> p_3$
or $2j_1,$ $2j_1-1> p_1$ and

\vspace{1mm}
\begin{equation}
\label{ogo9010}
C_{2l,2l,0}=\frac{(T-t)^{3/2}}{8\pi^2l^2},\ \ \
C_{2l-1,2l-1,0}=\frac{3(T-t)^{3/2}}{8\pi^2l^2},\ \ \
C_{0,0,2r}=\frac{\sqrt{2}(T-t)^{3/2}}{4\pi^2r^2},
\end{equation}

\vspace{3mm}
\begin{equation}
\label{ogo9030}
C_{2l-1,2l-1,2r-1}=0,\ \ \
C_{0,0,2r-1}=\frac{\sqrt{2}(T-t)^{3/2}}{4\pi r},\ \ \
C_{2l,2l,2r-1}=0,
\end{equation}

\vspace{3mm}
\begin{equation}
\label{ogo9020}
C_{2l,2l,2r}=
\begin{cases}
-\sqrt{2}(T-t)^{3/2}/(16\pi^2l^2),\ &r=2l
\cr
\cr
0,\  &r\ne 2l
\end{cases},
\end{equation}

\vspace{3mm}
\begin{equation}
\label{ogo9020a}
C_{2l-1,2l-1,2r}=
\begin{cases}
\sqrt{2}(T-t)^{3/2}/(16\pi^2l^2),\ &r=2l
\cr
\cr
-\sqrt{2}(T-t)^{3/2}/(4\pi^2l^2),\ &r=l
\cr
\cr
0,\  &r\ne l,\ r\ne 2l
\end{cases}.
\end{equation}

\vspace{6mm}

Let us show that

\vspace{-1mm}
\begin{equation}
\label{agenty1000}
\hbox{\vtop{\offinterlineskip\halign{
\hfil#\hfil\cr
{\rm l.i.m.}\cr
$\stackrel{}{{}_{p_1, p_3\to \infty}}$\cr
}} }S'_{2p_1,2p_3}=
\hbox{\vtop{\offinterlineskip\halign{
\hfil#\hfil\cr
{\rm l.i.m.}\cr
$\stackrel{}{{}_{p_1, p_3\to \infty}}$\cr
}} }S'_{2p_1,2p_3-1}=
\hbox{\vtop{\offinterlineskip\halign{
\hfil#\hfil\cr
{\rm l.i.m.}\cr
$\stackrel{}{{}_{p_1, p_3\to \infty}}$\cr
}} }S'_{2p_1-1,2p_3-1}=
\hbox{\vtop{\offinterlineskip\halign{
\hfil#\hfil\cr
{\rm l.i.m.}\cr
$\stackrel{}{{}_{p_1, p_3\to \infty}}$\cr
}} }S'_{2p_1-1,2p_3}.
\end{equation}

\vspace{4mm}

We have

\vspace{-2mm}
\begin{equation}
\label{agenty1010}
S'_{2p_1,2p_3}=S'_{2p_1-1,2p_3}+
\sum\limits_{j_3=0}^{2p_3}
C_{j_3, j_3, 2p_1}\zeta_{2p_1}^{(i_1)}.
\end{equation}

\vspace{5mm}

Using the relations (\ref{ogo9010}), (\ref{ogo9020}), and (\ref{ogo9020a}), we obtain

\vspace{-1mm}
$$
\sum\limits_{j_1=0}^{2p_3}
C_{j_3, j_3, 2p_1}=
C_{0, 0, 2p_1}+\sum\limits_{j_3=1}^{2p_3}
C_{j_3, j_3, 2p_1}=
$$

\vspace{3mm}
$$
=C_{0, 0, 2p_1}+\sum\limits_{j_3=1}^{p_3}
\biggl(C_{2j_3-1, 2j_3-1, 2p_1}+C_{2j_3, 2j_3, 2p_1}\biggr)=
$$

\vspace{3mm}
\begin{equation}
\label{agenty1020}
=\frac{\sqrt{2}(T-t)^{3/2}}{4\pi^2 p_1^2}\bigl(1-{\bf 1}_{\{p_3\ge p_1\}}\bigr).
\end{equation}

\vspace{5mm}

From (\ref{agenty1010}), (\ref{agenty1020}) we obtain

\vspace{-1mm}
\begin{equation}
\label{agenty1030}
\hbox{\vtop{\offinterlineskip\halign{
\hfil#\hfil\cr
{\rm l.i.m.}\cr
$\stackrel{}{{}_{p_1, p_3\to \infty}}$\cr
}} }S'_{2p_1,2p_3}=
\hbox{\vtop{\offinterlineskip\halign{
\hfil#\hfil\cr
{\rm l.i.m.}\cr
$\stackrel{}{{}_{p_1, p_3\to \infty}}$\cr
}} }S'_{2p_1-1,2p_3}.
\end{equation}

\vspace{5mm}

Further, we get (see (\ref{ogo9010})--(\ref{ogo9020}))

\vspace{-1mm}
\begin{equation}
\label{agenty1040}
S'_{2p_1-1,2p_3}=S'_{2p_1-1,2p_3-1}+
\sum\limits_{j_1=0}^{2p_1-1}
C_{2p_3, 2p_3, j_1}\zeta_{j_1}^{(i_1)},
\end{equation}

\vspace{3mm}
$$
\sum\limits_{j_1=0}^{2p_1-1}
C_{2p_3, 2p_3, j_1}\zeta_{j_1}^{(i_1)}=
C_{2p_3, 2p_3, 0}\zeta_{0}^{(i_1)}+
\sum\limits_{j_1=1}^{2p_1}
C_{2p_3, 2p_3, j_1}\zeta_{j_1}^{(i_1)}-
C_{2p_3, 2p_3, 2p_1}\zeta_{2p_1}^{(i_1)}=
$$

\vspace{3mm}
$$
=
C_{2p_3, 2p_3, 0}\zeta_{0}^{(i_1)}+
\sum\limits_{j_1=1}^{p_1}\biggl(
C_{2p_3, 2p_3, 2j_1-1}\zeta_{2j_1-1}^{(i_1)}
+C_{2p_3, 2p_3, 2j_1}\zeta_{2j_1}^{(i_1)}\biggr)
-
C_{2p_3, 2p_3, 2p_1}\zeta_{2p_1}^{(i_1)}=
$$

\vspace{3mm}
\begin{equation}
\label{agenty1050}
=\frac{(T-t)^{3/2}}{8\pi^2 p_3^2}\zeta_{0}^{(i_1)}+\frac{\sqrt{2}(T-t)^{3/2}}{16\pi^2 p_3^2}
\bigl({\bf 1}_{\{p_1=2p_3\}}-{\bf 1}_{\{p_1\ge 2p_3\}}\bigr)\zeta_{4p_3}^{(i_1)}.
\end{equation}

\vspace{6mm}

From (\ref{agenty1040}), (\ref{agenty1050}) we obtain

\vspace{-1mm}
\begin{equation}
\label{agenty1060}
\hbox{\vtop{\offinterlineskip\halign{
\hfil#\hfil\cr
{\rm l.i.m.}\cr
$\stackrel{}{{}_{p_1, p_3\to \infty}}$\cr
}} }S'_{2p_1-1,2p_3}=
\hbox{\vtop{\offinterlineskip\halign{
\hfil#\hfil\cr
{\rm l.i.m.}\cr
$\stackrel{}{{}_{p_1, p_3\to \infty}}$\cr
}} }S'_{2p_1-1,2p_3-1}.
\end{equation}

\vspace{4mm}

Further, we have

\vspace{-1mm}
\begin{equation}
\label{agenty1070}
S'_{2p_1,2p_3}=S'_{2p_1,2p_3-1}+
\sum\limits_{j_1=0}^{2p_1}
C_{2p_3, 2p_3, j_1}\zeta_{j_1}^{(i_1)},
\end{equation}

\vspace{3mm}
$$
\sum\limits_{j_1=0}^{2p_1}
C_{2p_3, 2p_3, j_1}\zeta_{j_1}^{(i_1)}=
C_{2p_3, 2p_3, 0}\zeta_{0}^{(i_1)}+
\sum\limits_{j_1=1}^{2p_1}
C_{2p_3, 2p_3, j_1}\zeta_{j_1}^{(i_1)}=
$$

\vspace{3mm}
\begin{equation}
\label{agenty1080}
=
C_{2p_3, 2p_3, 0}\zeta_{0}^{(i_1)}+
\sum\limits_{j_1=1}^{p_1}
\biggl(C_{2p_3, 2p_3, 2j_1-1}\zeta_{2j_1-1}^{(i_1)}+
C_{2p_3, 2p_3, 2j_1}\zeta_{2j_1}^{(i_1)}\biggr).
\end{equation}

\vspace{6mm}

From (\ref{agenty1080}), (\ref{ogo9010})--(\ref{ogo9020}) we obtain

\vspace{-1mm}
\begin{equation}
\label{agenty1090}
\sum\limits_{j_1=0}^{2p_1}
C_{2p_3, 2p_3, j_1}\zeta_{j_1}^{(i_1)}=
\frac{(T-t)^{3/2}}{8\pi^2 p_3^2}
\zeta_{0}^{(i_1)}-
\frac{\sqrt{2}(T-t)^{3/2}}{16\pi^2 p_3^2}
{\bf 1}_{\{p_1\ge 2p_3\}}\zeta_{4p_3}^{(i_1)}.
\end{equation}

\vspace{4mm}

The relations (\ref{agenty1070}), (\ref{agenty1090}) mean that

\vspace{-1mm}
\begin{equation}
\label{agenty1100}
\hbox{\vtop{\offinterlineskip\halign{
\hfil#\hfil\cr
{\rm l.i.m.}\cr
$\stackrel{}{{}_{p_1, p_3\to \infty}}$\cr
}} }S'_{2p_1,2p_3}=
\hbox{\vtop{\offinterlineskip\halign{
\hfil#\hfil\cr
{\rm l.i.m.}\cr
$\stackrel{}{{}_{p_1, p_3\to \infty}}$\cr
}} }S'_{2p_1,2p_3-1}.
\end{equation}

\vspace{4mm}

The equalities (\ref{agenty1030}), (\ref{agenty1060}), and (\ref{agenty1100})
imply (\ref{agenty1000}). This means that instead of (\ref{ogo1399}) it is enough
to prove the following equality

\vspace{-1mm}
\begin{equation}
\label{agenty1110}
\hbox{\vtop{\offinterlineskip\halign{
\hfil#\hfil\cr
{\rm l.i.m.}\cr
$\stackrel{}{{}_{p_1, p_3\to \infty}}$\cr
}} }
\sum\limits_{j_1=0}^{2p_1}\sum\limits_{j_3=0}^{2p_3}
C_{j_3 j_3 j_1}\zeta_{j_1}^{(i_1)}
=
\frac{1}{2}\int\limits_t^T\int\limits_t^{\tau}d{\bf f}_{s}^{(i_1)}d\tau\ \ \ \hbox{w.\ p.\ 1}.
\end{equation}

\vspace{4mm}

We have

\vspace{-1mm}
$$
S'_{2p_1,2p_3}=\sum\limits_{j_1=0}^{2p_1}\sum\limits_{j_3=0}^{2p_3}
C_{j_3 j_3 j_1}\zeta_{j_1}^{(i_1)}=
\frac{(T-t)^{3/2}}{6}\zeta_0^{(i_1)}+
$$

\vspace{3mm}
$$
+
\sum\limits_{j_3=1}^{p_3}C_{2j_3,2j_3,0}
\zeta_0^{(i_1)}+\sum\limits_{j_3=1}^{p_3}C_{2j_3-1,2j_3-1,0}\zeta_0^{(i_1)}+
\sum\limits_{j_1=1}^{p_1}\sum\limits_{j_3=1}^{p_3}
C_{2j_3,2j_3,2j_1-1}\zeta_{2j_1-1}^{(i_1)}+
$$

\vspace{3mm}
$$
+
\sum\limits_{j_1=1}^{p_1}\sum_{j_3=1}^{p_3}
C_{2j_3-1,2j_3-1,2j_1-1}\zeta_{2j_1-1}^{(i_1)}+
\sum\limits_{j_1=1}^{p_1}
C_{0,0,2j_1-1}\zeta_{2j_1-1}^{(i_1)}+
\sum\limits_{j_1=1}^{p_1}\sum\limits_{j_3=1}^{p_3}
C_{2j_3,2j_3,2j_1}\zeta_{2j_1}^{(i_1)}+
$$

\vspace{3mm}
\begin{equation}
\label{agenty1120}
+\sum\limits_{j_1=1}^{p_1}\sum_{j_3=1}^{p_3}
C_{2j_3-1,2j_3-1,2j_1}\zeta_{2j_1}^{(i_1)}+
\sum_{j_1=1}^{p_1}
C_{0,0,2j_1}\zeta_{2j_1}^{(i_1)}.
\end{equation}

\vspace{7mm}

After 
substituting 
(\ref{ogo9010})--(\ref{ogo9020a})
into (\ref{agenty1120}), we obtain

\vspace{1mm}
$$
\sum\limits_{j_1=0}^{2p_1}\sum\limits_{j_3=0}^{2p_3}
C_{j_3 j_3 j_1}\zeta_{j_1}^{(i_1)}=(T-t)^{3/2}
\left(\frac{1}{6}\zeta_{0}^{(i_1)} + \frac{1}{2\pi^2}
\sum\limits_{j_3=1}^{p_3}\frac{1}{j_3^2}\ \zeta_{0}^{(i_1)}+\right.
$$

\vspace{2mm}
\begin{equation}
\label{agenty1130}
\left.+
\frac{\sqrt{2}}{4\pi}\sum\limits_{j_1=1}^{p_1}\frac{1}{j_1}
\zeta_{2j_1-1}^{(i_1)}-
\frac{\sqrt{2}}{4\pi^2}\sum\limits_{j_1=1}^{\min\{p_1,p_3\}}\frac{1}{j_1^2}
\zeta_{2j_1}^{(i_1)}+
\frac{\sqrt{2}}{4\pi^2}\sum\limits_{j_1=1}^{p_1}\frac{1}{j_1^2}
\zeta_{2j_1}^{(i_1)}\right).
\end{equation}

\vspace{5mm}

From (\ref{agenty1130}) we have w.~p.~1

$$
\hbox{\vtop{\offinterlineskip\halign{
\hfil#\hfil\cr
{\rm l.i.m.}\cr
$\stackrel{}{{}_{p_1, p_3\to \infty}}$\cr
}} }
\sum\limits_{j_1=0}^{2p_1}\sum\limits_{j_3=0}^{2p_3}
C_{j_3 j_3 j_1}\zeta_{j_1}^{(i_1)}=(T-t)^{3/2}
\left(\frac{1}{6}\zeta_{0}^{(i_3)} + \frac{1}{2\pi^2}
\sum\limits_{j_3=1}^{\infty}\frac{1}{j_3^2}\ \zeta_{0}^{(i_1)}+\right.
$$

\vspace{2mm}
\begin{equation}
\label{ogo9050}
\left.+\hbox{\vtop{\offinterlineskip\halign{
\hfil#\hfil\cr
{\rm l.i.m.}\cr
$\stackrel{}{{}_{p_1\to \infty}}$\cr
}} }
\frac{\sqrt{2}}{4\pi}\sum\limits_{j_1=1}^{p_1}\frac{1}{j_1}
\zeta_{2j_1-1}^{(i_1)}\right).
\end{equation}

\vspace{5mm}

Using the Ito formula and Theorems 1, 2 for the case of 
trigonometric system of 
functions, we obtain w.~p.~1

$$
\frac{1}{2}\int\limits_t^T\int\limits_t^{\tau}d{\bf f}_{s}^{(i_1)}d\tau=
\frac{1}{2}\left((T-t)
\int\limits_t^Td{\bf f}_{s}^{(i_1)}+
\int\limits_t^T(t-s)d{\bf f}_{s}^{(i_1)}\right)=
$$

\vspace{2mm}
\begin{equation}
\label{ogo9060}
=\frac{1}{4}(T-t)^{3/2}
\Biggl(\zeta_0^{(i_1)}+\hbox{\vtop{\offinterlineskip\halign{
\hfil#\hfil\cr
{\rm l.i.m.}\cr
$\stackrel{}{{}_{p_1\to \infty}}$\cr
}} }
\frac{\sqrt{2}}{\pi}\sum_{j_1=1}^{p_1}
\frac{1}{j_1}
\zeta_{2j_1-1}^{(i_1)}
\Biggr).
\end{equation}

\vspace{5mm}

From (\ref{ogo9050}) and (\ref{ogo9060}) it follows that 

\vspace{1mm}
$$
\hbox{\vtop{\offinterlineskip\halign{
\hfil#\hfil\cr
{\rm l.i.m.}\cr
$\stackrel{}{{}_{p_1, p_3\to \infty}}$\cr
}} }
\sum\limits_{j_1=0}^{2p_1}\sum\limits_{j_3=0}^{2p_3}
C_{j_3 j_3 j_1}\zeta_{j_1}^{(i_1)}=
$$

\vspace{3mm}
$$
=(T-t)^{3/2}\left(
\frac{1}{6}\zeta_{0}^{(i_1)}+\frac{1}{12}
\zeta_{0}^{(i_1)}+
\hbox{\vtop{\offinterlineskip\halign{
\hfil#\hfil\cr
{\rm l.i.m.}\cr
$\stackrel{}{{}_{p_1\to \infty}}$\cr
}} }
\frac{\sqrt{2}}{4\pi}\sum\limits_{j_1=1}^{p_1}\frac{1}{j_1}
\zeta_{2j_1-1}^{(i_1)}\right)=
$$

\vspace{3mm}
$$
=(T-t)^{3/2}\left(
\frac{1}{4}\zeta_{0}^{(i_1)}+
\hbox{\vtop{\offinterlineskip\halign{
\hfil#\hfil\cr
{\rm l.i.m.}\cr
$\stackrel{}{{}_{p_1\to \infty}}$\cr
}} }
\frac{\sqrt{2}}{4\pi}\sum\limits_{j_1=1}^{p_1}\frac{1}{j_1}
\zeta_{2j_1-1}^{(i_1)}\right)=
$$

\vspace{3mm}
$$
=
\frac{1}{2}\int\limits_t^T\int\limits_t^{\tau}d{\bf f}_{s}^{(i_1)}d\tau,
$$

\vspace{4mm}
\noindent
where the equality is fulfilled w.~p.~1.

So, the relations (\ref{agenty1110}) and (\ref{ogo1399}) are proved for the case 
of trigonometric system of functions.

Let us prove the equality (\ref{ogo13a99}).
Since $\psi_1(\tau),$ $\psi_2(\tau),$ $\psi_3(\tau)\equiv 1,$
then the following 
relation 
for the Fourier coefficients is correct

\vspace{-1mm}
$$
C_{j_1 j_1 j_3}+C_{j_1 j_3 j_1}+C_{j_3 j_1 j_1}=\frac{1}{2}
C_{j_1}^2 C_{j_3}.
$$ 

\vspace{4mm}

Then w.\ p.\ 1 

\vspace{-1mm}
$$
\hbox{\vtop{\offinterlineskip\halign{
\hfil#\hfil\cr
{\rm l.i.m.}\cr
$\stackrel{}{{}_{p_1,p_3\to \infty}}$\cr
}} }\sum_{j_1=0}^{p_1}\sum_{j_3=0}^{p_3}
C_{j_1 j_3 j_1}\zeta_{j_3}^{(i_2)}=
$$

\vspace{3mm}
\begin{equation}
\label{ogo2010}
=
\hbox{\vtop{\offinterlineskip\halign{
\hfil#\hfil\cr
{\rm l.i.m.}\cr
$\stackrel{}{{}_{p_1,p_3\to \infty}}$\cr
}} }\sum_{j_1=0}^{p_1}\sum_{j_3=0}^{p_3}
\left(\frac{1}{2}C_{j_1}^2 C_{j_3}-C_{j_1 j_1 j_3}-C_{j_3 j_1 j_1}
\right)\zeta_{j_3}^{(i_2)}.
\end{equation}

\vspace{5mm}

Taking into account (\ref{ogo1299}) and (\ref{ogo1399}), 
we can write w.\ p.\ 1 

\vspace{-1mm}
$$
\hbox{\vtop{\offinterlineskip\halign{
\hfil#\hfil\cr
{\rm l.i.m.}\cr
$\stackrel{}{{}_{p_1,p_3\to \infty}}$\cr
}} }\sum_{j_1=0}^{p_1}\sum_{j_3=0}^{p_3}
C_{j_1 j_3 j_1}\zeta_{j_3}^{(i_2)}=
$$

\vspace{3mm}
$$
=
\frac{1}{2}C_0^3\zeta_0^{(i_2)}
-
\hbox{\vtop{\offinterlineskip\halign{
\hfil#\hfil\cr
{\rm l.i.m.}\cr
$\stackrel{}{{}_{p_1,p_3\to \infty}}$\cr
}} }\sum_{j_1=0}^{p_1}\sum_{j_3=0}^{p_3}
C_{j_1 j_1 j_3}\zeta_{j_3}^{(i_2)}-
$$

\vspace{3mm}
$$
-
\hbox{\vtop{\offinterlineskip\halign{
\hfil#\hfil\cr
{\rm l.i.m.}\cr
$\stackrel{}{{}_{p_1,p_3\to \infty}}$\cr
}} }\sum_{j_1=0}^{p_1}\sum_{j_3=0}^{p_3}
C_{j_3 j_1 j_1}\zeta_{j_3}^{(i_2)}=
$$

\vspace{3mm}
$$
=
\frac{1}{2}(T-t)^{3/2}
\zeta_0^{(i_2)}-
\frac{1}{4}(T-t)^{3/2}
\Biggl(\zeta_0^{(i_2)}+
\hbox{\vtop{\offinterlineskip\halign{
\hfil#\hfil\cr
{\rm l.i.m.}\cr
$\stackrel{}{{}_{p_1\to \infty}}$\cr
}} }
\frac{\sqrt{2}}{\pi}\sum_{j_1=1}^{p_1}
\frac{1}{j_1}
\zeta_{2j_1-1}^{(i_2)}
\Biggr)
-
$$

\vspace{3mm}
$$
-\frac{1}{4}(T-t)^{3/2}
\Biggl(\zeta_0^{(i_2)}-\hbox{\vtop{\offinterlineskip\halign{
\hfil#\hfil\cr
{\rm l.i.m.}\cr
$\stackrel{}{{}_{p_1\to \infty}}$\cr
}} }\frac{\sqrt{2}}{\pi}
\sum_{j_1=1}^{p_1}
\frac{1}{j_1}
\zeta_{2j_1-1}^{(i_2)}
\Biggr)=0.
$$

\vspace{5mm}

From Theorems 1, 2 and (\ref{ogo1299})--(\ref{ogo13a99}) 
we obtain the expansion
(\ref{feto19001ee}). Theorem 5 is proved.

\vspace{5mm}

\section{Modifications of Theorems 3--5}

\vspace{5mm}

Let us consider the following modification of Theorem 4.

\vspace{2mm}

{\bf Theorem 6}\ \cite{14}-\cite{16}, \cite{19}, \cite{20}-\cite{20xx2}. 
{\it Assume that
$\{\phi_j(x)\}_{j=0}^{\infty}$ is a complete orthonormal
system of Legendre polynomials
in the space $L_2([t, T])$
and $\psi_1(\tau),$ $\psi_2(\tau),$ $\psi_3(\tau)$ are continuously
differentiable functions at the interval $[t, T]$.
Then, for the 
iterated Stratonovich stochastic integral of third multiplicity

$$
J^{*}[\psi^{(3)}]_{T,t}={\int\limits_t^{*}}^T\psi_3(t_3)
{\int\limits_t^{*}}^{t_3}\psi_2(t_2)
{\int\limits_t^{*}}^{t_2}\psi_1(t_1)
d{\bf f}_{t_1}^{(i_1)}
d{\bf f}_{t_2}^{(i_2)}d{\bf f}_{t_3}^{(i_3)}\ \ \ (i_1, i_2, i_3=1,\ldots,m)
$$

\vspace{3mm}
\noindent
the following 
expansion 

\begin{equation}
\label{feto19000}
J^{*}[\psi^{(3)}]_{T,t}=
\hbox{\vtop{\offinterlineskip\halign{
\hfil#\hfil\cr
{\rm l.i.m.}\cr
$\stackrel{}{{}_{p\to \infty}}$\cr
}} }
\sum\limits_{j_1, j_2, j_3=0}^{p}
C_{j_3 j_2 j_1}\zeta_{j_1}^{(i_1)}\zeta_{j_2}^{(i_2)}\zeta_{j_3}^{(i_3)}
\end{equation}

\vspace{3mm}
\noindent
converging in the mean-square sense 
is valid for each of the following cases

\vspace{2mm}
\noindent
{\rm 1}.\ $i_1\ne i_2,\ i_2\ne i_3,\ i_1\ne i_3,$\\
{\rm 2}.\ $i_1=i_2\ne i_3$ and
$\psi_1(\tau)\equiv\psi_2(\tau),$\\
{\rm 3}.\ $i_1\ne i_2=i_3$ and
$\psi_2(\tau)\equiv\psi_3(\tau),$\\
{\rm 4}.\ $i_1, i_2, i_3=1,\ldots,m$
and $\psi_1(\tau)\equiv\psi_2(\tau)\equiv\psi_3(\tau),$\

\vspace{3mm}
\noindent
where
$$
C_{j_3 j_2 j_1}=\int\limits_t^T\psi_3(s)\phi_{j_3}(s)
\int\limits_t^s\psi_2(s_1)\phi_{j_2}(s_1)
\int\limits_t^{s_1}\psi_1(s_2)\phi_{j_1}(s_2)ds_2ds_1ds.
$$
}

\vspace{2mm}

{\bf Proof.}\ Case 1
directly follows from Theorems 1, 2. Let us consider Case 2.
We will prove w.~p.~1 the following relation

\vspace{-1mm}
$$
\hbox{\vtop{\offinterlineskip\halign{
\hfil#\hfil\cr
{\rm l.i.m.}\cr
$\stackrel{}{{}_{p\to \infty}}$\cr
}} }
\sum\limits_{j_1=0}^{p}\sum\limits_{j_3=0}^{p}
C_{j_3 j_1 j_1}\zeta_{j_3}^{(i_3)}=
\frac{1}{2}\int\limits_t^T\psi_3(s)
\int\limits_t^s\psi^2(s_1)ds_1d{\bf f}_s^{(i_3)},
$$

\vspace{3mm}
\noindent
where
$$
C_{j_3 j_1 j_1}=\int\limits_t^T\psi_3(s)\phi_{j_3}(s)
\int\limits_t^s\psi(s_1)\phi_{j_1}(s_1)
\int\limits_t^{s_1}\psi(s_2)\phi_{j_1}(s_2)ds_2ds_1ds.
$$

\vspace{5mm}

Using Theorems 
1, 2 we can write the following

\vspace{1mm}
$$
\frac{1}{2}\int\limits_t^T\psi_3(s)
\int\limits_t^s\psi^2(s_1)ds_1d{\bf f}_s^{(i_3)}=
\frac{1}{2}
\hbox{\vtop{\offinterlineskip\halign{
\hfil#\hfil\cr
{\rm l.i.m.}\cr
$\stackrel{}{{}_{p_3\to \infty}}$\cr
}} }\sum_{j_3=0}^{p_3}
\tilde C_{j_3}\zeta_{j_3}^{(i_3)},
$$

\vspace{5mm}
\noindent
where 
$$
\tilde C_{j_3}=
\int\limits_t^T
\phi_{j_3}(s)\psi_3(s)\int\limits_t^s\psi^2(s_1)ds_1ds.
$$

\vspace{3mm}

We have

\vspace{1mm}
$$
{\sf M}\left\{\left(\sum_{j_3=0}^p\left(\sum_{j_1=0}^p
C_{j_3j_1j_1}-\frac{1}{2}\tilde C_{j_3}\right)
\zeta_{j_3}^{(i_3)}\right)^2\right\}=
\sum_{j_3=0}^p\left(\sum\limits_{j_1=0}^{p}C_{j_3j_1 j_1}-
\frac{1}{2}\tilde C_{j_3}\right)^2=
$$

\vspace{1mm}
$$
=\sum_{j_3=0}^p\left(\frac{1}{2}\sum\limits_{j_1=0}^{p}
\int\limits_t^T\phi_{j_3}(s)\psi_3(s)
\left(\int\limits_t^s\phi_{j_1}(s_1)\psi(s_1)ds_1\right)^2ds-
\frac{1}{2}
\int\limits_t^T
\phi_{j_3}(s)\psi_3(s)\int\limits_t^s\psi^2(s_1)ds_1ds\right)^2=
$$

\vspace{1mm}
$$
=\frac{1}{4}\sum_{j_3=0}^p\left(
\int\limits_t^T\phi_{j_3}(s)\psi_3(s)\left(
\sum\limits_{j_1=0}^{p}
\left(\int\limits_t^s\phi_{j_1}(s_1)\psi(s_1)ds_1\right)^2
-\int\limits_t^s\psi^2(s_1)ds_1\right)ds\right)^2=
$$

\vspace{1mm}

\begin{equation}
\label{otit5000}
=\frac{1}{4}\sum_{j_3=0}^p\left(
\int\limits_t^T\phi_{j_3}(s)\psi_3(s)
\sum\limits_{j_1=p+1}^{\infty}
\left(\int\limits_t^s\phi_{j_1}(s_1)\psi(s_1)ds_1\right)^2
ds\right)^2.
\end{equation}

\vspace{6mm}

In order to get (\ref{otit5000}) we used the Parseval equality
in the form

\vspace{-1mm}
$$
\sum_{j_1=0}^{\infty}\left(\int\limits_t^s\phi_{j_1}(s_1)
\psi(s_1)ds_1\right)^2=
\int\limits_t^T K^2(s,s_1)ds_1=\int\limits_t^s\psi^2(s_1)ds_1,
$$

\vspace{3mm}
\noindent
where
$$
K(s,s_1)=\psi(s_1){\bf 1}_{\{s_1< s\}},\ \ \ s, s_1\in [t, T].
$$

\vspace{4mm}

We have for $j_1\in\mathbb{N}$

\vspace{-2mm}
$$
\left(\int\limits_t^s\psi(s_1)\phi_{j_1}(s_1)ds_1\right)^2=
$$

$$
=
\frac{(T-t)(2j_1+1)}{4}
\left(\int\limits_{-1}^{z(s)}P_{j_1}(y)
\psi\left(\frac{T-t}{2}y+\frac{T+t}{2}\right)dy\right)^2=
$$

$$
=\frac{T-t}{4(2j_1+1)}\Biggl(\left(P_{j_1+1}(z(s))-
P_{j_1-1}(z(s))\right)\psi(s)-\Biggr.
$$

\begin{equation}
\label{otit6000}
\Biggl.-\frac{T-t}{2}
\int\limits_{-1}^{z(s)}(\left(P_{j_1+1}(y)-P_{j_1-1}(y)\right)
\psi'\left(\frac{T-t}{2}y+\frac{T+t}{2}\right)dy\Biggr)^2,
\end{equation}

\vspace{5mm}
\noindent
where 
$$
z(s)=\left(s-\frac{T+t}{2}\right)\frac{2}{T-t},
$$

\vspace{5mm}
\noindent
and $\psi'$ is a derivative of the function $\psi(s)$
with respect to the variable 

\vspace{1mm}
$$
\frac{T-t}{2}y+\frac{T+t}{2}.
$$

\vspace{4mm}

Further consideration is similar to
the proof of Case 2 from Theorem 4. 
Finally, from (\ref{otit5000}) and (\ref{otit6000})
we obtain

$$
{\sf M}\left\{\left(\sum_{j_3=0}^p\left(\sum_{j_1=0}^p
C_{j_3j_1j_1}-\frac{1}{2}\tilde C_{j_3}\right)
\zeta_{j_3}^{(i_3)}\right)^2\right\}<
$$

\vspace{1mm}
$$
< K\frac{p}{p^2}\left(
\int\limits_{-1}^1\frac{dy}{(1-y^2)^{3/4}}
+\int\limits_{-1}^1\frac{dy}{(1-y^2)^{1/4}}\right)^2\le
$$

\vspace{2mm}
$$
\le \frac{K_1}{p}\ \to 0\ \ \  \hbox{if}\ \ \ p\ \to \infty,
$$

\vspace{5mm}
\noindent
where $K, K_1$ are constants. Case 2 is proved.

Let us consider Case 3. In this case we will prove w.~p.~1 the following
relation

\vspace{1mm}
$$
\hbox{\vtop{\offinterlineskip\halign{
\hfil#\hfil\cr
{\rm l.i.m.}\cr
$\stackrel{}{{}_{p\to \infty}}$\cr
}} }
\sum\limits_{j_1=0}^{p}\sum\limits_{j_3=0}^{p}
C_{j_3 j_3 j_1}\zeta_{j_1}^{(i_1)}=
\frac{1}{2}\int\limits_t^T\psi^2(s)
\int\limits_t^s\psi_1(s_1)d{\bf f}_{s_1}^{(i_1)}ds,
$$

\vspace{3mm}
\noindent
where
$$
C_{j_3 j_3 j_1}=\int\limits_t^T\psi(s)\phi_{j_3}(s)
\int\limits_t^s\psi(s_1)\phi_{j_3}(s_1)
\int\limits_t^{s_1}\psi_1(s_2)\phi_{j_1}(s_2)ds_2ds_1ds.
$$

\vspace{5mm}

Using the Ito formula, we obtain w.~p.~1

\begin{equation}
\label{otit7000}
\frac{1}{2}\int\limits_t^T\psi^2(s)
\int\limits_t^s\psi_1(s_1)d{\bf f}_{s_1}^{(i_1)}ds=
\frac{1}{2}\int\limits_t^T\psi_1(s_1)
\int\limits_{s_1}^T\psi^2(s)dsd{\bf f}_{s_1}^{(i_1)}.
\end{equation}

\vspace{5mm}

Applying Theorems 
1 and 2, we have

\begin{equation}
\label{otit9000}
\frac{1}{2}\int\limits_t^T\psi_1(s_1)
\int\limits_{s_1}^T\psi^2(s)dsd{\bf f}_{s_1}^{(i_1)}=
\frac{1}{2}
\hbox{\vtop{\offinterlineskip\halign{
\hfil#\hfil\cr
{\rm l.i.m.}\cr
$\stackrel{}{{}_{p_1\to \infty}}$\cr
}} }\sum_{j_1=0}^{p_1}
C_{j_1}^{*}\zeta_{j_1}^{(i_1)},
\end{equation}

\vspace{4mm}
\noindent
where
$$
C_{j_1}^{*}=
\int\limits_t^T
\phi_{j_1}(s_1)\psi_1(s_1)\int\limits_{s_1}^T\psi^2(s)dsds_1.
$$

\vspace{4mm}

Moreover,
$$
C_{j_3 j_3 j_1}=\int\limits_t^T\psi(s)\phi_{j_3}(s)
\int\limits_t^s\psi(s_1)\phi_{j_3}(s_1)
\int\limits_t^{s_1}\psi_1(s_2)\phi_{j_1}(s_2)ds_2ds_1ds=
$$

\vspace{1mm}
$$
=\int\limits_t^T\psi_1(s_2)\phi_{j_1}(s_2)
\int\limits_{s_2}^T\psi(s_1)\phi_{j_3}(s_1)
\int\limits_{s_1}^{T}\psi(s)\phi_{j_3}(s)dsds_1ds_2=
$$

\vspace{1mm}
\begin{equation}
\label{otit8000}
=\frac{1}{2}\int\limits_t^T\psi_1(s_2)\phi_{j_1}(s_2)
\left(\int\limits_{s_2}^T\psi(s_1)\phi_{j_3}(s_1)ds_1\right)^2
ds_2.
\end{equation}

\vspace{5mm}

From (\ref{otit7000})--(\ref{otit8000})
we obtain

\vspace{1mm}
$$
{\sf M}\left\{\left(\sum_{j_1=0}^p\left(\sum_{j_3=0}^p
C_{j_3j_3j_1}-\frac{1}{2} C_{j_1}^{*}\right)
\zeta_{j_1}^{(i_1)}\right)^2\right\}=
\sum_{j_1=0}^p\left(\sum\limits_{j_3=0}^{p}C_{j_3j_3 j_1}-
\frac{1}{2} C_{j_1}^{*}\right)^2=
$$

\vspace{2mm}
$$
=\frac{1}{4}\sum_{j_1=0}^p\left(
\int\limits_t^T\phi_{j_1}(s_1)\psi_1(s_1)\left(
\sum\limits_{j_3=0}^{p}
\left(\int\limits_{s_1}^T\phi_{j_3}(s)\psi(s)ds_1\right)^2
-
\int\limits_{s_1}^T\psi^2(s)ds\right)ds_1\right)^2
$$

\vspace{2mm}
\begin{equation}
\label{otit9500}
=\frac{1}{4}\sum_{j_1=0}^p\left(
\int\limits_t^T\phi_{j_1}(s_1)\psi_1(s_1)
\sum\limits_{j_3=p+1}^{\infty}
\left(\int\limits_{s_1}^T\phi_{j_3}(s)\psi(s)ds\right)^2
ds_1\right)^2.
\end{equation}

\vspace{6mm}

In order to get (\ref{otit9500}) we used the Parseval equality
in the form

$$
\sum_{j_3=0}^{\infty}\left(\int\limits_{s_1}^T\phi_{j_3}(s)
\psi(s)ds\right)^2=
\int\limits_t^T K^2(s,s_1)ds=\int\limits_{s_1}^T\psi^2(s)ds,
$$

\vspace{3mm}
\noindent
where
$$
K(s,s_1)=\psi(s){\bf 1}_{\{s> s_1\}},\ \ \ s, s_1\in [t, T].
$$

\vspace{5mm}

Further consideration is similar to the proof of Case 3
from Theorem 4. 
Finally, from (\ref{otit9500})
we get

$$
{\sf M}\left\{\left(\sum_{j_1=0}^p\left(\sum_{j_3=0}^p
C_{j_3j_3j_1}-\frac{1}{2} C_{j_1}^{*}\right)
\zeta_{j_1}^{(i_1)}\right)^2\right\}<
$$

\vspace{1mm}
$$
< K\frac{p}{p^2}\left(
\int\limits_{-1}^1\frac{dy}{(1-y^2)^{3/4}}
+\int\limits_{-1}^1\frac{dy}{(1-y^2)^{1/4}}\right)^2\le
$$

\vspace{2mm}
$$
\le \frac{K_1}{p}\ \ \ \to 0\ \ \ \hbox{if}\ p\ \to \infty,
$$

\vspace{5mm}
\noindent 
where $K, K_1$ are constants. Case 3 is proved.

Let us consider Case 4. We will prove w.~p.~1 the following
relation

\vspace{1mm}
$$
\hbox{\vtop{\offinterlineskip\halign{
\hfil#\hfil\cr
{\rm l.i.m.}\cr
$\stackrel{}{{}_{p\to \infty}}$\cr
}} }
\sum\limits_{j_1=0}^{p}\sum\limits_{j_3=0}^{p}
C_{j_1 j_3 j_1}\zeta_{j_3}^{(i_2)}=0\ \ \ 
(\psi_1(s), \psi_2(s), \psi_3(s)\equiv \psi(s)).
$$

\vspace{5mm}

In Case 4 we obtain w.~p.~1

$$
\hbox{\vtop{\offinterlineskip\halign{
\hfil#\hfil\cr
{\rm l.i.m.}\cr
$\stackrel{}{{}_{p\to \infty}}$\cr
}} }
\sum\limits_{j_1, j_3=0}^{p}
C_{j_1 j_3 j_1}\zeta_{j_3}^{(i_2)}=
$$

\vspace{2mm}
$$
=\hbox{\vtop{\offinterlineskip\halign{
\hfil#\hfil\cr
{\rm l.i.m.}\cr
$\stackrel{}{{}_{p\to \infty}}$\cr
}} }
\sum\limits_{j_1, j_3=0}^{p}
\left(\frac{1}{2}C_{j_1}^2 C_{j_3}-C_{j_1 j_1 j_3}-C_{j_3 j_1 j_1}
\right)\zeta_{j_3}^{(i_2)}=
$$

\vspace{2mm}
$$
=
\hbox{\vtop{\offinterlineskip\halign{
\hfil#\hfil\cr
{\rm l.i.m.}\cr
$\stackrel{}{{}_{p\to \infty}}$\cr
}} }
\frac{1}{2}\sum\limits_{j_1=0}^{p}
C_{j_1}^2\sum\limits_{j_3=0}^{p}C_{j_3}\zeta_{j_3}^{(i_2)}-
\hbox{\vtop{\offinterlineskip\halign{
\hfil#\hfil\cr
{\rm l.i.m.}\cr
$\stackrel{}{{}_{p\to \infty}}$\cr
}} }
\sum\limits_{j_1, j_3=0}^{p}
C_{j_1 j_1 j_3}\zeta_{j_3}^{(i_2)}-
$$

\vspace{2mm}
$$
-
\hbox{\vtop{\offinterlineskip\halign{
\hfil#\hfil\cr
{\rm l.i.m.}\cr
$\stackrel{}{{}_{p\to \infty}}$\cr
}} }
\sum\limits_{j_1, j_3=0}^{p}
C_{j_3 j_1 j_1}\zeta_{j_3}^{(i_2)}
=
$$

\vspace{2mm}
$$
=\frac{1}{2}\sum\limits_{j_1=0}^{\infty}
C_{j_1}^2\int\limits_t^T\psi(s)d{\bf f}_s^{(i_2)}
-\frac{1}{2}\int\limits_t^T\psi^2(s)
\int\limits_t^s\psi(s_1)d{\bf f}_{s_1}^{(i_2)}ds-
$$

\vspace{2mm}
$$
-\frac{1}{2}\int\limits_t^T\psi(s)
\int\limits_t^s\psi^2(s_1)ds_1d{\bf f}_{s}^{(i_2)}=
\frac{1}{2}\int\limits_t^T\psi^2(s)ds
\int\limits_t^T\psi(s)d{\bf f}_s^{(i_2)}-
$$

\vspace{2mm}
$$
-\frac{1}{2}\int\limits_t^T\psi(s_1)
\int\limits_{s_1}^T\psi^2(s)dsd{\bf f}_{s_1}^{(i_2)}-
\frac{1}{2}\int\limits_t^T\psi(s_1)
\int\limits_t^{s_1}\psi^2(s)dsd{\bf f}_{s_1}^{(i_2)}=
$$

\vspace{2mm}
$$
=\frac{1}{2}\int\limits_t^T\psi^2(s)ds
\int\limits_t^T\psi(s)d{\bf f}_s^{(i_2)}-
\frac{1}{2}\int\limits_t^T\psi(s_1)
\int\limits_t^T\psi^2(s)dsd{\bf f}_{s_1}^{(i_2)}=0,
$$

\vspace{6mm}
\noindent
where we used the Parseval equality in the form

\vspace{1mm}
$$
\sum\limits_{j_1=0}^{\infty}
C_{j}^2=
\sum\limits_{j=0}^{\infty}
\left(\int\limits_t^T\psi(s)\phi_j(s)ds\right)^2
=\int\limits_t^T\psi^2(s)ds.
$$

\vspace{5mm}

Case 4 and Theorem 6 are proved.

Let us consider the trigonometric version of Theorem 6.

\vspace{2mm}

{\bf Theorem 7}\ \cite{16}, \cite{19}, \cite{20}-\cite{20xx2}.  
{\it Suppose that
$\{\phi_j(x)\}_{j=0}^{\infty}$ is a complete orthonormal
system of trigonometric functions
in the space $L_2([t, T])$
and $\psi_1(\tau),$ $\psi_2(\tau),$ $\psi_3(\tau)$ are continuously
differentiable functions at the interval $[t, T]$.
Then, for 
the iterated Stratonovich stochastic integral of third multiplicity

$$
J^{*}[\psi^{(3)}]_{T,t}={\int\limits_t^{*}}^T\psi_3(t_3)
{\int\limits_t^{*}}^{t_3}\psi_2(t_2)
{\int\limits_t^{*}}^{t_2}\psi_1(t_1)
d{\bf f}_{t_1}^{(i_1)}
d{\bf f}_{t_2}^{(i_2)}d{\bf f}_{t_3}^{(i_3)}\ \ \ (i_1, i_2, i_3=1,\ldots,m)
$$

\vspace{3mm}
\noindent
the following 
expansion 

$$
J^{*}[\psi^{(3)}]_{T,t}=
\hbox{\vtop{\offinterlineskip\halign{
\hfil#\hfil\cr
{\rm l.i.m.}\cr
$\stackrel{}{{}_{p\to \infty}}$\cr
}} }
\sum\limits_{j_1, j_2, j_3=0}^{p}
C_{j_3 j_2 j_1}\zeta_{j_1}^{(i_1)}\zeta_{j_2}^{(i_2)}\zeta_{j_3}^{(i_3)}
$$

\vspace{4mm}
\noindent
converging in the mean-square sense 
is valid for each of the following cases

\vspace{2mm}
\noindent
{\rm 1}.\ $i_1\ne i_2,\ i_2\ne i_3,\ i_1\ne i_3,$\\
{\rm 2}.\ $i_1=i_2\ne i_3$ and
$\psi_1(\tau)\equiv\psi_2(\tau)$,\\
{\rm 3}.\ $i_1\ne i_2=i_3$ and
$\psi_2(\tau)\equiv\psi_3(\tau)$,\\
{\rm 4}.\ $i_1, i_2, i_3=1,\ldots,m$
and $\psi_1(\tau)\equiv\psi_2(\tau)\equiv\psi_3(\tau)$,\

\vspace{2mm}
\noindent
where
$$
C_{j_3 j_2 j_1}=\int\limits_t^T\psi_3(s)\phi_{j_3}(s)
\int\limits_t^s\psi_2(s_1)\phi_{j_2}(s_1)
\int\limits_t^{s_1}\psi_1(s_2)\phi_{j_1}(s_2)ds_2ds_1ds.
$$
}

\vspace{4mm}

{\bf Proof.}\ We have

\vspace{1mm}
$$
\int\limits_t^{s}\phi_{j_1}(\theta)\psi(\theta)d\theta=
\frac{\sqrt{2}}{\sqrt{T-t}}
\int\limits_t^{s}
\begin{cases}
\psi(\theta)\ {\rm sin}\left(\left(2\pi j_1(\theta-t)\right)/(T-t)\right)
d\theta\cr\cr
\psi(\theta)\ {\rm cos}\left(\left(2\pi j_1(\theta-t)\right)/(T-t)\right)
d\theta
\end{cases}
=
$$

\vspace{5mm}
$$
=\sqrt{\frac{T-t}{2}}\frac{1}{\pi j_1}\left(\
\begin{cases}
-\psi(s)\ {\rm cos}\left(\left(2\pi j_1(s-t)\right)/(T-t)\right)+\psi(t)
\cr\cr
\psi(s)\ {\rm sin}\left(\left(2\pi j_1(s-t)\right)/(T-t)\right)
\end{cases}+\right.
$$

\vspace{5mm}

$$
\left.+\int\limits_t^{s}
\begin{cases}
\psi'(\theta)\ {\rm cos}\left(\left(2\pi j_1(\theta-t)\right)/(T-t)\right)
d\theta\cr\cr
-\psi'(\theta)\ {\rm sin}\left(\left(2\pi j_1(\theta-t)\right)/(T-t)\right)
d\theta
\end{cases}\right),
$$

\vspace{7mm}
\noindent
where $j_1\ne 0$ and 
$\{\phi_j(x)\}_{j=0}^{\infty}$ is a complete orthonormal
system of trigonometric functions
in $L_2([t, T])$.

Then
\begin{equation}
\label{2017x11}
\left|\int\limits_t^{s}\phi_{j_1}(\theta)\psi(\theta)d\theta\right|\le
\frac{K}{j_1}\ \ \ (j_1\ne 0).
\end{equation}

\vspace{3mm}

Analogously, we get

\vspace{-2mm}
\begin{equation}
\label{2017x12}
\left|\int\limits_s^{T}\phi_{j_1}(\theta)\psi(\theta)d\theta\right|\le
\frac{K}{j_1}\ \ \ (j_1\ne 0).
\end{equation}

\vspace{6mm}

Using (\ref{otit5000}), (\ref{otit9500})--(\ref{2017x12}), 
we obtain
 
$$
{\sf M}\left\{\left(\sum_{j_3=0}^p\left(\sum_{j_1=0}^p
C_{j_3j_1j_1}-\frac{1}{2}\tilde C_{j_3}\right)
\zeta_{j_3}^{(i_3)}\right)^2\right\}\le \frac{K_1}{p}\ \to 0\ \ \
\hbox{if}\ \ \  p\ \to \infty,
$$

\vspace{1mm}
$$
{\sf M}\left\{\left(\sum_{j_1=0}^p\left(\sum_{j_3=0}^p
C_{j_3j_3j_1}-\frac{1}{2} C_{j_1}^{*}\right)
\zeta_{j_1}^{(i_1)}\right)^2\right\}
\le \frac{K_1}{p}\ \to 0\ \ \ \hbox{if}\ \ \ p\ \to \infty,
$$

\vspace{5mm}
\noindent
where constant $K_1$ does not depend on $p.$

The consideration of Case 4 is similar to the case of Legendre 
polynomials (see Theorem 6). Theorem 7 is proved.

Note that the analogues of Theorems 6 and 7 have been proved in
\cite{21} without the restrictions 1--4 (see the formulations of
Theorems 6 and 7). However, in \cite{21} the additional
smoothness assumptions were used.

\vspace{5mm}

\section{Expansions of Iterated Stratonovich Stochastic Integrals
of Multiplicities 3 to 6. Some Recent Results}

\vspace{5mm}

Recently, a new approach to the expansion and mean-square 
approximation of iterated Stratonovich stochastic integrals has been obtained
\cite{20xx} (Sect.~2.10--2.19), \cite{21} (Sect.~13--21), 
\cite{25} (Sect.~5--12), \cite{arxiv-11} (Sect.~7--14), \cite{new-art-1xxy},
\cite{new-art-1xxys}.
Let us formulate five theorems that were proved using this approach.

\vspace{2mm}

{\bf Theorem 8}\ \cite{20xx}, \cite{21}, \cite{25}, \cite{arxiv-11}, \cite{new-art-1xxy}.\
{\it Suppose 
that $\{\phi_j(x)\}_{j=0}^{\infty}$ is a complete orthonormal system of 
Legendre polynomials or trigonometric functions in the space $L_2([t, T]).$
Furthermore, let $\psi_1(\tau), \psi_2(\tau),$ $\psi_3(\tau)$ are continuously dif\-ferentiable 
nonrandom functions on $[t, T].$ 
Then, for the 
iterated Stra\-to\-no\-vich stochastic integral of third multiplicity

$$
J^{*}[\psi^{(3)}]_{T,t}={\int\limits_t^{*}}^T\psi_3(t_3)
{\int\limits_t^{*}}^{t_3}\psi_2(t_2)
{\int\limits_t^{*}}^{t_2}\psi_1(t_1)
d{\bf w}_{t_1}^{(i_1)}
d{\bf w}_{t_2}^{(i_2)}d{\bf w}_{t_3}^{(i_3)}\ \ \ (i_1,i_2,i_3=0,1,\ldots,m)
$$

\vspace{4mm}
\noindent
the following 
relations

\vspace{-1mm}
\begin{equation}
\label{fin1}
J^{*}[\psi^{(3)}]_{T,t}
=\hbox{\vtop{\offinterlineskip\halign{
\hfil#\hfil\cr
{\rm l.i.m.}\cr
$\stackrel{}{{}_{p\to \infty}}$\cr
}} }
\sum\limits_{j_1, j_2, j_3=0}^{p}
C_{j_3 j_2 j_1}\zeta_{j_1}^{(i_1)}\zeta_{j_2}^{(i_2)}\zeta_{j_3}^{(i_3)},
\end{equation}

\vspace{3mm}
\begin{equation}
\label{fin2}
{\sf M}\left\{\left(
J^{*}[\psi^{(3)}]_{T,t}-
\sum\limits_{j_1, j_2, j_3=0}^{p}
C_{j_3 j_2 j_1}\zeta_{j_1}^{(i_1)}\zeta_{j_2}^{(i_2)}\zeta_{j_3}^{(i_3)}\right)^2\right\}
\le \frac{C}{p}
\end{equation}

\vspace{5mm}
\noindent
are fulfilled, where $i_1, i_2, i_3=0,1,\ldots,m$ in {\rm (\ref{fin1})} and 
$i_1, i_2, i_3=1,\ldots,m$ in {\rm (\ref{fin2})},
constant $C$ is independent of $p,$

$$
C_{j_3 j_2 j_1}=\int\limits_t^T\psi_3(t_3)\phi_{j_3}(t_3)
\int\limits_t^{t_3}\psi_2(t_2)\phi_{j_2}(t_2)
\int\limits_t^{t_2}\psi_1(t_1)\phi_{j_1}(t_1)dt_1dt_2dt_3
$$

\vspace{4mm}
\noindent
and
$$
\zeta_{j}^{(i)}=
\int\limits_t^T \phi_{j}(\tau) d{\bf f}_{\tau}^{(i)}
$$ 

\vspace{2mm}
\noindent
are independent standard Gaussian random variables for various 
$i$ or $j$ {\rm (}in the case when $i\ne 0${\rm );} 
another notations are the same as in Theorems~{\rm 1, 2}.}

\vspace{2mm}

{\bf Theorem 9}\ \cite{20xx}, \cite{21}, \cite{25}, \cite{arxiv-11}, \cite{new-art-1xxy}.\
{\it Let
$\{\phi_j(x)\}_{j=0}^{\infty}$ be a complete orthonormal system of 
Legendre polynomials or trigonometric functions in the space $L_2([t, T]).$
Furthermore, let $\psi_1(\tau), \ldots,$ $\psi_4(\tau)$ be continuously dif\-ferentiable 
nonrandom functions on $[t, T].$ 
Then, for the 
iterated Stra\-to\-no\-vich stochastic integral of fourth multiplicity

\begin{equation}
\label{fin0}
J^{*}[\psi^{(4)}]_{T,t}={\int\limits_t^{*}}^T\psi_4(t_4)
{\int\limits_t^{*}}^{t_4}\psi_3(t_3)
{\int\limits_t^{*}}^{t_3}\psi_2(t_2)
{\int\limits_t^{*}}^{t_2}\psi_1(t_1)
d{\bf w}_{t_1}^{(i_1)}
d{\bf w}_{t_2}^{(i_2)}d{\bf w}_{t_3}^{(i_3)}d{\bf w}_{t_4}^{(i_4)}
\end{equation}

\vspace{4mm}
\noindent
the following 
relations

\begin{equation}
\label{fin3}
J^{*}[\psi^{(4)}]_{T,t}
=\hbox{\vtop{\offinterlineskip\halign{
\hfil#\hfil\cr
{\rm l.i.m.}\cr
$\stackrel{}{{}_{p\to \infty}}$\cr
}} }
\sum\limits_{j_1, j_2, j_3,j_4=0}^{p}
C_{j_4j_3 j_2 j_1}\zeta_{j_1}^{(i_1)}\zeta_{j_2}^{(i_2)}\zeta_{j_3}^{(i_3)}\zeta_{j_4}^{(i_4)},
\end{equation}

\vspace{3mm}

\begin{equation}
\label{fin4}
{\sf M}\left\{\left(
J^{*}[\psi^{(4)}]_{T,t}-
\sum\limits_{j_1, j_2, j_3, j_4=0}^{p}
C_{j_4 j_3 j_2 j_1}\zeta_{j_1}^{(i_1)}\zeta_{j_2}^{(i_2)}\zeta_{j_3}^{(i_3)}
\zeta_{j_4}^{(i_4)}
\right)^2\right\}
\le \frac{C}{p^{1-\varepsilon}}
\end{equation}

\vspace{5mm}
\noindent
are fulfilled, where $i_1, \ldots , i_4=0,1,\ldots,m$ in {\rm (\ref{fin0}),} {\rm (\ref{fin3})} 
and $i_1, \ldots, i_4=1,\ldots,m$ in {\rm (\ref{fin4}),}
constant $C$ does not depend on $p,$
$\varepsilon$ is an arbitrary
small positive real number 
for the case of complete orthonormal system of 
Legendre polynomials in the space $L_2([t, T])$
and $\varepsilon=0$ for the case of
complete orthonormal system of 
trigonometric functions in the space $L_2([t, T]),$

$$
C_{j_4 j_3 j_2 j_1}=
$$

$$
=
\int\limits_t^T\psi_4(t_4)\phi_{j_4}(t_4)
\int\limits_t^{t_4}\psi_3(t_3)\phi_{j_3}(t_3)
\int\limits_t^{t_3}\psi_2(t_2)\phi_{j_2}(t_2)
\int\limits_t^{t_2}\psi_1(t_1)\phi_{j_1}(t_1)dt_1dt_2dt_3dt_4;
$$

\vspace{4mm}
\noindent
another notations are the same as in Theorem~{\rm 8}.}

\vspace{2mm}

{\bf Theorem 10}\ \cite{20xx}, \cite{21}, \cite{25}, \cite{arxiv-11}, \cite{new-art-1xxy}.\
{\it Assume 
that $\{\phi_j(x)\}_{j=0}^{\infty}$ is a complete orthonormal system of 
Legendre polynomials or trigonometric functions in the space $L_2([t, T])$
and $\psi_1(\tau), \ldots,$ $\psi_5(\tau)$ are continuously dif\-ferentiable 
nonrandom functions on $[t, T].$ 
Then, for the 
iterated Stra\-to\-no\-vich stochastic integral of fifth multiplicity

\begin{equation}
\label{fin7}
J^{*}[\psi^{(5)}]_{T,t}={\int\limits_t^{*}}^T\psi_5(t_5)
\ldots
{\int\limits_t^{*}}^{t_2}\psi_1(t_1)
d{\bf w}_{t_1}^{(i_1)}
\ldots d{\bf w}_{t_5}^{(i_5)}
\end{equation}

\vspace{4mm}
\noindent
the following 
relations

\begin{equation}
\label{fin8}
J^{*}[\psi^{(5)}]_{T,t}
=\hbox{\vtop{\offinterlineskip\halign{
\hfil#\hfil\cr
{\rm l.i.m.}\cr
$\stackrel{}{{}_{p\to \infty}}$\cr
}} }
\sum\limits_{j_1,\ldots,j_5=0}^{p}
C_{j_5 \ldots j_1}\zeta_{j_1}^{(i_1)}\ldots \zeta_{j_5}^{(i_5)},
\end{equation}

\vspace{3mm}

\begin{equation}
\label{fin9}
{\sf M}\left\{\left(
J^{*}[\psi^{(5)}]_{T,t}-
\sum\limits_{j_1, \ldots, j_5=0}^{p}
C_{j_5 \ldots j_1}\zeta_{j_1}^{(i_1)}\ldots
\zeta_{j_5}^{(i_5)}
\right)^2\right\}
\le \frac{C}{p^{1-\varepsilon}}
\end{equation}

\vspace{5mm}
\noindent
are fulfilled, where $i_1, \ldots , i_5=0,1,\ldots,m$ in {\rm (\ref{fin7}),} {\rm (\ref{fin8})} 
and $i_1, \ldots, i_5=1,\ldots,m$ in {\rm (\ref{fin9}),}
constant $C$ is independent of $p,$
$\varepsilon$ is an arbitrary
small positive real number 
for the case of complete orthonormal system of 
Legendre polynomials in the space $L_2([t, T])$
and $\varepsilon=0$ for the case of
complete orthonormal system of 
trigonometric functions in the space $L_2([t, T]),$

$$
C_{j_5 \ldots j_1}=
\int\limits_t^T\psi_5(t_5)\phi_{j_5}(t_5)\ldots
\int\limits_t^{t_2}\psi_1(t_1)\phi_{j_1}(t_1)dt_1\ldots dt_5;
$$

\vspace{3mm}
\noindent
another notations are the same as in Theorems~{\rm 8, 9}.}

\vspace{2mm}

{\bf Theorem 11}\ \cite{20xx}, \cite{21}, \cite{25}, \cite{arxiv-11}, \cite{new-art-1xxy}.\
{\it Suppose that 
$\{\phi_j(x)\}_{j=0}^{\infty}$ is a complete orthonormal system of 
Legendre polynomials or trigonometric functions in the space $L_2([t, T]).$
Then, for the 
iterated Stratonovich stochastic integral of sixth multiplicity

\begin{equation}
\label{after10001qu1}
J_{T,t}^{*(i_1\ldots i_6)}={\int\limits_t^{*}}^T
\ldots
{\int\limits_t^{*}}^{t_2}
d{\bf w}_{t_1}^{(i_1)}
\ldots d{\bf w}_{t_6}^{(i_6)}
\end{equation}

\vspace{3mm}
\noindent
the following 
expansion 

\vspace{-1mm}
$$
J_{T,t}^{*(i_1\ldots i_6)}
=\hbox{\vtop{\offinterlineskip\halign{
\hfil#\hfil\cr
{\rm l.i.m.}\cr
$\stackrel{}{{}_{p\to \infty}}$\cr
}} }
\sum\limits_{j_1, \ldots, j_6=0}^{p}
C_{j_6 \ldots j_1}\zeta_{j_1}^{(i_1)}\ldots
\zeta_{j_6}^{(i_6)}
$$

\vspace{4mm}
\noindent
that converges in the mean-square sense is valid, where
$i_1, \ldots, i_6=0, 1,\ldots,m,$

$$
C_{j_6 \ldots j_1}=
\int\limits_t^T\phi_{j_6}(t_6)\ldots
\int\limits_t^{t_2}\phi_{j_1}(t_1)dt_1\ldots dt_6;
$$

\vspace{3mm}
\noindent
another notations are the same as in Theorems~{\rm 8--10}.}

\vspace{2mm}

{\bf Theorem~12}\ \ \cite{20xx}, \cite{21}, \cite{25}, \cite{arxiv-11}.\
{\it Suppose that 
$\{\phi_j(x)\}_{j=0}^{\infty}$ is a complete orthonormal system of 
Legendre polynomials or trigonometric functions in the space $L_2([t, T]).$
Furthermore, let $\psi_1(\tau), \psi_2(\tau), \psi_3(\tau)$ are continuously dif\-ferentiable 
nonrandom functions on $[t, T]$.
Then, for the 
iterated Stra\-to\-no\-vich stochastic integral of third multiplicity

$$
J^{*}[\psi^{(3)}]_{T,t}^{(i_1 i_2 i_3)}={\int\limits_t^{*}}^T\psi_3(t_3)
{\int\limits_t^{*}}^{t_3}\psi_2(t_2)
{\int\limits_t^{*}}^{t_2}\psi_1(t_1)
d{\bf w}_{t_1}^{(i_1)}
d{\bf w}_{t_2}^{(i_2)}d{\bf w}_{t_3}^{(i_3)}
$$

\vspace{3mm}
\noindent
the following 
expansion 

\vspace{-1mm}
$$
J^{*}[\psi^{(3)}]_{T,t}^{(i_1 i_2 i_3)}
=\hbox{\vtop{\offinterlineskip\halign{
\hfil#\hfil\cr
{\rm l.i.m.}\cr
$\stackrel{}{{}_{p_1,p_2,p_3\to \infty}}$\cr
}} }
\sum\limits_{j_1=0}^{p_1}\sum\limits_{j_2=0}^{p_2}\sum\limits_{j_3=0}^{p_3}
C_{j_3 j_2 j_1}\zeta_{j_1}^{(i_1)}\zeta_{j_2}^{(i_2)}\zeta_{j_3}^{(i_3)}
$$

\vspace{4mm}
\noindent
that converges in the mean-square sense is valid, where
$i_1, i_2, i_3=0, 1,\ldots,m,$

$$
C_{j_3 j_2 j_1}=\int\limits_t^T\psi_3(t_3)\phi_{j_3}(t_3)
\int\limits_t^{t_3}\psi_2(t_2)\phi_{j_2}(t_2)
\int\limits_t^{t_2}\psi_1(t_1)\phi_{j_1}(t_1)dt_1dt_2dt_3
$$

\vspace{3mm}
\noindent
another notations are the same as in Theorems~{\rm 8--11}.}

\vspace{2mm}

Obviously, Theorem~12 generalizes the main results of this article
for iterated Stratonovich stochastic integrals of third multiplicity.

The results of Theorems 8--11 were developed in \cite{20xx} (Chapter~2), 
\cite{21}, \cite{25}, \cite{arxiv-11}.
In particular, analogues of Theorem~11 for iterated Stratonovich stochastic
integrals of multiplicities 7 and 8 were obtained in \cite{20xx} (Sect.~2.36, 2.37).
In addition, the variants of Thorems 8--11
were obtained
for the case when $\{\phi_j(x)\}_{j=0}^{\infty}$ is an arbitrary complete orthonormal system
of functions in $L_2([t, T])$ \cite{20xx} (Sect.~2.1.4, 2.23, 2.24, 2.31--2.34),
\cite{21}, \cite{25}, \cite{arxiv-11}.

\vspace{5mm}

\section{Theorems 1--12 from Point
of View of the Wong--Zakai Approximation}

\vspace{5mm}

The iterated Ito stochastic integrals and solutions
of Ito SDEs are complex and important functionals
from the independent components ${\bf f}_{s}^{(i)},$
$i=1,\ldots,m$ of the multidimensional
Wiener process ${\bf f}_{s},$ $s\in[0, T].$
Let ${\bf f}_{s}^{(i)p},$ $p\in\mathbb{N}$ 
be some approximation of
${\bf f}_{s}^{(i)},$
$i=1,\ldots,m$.
Suppose that 
${\bf f}_{s}^{(i)p}$
converges to
${\bf f}_{s}^{(i)},$
$i=1,\ldots,m$ if $p\to\infty$ in some sense and has
differentiable sample trajectories.

A natural question arises: if we replace 
${\bf f}_{s}^{(i)}$
by ${\bf f}_{s}^{(i)p},$
$i=1,\ldots,m$ in the functionals
mentioned above, will the resulting
functionals converge to the original
functionals from the components 
${\bf f}_{s}^{(i)},$
$i=1,\ldots,m$ of the multidimentional
Wiener process ${\bf f}_{s}$?
The answere to this question is negative 
in the general case. However, 
in the pioneering works of Wong E. and Zakai M. \cite{W-Z-1},
\cite{W-Z-2},
it was shown that under the special conditions and 
for some types of approximations 
of the Wiener process the answere is affirmative
with one peculiarity: the convergence takes place 
to the iterated Stratonovich stochastic integrals
and solutions of Stratonovich SDEs and not to iterated 
Ito stochastic integrals and solutions
of Ito SDEs.
The piecewise 
linear approximation 
as well as the regularization by convolution 
\cite{W-Z-1}-\cite{Watanabe} relate the 
mentioned types of approximations
of the Wiener process. The above approximation 
of stochastic integrals and solutions of SDEs 
is often called the Wong--Zakai approximation.

Let ${\bf w}_{\tau},$ $\tau\in[0, T]$ is a random vector with 
an $m+1$ components: ${\bf w}_{\tau}^{(i)}={\bf f}_{\tau}^{(i)}$ 
for $i=1,\ldots,m$ and 
${\bf w}_{\tau}^{(0)}=\tau,$\ 
${\bf f}_{\tau}^{(i)}$ $(i=1,\ldots,m)$
are independent standard Wiener processes.

It is well known that the following representation 
takes place \cite{Lipt}, \cite{7e}

\begin{equation}
\label{um1x}
{\bf w}_{\tau}^{(i)}-{\bf w}_{t}^{(i)}=
\sum_{j=0}^{\infty}\int\limits_t^{\tau}
\phi_j(s)ds\ \zeta_j^{(i)},\ \ \ \zeta_j^{(i)}=
\int\limits_t^T \phi_j(s)d{\bf w}_s^{(i)},
\end{equation}

\vspace{4mm}
\noindent
where $\tau\in[t, T],$ $t\ge 0,$
$\{\phi_j(x)\}_{j=0}^{\infty}$ is an arbitrary complete 
orthonormal system of functions in the space $L_2([t, T]),$ and
$\zeta_j^{(i)}$ are independent standard Gaussian 
random variables for various $i$ or $j.$
Moreover, the series (\ref{um1x}) converges for any $\tau\in [t, T]$
in the mean-square sense.

Let ${\bf w}_{\tau}^{(i)p}-{\bf w}_{t}^{(i)p}$ be 
the mean-square approximation of the process
${\bf w}_{\tau}^{(i)}-{\bf w}_{t}^{(i)},$
which has the following form

\vspace{-3mm}
\begin{equation}
\label{um1xx}
{\bf w}_{\tau}^{(i)p}-{\bf w}_{t}^{(i)p}=
\sum_{j=0}^{p}\int\limits_t^{\tau}
\phi_j(s)ds\ \zeta_j^{(i)}.
\end{equation}

\vspace{3mm}

From (\ref{um1xx}) we obtain

\vspace{-4mm}
\begin{equation}
\label{um1xxx}
d{\bf w}_{\tau}^{(i)p}=
\sum_{j=0}^{p}
\phi_j(\tau)\zeta_j^{(i)} d\tau.
\end{equation}

\vspace{4mm}

Consider the following iterated Riemann--Stieltjes
integral

\begin{equation}
\label{um1xxxx}
\int\limits_t^T
\psi_k(t_k)\ldots \int\limits_t^{t_2}\psi_1(t_1)
d{\bf w}_{t_1}^{(i_1)p_1}\ldots d{\bf w}_{t_k}^{(i_k)p_k},
\end{equation}

\vspace{4mm}
\noindent
where $p_1,\ldots,p_k\in\mathbb{N},$\ \ $i_1,\ldots,i_k=0,1,\ldots,m,$ 

\begin{equation}
\label{um1xxx1}
d{\bf w}_{\tau}^{(i)p}=
\left\{\begin{matrix}
d{\bf f}_{\tau}^{(i)p}\ &\hbox{\rm for}\ \ \ i=1,\ldots,m\cr\cr\cr
d\tau^p\ &\hbox{\rm for}\ \ \ i=0
\end{matrix}
,\right.
\end{equation}

\vspace{4mm}
\noindent
and $d{\bf f}_{\tau}^{(i)p},$ $d\tau^p$ are defined by the relation (\ref{um1xxx}).

Let us substitute (\ref{um1xxx}) into (\ref{um1xxxx})

\begin{equation}
\label{um1xxxx1}
\int\limits_t^T
\psi_k(t_k)\ldots \int\limits_t^{t_2}\psi_1(t_1)
d{\bf w}_{t_1}^{(i_1)p_1}\ldots d{\bf w}_{t_k}^{(i_k)p_k}=
\sum\limits_{j_1=0}^{p_1}\ldots \sum\limits_{j_k=0}^{p_k}
C_{j_k \ldots j_1}\prod\limits_{l=1}^k \zeta_{j_l}^{(i_l)},
\end{equation}

\vspace{4mm}
\noindent
where 
$$
\zeta_j^{(i)}=\int\limits_t^T \phi_j(s)d{\bf w}_s^{(i)}
$$ 

\vspace{2mm}
\noindent
are independent standard Gaussian random variables for various 
$i$ or $j$ (in the case when $i\ne 0$),
${\bf w}_{s}^{(i)}={\bf f}_{s}^{(i)}$ for
$i=1,\ldots,m$ and 
${\bf w}_{s}^{(0)}=s,$

$$
C_{j_k \ldots j_1}=\int\limits_t^T\psi_k(t_k)\phi_{j_k}(t_k)\ldots
\int\limits_t^{t_2}
\psi_1(t_1)\phi_{j_1}(t_1)
dt_1\ldots dt_k
$$

\vspace{4mm}
\noindent
is the Fourier coefficient.

To best of our knowledge \cite{W-Z-1}-\cite{Watanabe}
the approximations of the Wiener process
in the Wong--Zakai approximation must satisfy fairly strong
restrictions
\cite{Watanabe}
(see Definition 7.1, pp.~480--481).
Moreover, approximations of the Wiener process that are
similar to (\ref{um1xx})
were not considered in \cite{W-Z-1}, \cite{W-Z-2}
(also see \cite{Watanabe}, Theorems 7.1, 7.2).
Therefore, the proof of analogs of Theorems 7.1 and 7.2 \cite{Watanabe}
for approximations of the Wiener 
process based on its series expansion (\ref{um1x})
should be carried out separately.
Thus, the mean-square convergence of the right-hand side
of (\ref{um1xxxx1}) to the iterated Stratonovich stochastic integral 
(\ref{str})
does not follow from the results of the papers
\cite{W-Z-1}, \cite{W-Z-2} (also see \cite{Watanabe},
Theorems 7.1, 7.2).

From the other hand, Theorems 1, 2 and Theorems 3--12 from this 
paper (also see Chapters 1 and 2 from \cite{20xx})  can be considered as the proof of the
Wong--Zakai approximation for the iterated 
Stratonovich stochastic integrals (\ref{str}) of multiplicities 1 to 6
based on the approximation (\ref{um1xx}) of the Wiener process.
At that, the iterated Riemann--Stieltjes integrals (\ref{um1xxxx}) converge
(according to Theorems 1--12 from this article and Chapters 1, 2 from \cite{20xx})
to the appropriate iterated Stratonovich 
stochastic integrals (\ref{str}). Recall that
$\{\phi_j(x)\}_{j=0}^{\infty}$ (see (\ref{um1x}), (\ref{um1xx}), and
Theorems 3--12)
is a complete 
orthonormal system of Legendre polynomials or 
trigonometric functions 
in the space $L_2([t, T])$.

To illustrate the above reasoning, 
consider two examples for the case $k=2,$
$\psi_1(s),$ $\psi_2(s)\equiv 1;$ $i_1, i_2=1,\ldots,m.$

The first example relates to the piecewise linear approximation
of the multidimensional Wiener process (these approximations 
were considered in \cite{W-Z-1}-\cite{Watanabe}).

Let ${\bf b}_{\Delta}^{(i)}(t),$ $t\in[0, T]$ be the piecewise
linear approximation of the $i$th component ${\bf f}_t^{(i)}$
of the multidimensional standard Wiener process ${\bf f}_t,$
$t\in [0, T]$ with independent components
${\bf f}_t^{(i)},$ $i=1,\ldots,m,$ i.e.

$$
{\bf b}_{\Delta}^{(i)}(t)={\bf f}_{k\Delta}^{(i)}+
\frac{t-k\Delta}{\Delta}\Delta{\bf f}_{k\Delta}^{(i)},
$$

\vspace{3mm}
\noindent
where 

\vspace{-2mm}
$$
\Delta{\bf f}_{k\Delta}^{(i)}={\bf f}_{(k+1)\Delta}^{(i)}-
{\bf f}_{k\Delta}^{(i)},\ \ \
t\in[k\Delta, (k+1)\Delta),\ \ \ k=0, 1,\ldots, N-1.
$$

\vspace{4mm}

Note that w.~p.~1

\vspace{-1mm}
\begin{equation}
\label{pridum}
\frac{d{\bf b}_{\Delta}^{(i)}}{dt}(t)=
\frac{\Delta{\bf f}_{k\Delta}^{(i)}}{\Delta},\ \ \
t\in[k\Delta, (k+1)\Delta),\ \ \ k=0, 1,\ldots, N-1.
\end{equation}

\vspace{4mm}

Consider the following iterated Riemann--Stieltjes
integral

\vspace{1mm}
$$
\int\limits_0^T
\int\limits_0^{s}
d{\bf b}_{\Delta}^{(i_1)}(\tau)d{\bf b}_{\Delta}^{(i_2)}(s),\ \ \ 
i_1,i_2=1,\ldots,m.
$$

\vspace{4mm}

Using (\ref{pridum}) and additive property of Riemann--Stieltjes integrals, 
we can write w.~p.~1

\vspace{2mm}
$$
\int\limits_0^T
\int\limits_0^{s}
d{\bf b}_{\Delta}^{(i_1)}(\tau)d{\bf b}_{\Delta}^{(i_2)}(s)=
\int\limits_0^T
\int\limits_0^{s}
\frac{d{\bf b}_{\Delta}^{(i_1)}}{d\tau}(\tau)d\tau
\frac{d {\bf b}_{\Delta}^{(i_2)}}{d s}(s)
ds =
$$

\vspace{3mm}
$$
=
\sum\limits_{l=0}^{N-1}\int\limits_{l\Delta}^{(l+1)\Delta}
\left(
\sum\limits_{q=0}^{l-1}\int\limits_{q\Delta}^{(q+1)\Delta}
\frac{\Delta{\bf f}_{q\Delta}^{(i_1)}}{\Delta}d\tau+
\int\limits_{l\Delta}^{s}
\frac{\Delta{\bf f}_{l\Delta}^{(i_1)}}{\Delta}d\tau\right)
\frac{\Delta{\bf f}_{l\Delta}^{(i_2)}}{\Delta}ds=
$$

\vspace{3mm}
$$
=\sum\limits_{l=0}^{N-1}\sum\limits_{q=0}^{l-1}
\Delta{\bf f}_{q\Delta}^{(i_1)}
\Delta{\bf f}_{l\Delta}^{(i_2)}+
\frac{1}{\Delta^2}\sum\limits_{l=0}^{N-1}
\Delta{\bf f}_{l\Delta}^{(i_1)}
\Delta{\bf f}_{l\Delta}^{(i_2)}
\int\limits_{l\Delta}^{(l+1)\Delta}
\int\limits_{l\Delta}^{s}d\tau ds=
$$

\vspace{3mm}
\begin{equation}
\label{oh-ty}
=\sum\limits_{l=0}^{N-1}\sum\limits_{q=0}^{l-1}
\Delta{\bf f}_{q\Delta}^{(i_1)}
\Delta{\bf f}_{l\Delta}^{(i_2)}+
\frac{1}{2}\sum\limits_{l=0}^{N-1}
\Delta{\bf f}_{l\Delta}^{(i_1)}
\Delta{\bf f}_{l\Delta}^{(i_2)}.
\end{equation}

\vspace{6mm}

Using (\ref{oh-ty}) and standard relations between 
Ito and Stratonovich stochastic integrals, it 
is not difficult to show 
that

\vspace{1mm}
$$
\hbox{\vtop{\offinterlineskip\halign{
\hfil#\hfil\cr
{\rm l.i.m.}\cr
$\stackrel{}{{}_{N\to \infty}}$\cr
}} }
\int\limits_0^T
\int\limits_0^{s}
d{\bf b}_{\Delta}^{(i_1)}(\tau)d{\bf b}_{\Delta}^{(i_2)}(s)=
\int\limits_0^T
\int\limits_0^{s}
d{\bf f}_{\tau}^{(i_1)}d{\bf f}_{s}^{(i_2)}+
\frac{1}{2}{\bf 1}_{\{i_1=i_2\}}\int\limits_0^T ds=
$$

\vspace{3mm}
\begin{equation}
\label{uh-111}
=
\int\limits_0^{*T}
\int\limits_0^{*s}
d{\bf f}_{\tau}^{(i_1)}d{\bf f}_{s}^{(i_2)},
\end{equation}

\vspace{5mm}
\noindent
where $\Delta\to 0$ if $N\to\infty$ ($N\Delta=T$).

Obviously, (\ref{uh-111}) agrees with Theorem 7.1 (see \cite{Watanabe},
p.~486).

The next example relates to the approximation
of the Wiener process based on its series expansion
(\ref{um1x}) for $t=0$, where
$\{\phi_j(x)\}_{j=0}^{\infty}$ 
is a complete 
orthonormal system of Legendre polynomials or 
trigonometric functions 
in the space $L_2([0, T])$.

Consider the following iterated Riemann--Stieltjes
integral

\vspace{-1mm}
\begin{equation}
\label{abcd1}
\int\limits_0^T
\int\limits_0^{s}
d{\bf f}_{\tau}^{(i_1)p}d{\bf f}_{s}^{(i_2)p},\ \ \ 
i_1,i_2=1,\ldots,m,
\end{equation}

\vspace{3mm}
\noindent
where $d{\bf f}_{\tau}^{(i)p}$ is defined by the
relation
(\ref{um1xxx}).

Let us substitute (\ref{um1xxx}) into (\ref{abcd1}) 

\vspace{-1mm}
\begin{equation}
\label{set18}
\int\limits_0^T
\int\limits_0^{s}
d{\bf f}_{\tau}^{(i_1)p}d{\bf f}_{s}^{(i_2)p}=
\sum\limits_{j_1,j_2=0}^p
C_{j_2 j_1} \zeta_{j_1}^{(i_1)}\zeta_{j_2}^{(i_2)},
\end{equation}

\vspace{3mm}
\noindent
where 
$$
C_{j_2 j_1}=
\int\limits_0^T \phi_{j_2}(s)\int\limits_0^s
\phi_{j_1}(\tau)d\tau ds
$$

\vspace{3mm}
\noindent
is the Fourier coefficient; another notations 
are the same as in (\ref{um1xxxx1}).

As we noted above, approximations of the Wiener process that are
similar to (\ref{um1xx})
were not considered in \cite{W-Z-1}, \cite{W-Z-2}
(also see Theorems 7.1, 7.2 in \cite{Watanabe}).
Furthermore, the extension of the results of Theorems 7.1 and 7.2
\cite{Watanabe} to the case under consideration is
not obvious.

On the other hand, we can apply the theory built in Chapters 1 and 2
of the monographs \cite{20xx}-\cite{20xx2}. More precisely, 
using 
Theorem 2.2 \cite{20xx},
we obtain from (\ref{set18}) the desired result

\vspace{-1mm}
$$
\hbox{\vtop{\offinterlineskip\halign{
\hfil#\hfil\cr
{\rm l.i.m.}\cr
$\stackrel{}{{}_{p\to \infty}}$\cr
}} }
\int\limits_0^T
\int\limits_0^{s}
d{\bf f}_{\tau}^{(i_1)p}d{\bf f}_{s}^{(i_2)p}
=
\hbox{\vtop{\offinterlineskip\halign{
\hfil#\hfil\cr
{\rm l.i.m.}\cr
$\stackrel{}{{}_{p\to \infty}}$\cr
}} }
\sum\limits_{j_1,j_2=0}^p
C_{j_2 j_1} \zeta_{j_1}^{(i_1)}\zeta_{j_2}^{(i_2)}=
$$

\vspace{2mm}
\begin{equation}
\label{umen-bl}
=
\int\limits_0^{*T}
\int\limits_0^{*s}
d{\bf f}_{\tau}^{(i_1)}d{\bf f}_{s}^{(i_2)}.
\end{equation}

\vspace{5mm}

From the other hand, by Theorems 1, 2
(see (\ref{a2})) for the case
$k=2$ we obtain from (\ref{set18}) the following relation

\vspace{-2mm}
$$
\hbox{\vtop{\offinterlineskip\halign{
\hfil#\hfil\cr
{\rm l.i.m.}\cr
$\stackrel{}{{}_{p\to \infty}}$\cr
}} }
\int\limits_0^T
\int\limits_0^{s}
d{\bf f}_{\tau}^{(i_1)p}d{\bf f}_{s}^{(i_2)p}=
\hbox{\vtop{\offinterlineskip\halign{
\hfil#\hfil\cr
{\rm l.i.m.}\cr
$\stackrel{}{{}_{p\to \infty}}$\cr
}} }
\sum\limits_{j_1,j_2=0}^p
C_{j_2 j_1} \zeta_{j_1}^{(i_1)}\zeta_{j_2}^{(i_2)}=
$$

\vspace{2mm}
$$
=
\hbox{\vtop{\offinterlineskip\halign{
\hfil#\hfil\cr
{\rm l.i.m.}\cr
$\stackrel{}{{}_{p\to \infty}}$\cr
}} }
\sum\limits_{j_1,j_2=0}^p
C_{j_2 j_1} \biggl(\zeta_{j_1}^{(i_1)}\zeta_{j_2}^{(i_2)}-
{\bf 1}_{\{i_1=i_2\}}{\bf 1}_{\{j_1=j_2\}}\biggr)+
{\bf 1}_{\{i_1=i_2\}}\sum\limits_{j_1=0}^{\infty}
C_{j_1 j_1}=
$$

\vspace{2mm}
\begin{equation}
\label{umen-blx}
=
\int\limits_0^T
\int\limits_0^{s}
d{\bf f}_{\tau}^{(i_1)}d{\bf f}_{s}^{(i_2)}+
{\bf 1}_{\{i_1=i_2\}}\sum\limits_{j_1=0}^{\infty}
C_{j_1 j_1}.
\end{equation}

\vspace{5mm}

Since
$$
\sum\limits_{j_1=0}^{\infty}
C_{j_1 j_1}=\frac{1}{2}\sum\limits_{j_1=0}^{\infty}
\left(\int\limits_0^T \phi_j(\tau)d\tau\right)^2
=\frac{1}{2}
\left(\int\limits_0^T \phi_0(\tau)d\tau\right)^2=\frac{1}{2}
\int\limits_0^T ds,
$$

\vspace{5mm}
\noindent
then using standard relations
between Ito and Stratonovich stochastic integrals
and (\ref{umen-blx}) we obtain (\ref{umen-bl}).

\end{document}